\documentclass[english]{article}
\usepackage[fontsize=11pt]{scrextend}
\usepackage[utf8]{inputenc}
\usepackage{graphicx}
\usepackage{amsmath}
\usepackage{empheq}
\usepackage{amsfonts}
\usepackage{enumitem}
\usepackage{tikz}
\usetikzlibrary{arrows.meta}   
\usepackage{amssymb}
\usepackage{comment}
\usepackage[english]{babel}
\usepackage{hyperref} 
\usepackage[top=2cm, bottom=3cm, left=2.5cm , right=2.5cm]{geometry}
\usepackage{mathrsfs}
\usepackage{dsfont}
\usepackage{color}
\usepackage{flushend}
\usepackage{mathtools}
\usepackage{fancyhdr}
\usepackage[thmmarks,amsmath]{ntheorem}
\usepackage{stmaryrd}
\usepackage[nottoc]{tocbibind}
\usepackage[toc,page]{appendix}
\usepackage{appendix}
\usepackage{caption}
\usepackage{xcolor}
\usepackage{float}
\usepackage{textcomp}
\usepackage{bm}  
\usepackage{lipsum}
\usepackage{titlesec}
\usepackage{braket} 
\usepackage{esint} 
\usepackage{todonotes}
\usepackage{centernot}

\titleformat{\paragraph}[block]
{\normalsize\bfseries}{}{0pt}{}

\titlespacing{\paragraph}{0pt}{*1.5}{*0.5}

\DeclareMathOperator{\re}{Re}
\DeclareMathOperator{\im}{Im}
\newcommand{\Ec}[1]{\mathbb{E}_N\left[#1\right]}

\renewcommand{\textbf}[1]{\begingroup\bfseries\mathversion{bold}#1\endgroup}

\newcommand{\N}{\mathbb{N}}

\newcommand{\Q}{\mathbb{Q}}
\newcommand{\E}{\mathbb{E}}
\newcommand{\R}{\mathbb{R}}
\newcommand{\C}{\mathbb{C}}

\renewcommand{\P}{\mathbb{P}}

\newcommand{\mcal}[1]{\mathcal{#1}}
\newcommand{\mc}[1]{\mathcal{#1}}
\newcommand{\mscr}[1]{\mathscr{#1}}
\newcommand{\mbb}[1]{\mathbb{#1}}

\newcommand{\ind}[1]{\mathbf{1}_{#1}}
\newcommand{\diff}{\mathop{}\mathopen{}\mathrm{d}}

\newcommand{\msf}[1]{\mathsf{#1}}
\newcommand{\mbf}[1]{\mathbf{#1}}
\newcommand{\mfrak}[1]{\mathfrak{#1}}
\newcommand{\mrm}[1]{\mathrm{#1}}

\newcommand{\ii}{\mathrm{i}}

\newcommand{\tend}[1]{\underset{#1}{\longrightarrow}}

\theoremstyle{plain}
\newtheorem{Thm}{Theorem}[section]
\newtheorem{Prop}[Thm]{Proposition}
\newtheorem{Lem}[Thm]{Lemma}

\newcommand{\vertiii}[1]{{\left\vert\kern-0.25ex\left\vert\kern-0.25ex\left\vert #1 
		\right\vert\kern-0.25ex\right\vert\kern-0.25ex\right\vert}}

{
	\theorembodyfont{\normalsize}
	\theoremsymbol{\ensuremath{\square}}
	\newtheorem*{Pro}{Proof}
	
}
\newtheorem{Def}[Thm]{Definition}
\newtheorem{Cor}[Thm]{Corollary}
{	
	\theorembodyfont{}
	\newtheorem{Rem}[Thm]{Remark}
}

{	
	\theorembodyfont{\normalsize}
	
}
\newtheorem{assu}[Thm]{Assumptions}
\newcommand{\defi }{\coloneqq}

\begin{document}
	\title{Mesoscopic transition for $\beta$-ensembles at intermediary temperature}
	\author{Charlie Dworaczek Guera\footnote{KTH Royal Institute of Technology, Department of Mathematics, 11428, Stockholm, Sweden. \newline
			\textit{email:} chadg@kth.se}\and Gaultier Lambert\footnote{KTH Royal Institute of Technology, Department of Mathematics, 11428, Stockholm, Sweden. \newline
			\textit{email:} glambert@kth.se}\and Luke Peilen\footnote{College of the Holy Cross, Department of Mathematics and Computer Science, Worcester, MA, United States. \newline
			\textit{email:} lpeilen@holycross.edu}}
	\date{}
	\maketitle
	
		\begin{abstract}
		This paper establishes a mesoscopic central limit theorem for linear statistics of $\beta$-ensembles or log-gas, as the dimension $N\to\infty$, in the temperature regime $1/N\ll\beta(N)\le 1$. For simplicity, we assume that the potential is one-cut regular and analytic. In this regime, the size of the fluctuations depends on $\beta(N)$ and the mesoscopic scale. We show that there is a transition at a critical $\eta\asymp 1/ N\beta(N)$ between a \textit{Random Matrix regime}, where the limit variance is given by the $\msf{H}^{1/2}$-norm and a \textit{Poisson regime} where the limit variance is give by the $L^2$-norm. We also describe the critical regime. 
		
		The proof of the CLT relies on \emph{optimal local laws} at intermediate temperatures and Stein's method for $\beta$-ensembles. In particular, in this regime, it is necessary to construct new correction terms to the classical equilibrium measure to obtain a suitable re-centring of linear statistics and describe their fluctuations.  We also obtain a free energy expansion.
	\end{abstract}

	\tableofcontents
	
	\section{Introduction}
	\subsection{Setting of the problem}
		We are interested in the following probability distribution on $\R^N$:
	\begin{equation}\label{def:mesureparticules}
		\diff\P_N(\bm{\lambda})\defi p_N (\bm{\lambda})\diff\lambda_1\dots \diff\lambda_N\hspace{0,5cm}\text{ with }\hspace{0,5cm}p_N (\bm{\lambda})\defi \dfrac{1}{\mcal{Z}_N[V]}\prod_{i<j}^{N}\left|\lambda_i-\lambda_j\right|^{\beta }.\prod_{i=1}^{N}e^{-\frac{N\beta}{2} V(\lambda_i)} .
	\end{equation}
	where $\beta=\beta_N>0$ and $V:\R\to\R$ is growing sufficiently rapidly so that the partition function
	\begin{equation}\label{def:partfunction}
		\mcal{Z}_N[V]\defi\int_{\R^N}\prod_{i<j}^{N}\left |\lambda_i-\lambda_j\right|^{\beta }.\prod_{i=1}^{N}e^{-\frac{N\beta}{2} V(\lambda_i)} \diff\lambda_i <+\infty.
	\end{equation}
	
	This probability measure is called a $\beta$-ensemble, or \textit{one dimensional log-gas}, and it may be viewed as a statistical physics ensemble of $N$ confined particles interacting through a $2$-body repulsive logarithmic interaction at the inverse temperature $\beta$.
	For $V(x)=x^2$, this corresponds to the joint-distribution of the eigenvalues of the \textit{Gaussian $\beta$-Ensemble} (G$\beta$E) which is a central model in random matrix theory. 
	In addition, for $\beta=2$, such measures are determinantal, and the correlation is expressed in terms of orthogonal polynomials. 
	In general, the regime $\beta>0$ fixed has been extensively studied in the literature.

	In this article, we are interested in the regime where the temperature $\beta = \beta_N$ depends on the dimension such that $1/N\ll\beta\ll1$ as $N\to\infty$, that is referred to as the \textit{intermediary temperature regime}.
	In the sequel, we denote by $\alpha \defi \beta N$ and we study the empirical measure: 
	$$\bm{\mu}_N\defi \dfrac{1}{N}\sum_{i=1}^N\delta_{\lambda_i}$$
	
	It is well known that for a wide class of potentials $V$, called the \textit{one-cut} class, which includes convex smooth potentials growing fast enough at infinity, that
	converges (in a large deviation sense \cite{arous1997large,anderson2010introduction,Garcia}) towards a deterministic measure $\mu_{\mrm{eq}}$, called the \textit{equilibrium measure}, which is supported on a single interval that we can always assume to be $[-1,1]$ (replacing $V(x)$ by $V(ax+b)$ for finely tuned $a,b$). This measure arises as the unique minimizer of the following energy functional:
	\begin{equation}\label{eq:energy functional}
		\mcal{E}(\mu)\defi\int_\R V(x)\diff\mu(x)-\iint_{\R^2}\log|x-y|\diff\mu(x)\diff\mu(y).
	\end{equation}
	As a consequence of the large deviation principle, the limit $\bm{\mu}_N(f)\rightarrow\mu_{\mrm{eq}}(f)$ when $N\rightarrow\infty$ follows. A natural question that follows is the study of the fluctuations (when appropriately rescaled) of $[\bm{\mu}_N-\mu_{\mrm{eq}}](f)$. In the fixed-$\beta$ case, it is well known since the seminal work \cite{johansson} that the following Gaussian convergence holds:
	$$N(\bm{\mu}_N-\mu_{\mrm{eq}})(f)\overset{\mrm{law}}{\tend{N\rightarrow\infty}}\left(\dfrac{1}{\beta}-\dfrac{1}{2}\right)\nu_1(f) +\dfrac{1}{\sqrt{\beta}}\;\mcal{N}\left(0,\sigma^{2}(f)\right) $$ for $f$ smooth enough where $\nu_1$ is a certain linear map and $\sigma^{2}$ a positive quadratic form. For $V(x)=x^{2}$, it is known that $\nu_1=\tfrac{1}{2}\left(\delta_{-1}+\delta_{1}\right)-\ind{|x|\leq1}\pi^{-1}(1-x^{2})^{-1/2}\diff x$.
	
	At intermediary temperature, the situation is slightly different. Assuming that there exists a smallest integer $p$ such that $\alpha^{p}\gg \sqrt{N\alpha}$ (\textit{i.e.} when $\beta\gg N^{-1+\kappa}$ for $\kappa>0$), we prove here that there exist linear functionals $(\nu_k)_{k\in\llbracket1,p\rrbracket}$ such that for all $f$ sufficiently smooth:
	$$\sqrt{N\alpha}\left(\bm{\mu}_N-\mu_{\mrm{eq}}-\sum_{k=1}^{p}\dfrac{\nu_k}{\alpha^{k}}\right)(f)\overset{\mrm{law}}{\tend{N\rightarrow\infty}}\mcal{N}\left(0,\sigma^{2}(f)\right).$$
	
\begin{figure}[h]
	\begin{center}
		\begin{tikzpicture}[scale=6]
			
			\draw[line width=1pt,-{Latex[length=3mm,width=2mm]}] 
			(0,0) -- (1.075,0) ;
			\draw[line width=1pt,-{Latex[length=3mm,width=2mm]}] 
			(0,0) -- (0,1.075);
			
			\draw[thick] (0,0) rectangle (1,1);
			
			\draw[dashed] (0,1) -- (1,0);
			
			\node[text width=3cm] at (0.45,0.35) {$\alpha\eta\gg1$};
			\node[text width=3cm] at (0.352,0.25) {Random Matrix};
			\node at (0.35,0.15) {$\mathsf{H}^{1/2}$};
			\node[text width=3cm] at (0.9,0.85) {$\alpha\eta\ll1$};
			\node[text width=3cm] at (0.9,0.75) {Poisson};
			\node at (0.76,0.65) {$L^2$};
			\node at (0.5,-0.07) {$\beta^{-1}$};
			\node at (1,-0.07) {$N$};
			\node at (-0.09,0.5) {$\eta^{-1}$};
			\node at (-0.09,1) {$N$};
			\node at (-0.0,-0.05) {$1$};
			\draw[dashed] (0,1) -- (1,0);
			
			\node[rotate=-45] at (0.53,0.53) {$\alpha\eta\tau\asymp1$\quad transition $\Sigma_\tau$};
			\node[rotate=90] at (1.05,0.5) {high temperature regime};
		\end{tikzpicture}
		
		\caption{At intermediary temperature, the mesoscopic fluctuations at scale $\eta\gg1/\alpha$, \textit{i.e.} in the Random matrix regime, are given by the usual $\msf{H}^{1/2}$-norm defined in \eqref{def: H12 norm}. Zooming in, for $\eta\ll1/\alpha$, the fluctuations in the Poisson regime, are given after a different rescaling by a $L^{2}$-norm. Finally, let $\tau>0$, an interpolating quadratic form $\Sigma_\tau$ arise as the variance at the transition, namely when $\alpha\eta\tau\asymp1$.}
		\label{fig:transition}
	\end{center}
\end{figure}
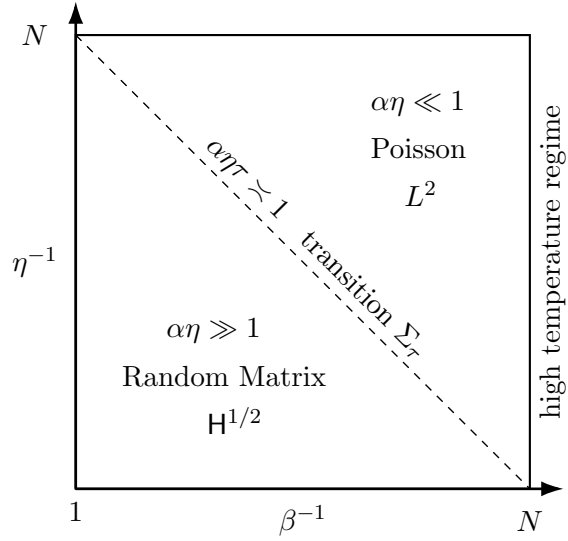

	For a point $E\in(-1,1)$, in the bulk, it is also interesting to look at mesoscopic fluctuations at scale $\eta\gg1/N$ \textit{i.e.} to understand the limit of $[\bm{\mu}_N-\mu_{\mrm{eq}}](f_\eta)$ for $f_\eta=f\left(\eta^{-1}(\cdot-E)\right)$ and $f$ a function decaying sufficiently fast at infinity. For the constant $\beta$ case, it is known \cite{bekerman2018mesoscopic,peilen2024local} that:
	$$N(\bm{\mu}_N-\mu_{\mrm{eq}})(f_\eta)\overset{\mrm{law}}{\tend{N\rightarrow\infty}}\dfrac{1}{\sqrt{\beta}}\;\mcal{N}\left(0,\|f\|_{\msf{H}^{1/2}}^{2}\right) $$
	for $\eta\gg1/N$, $f$ sufficiently smooth and decaying fast enough at infinity. The limiting variance is given by: 
	\begin{equation}\label{def: H12 norm}
		\|f\|_{\msf{H}^{1/2}}^{2}\defi\dfrac{1}{2\pi^{2}}\iint_{\R^2}\left(\dfrac{f(x)-f(y)}{x-y}\right)^{2}\diff x\diff y. 
	\end{equation}
	At intermediary temperature, assuming $\beta\geq N^{-1+\kappa}$ for some $\kappa>0$, the situation is also slightly different. We show that a transition occurs which is the main result of this article.
	This transition is summarized in Figure \ref{fig:transition}.

	\subsection{Main result and assumptions}
	In what follows, $x\gg y$ (resp. $x\asymp y$,) means for $x,y>0$ that $x/y\rightarrow+\infty$ (resp. there exists $c,c'>0$ such that $cx\leq y\leq c'x$). We recall that it is well-known that there exists
	$C_{\mrm{eq}}\in\R$ such that:
	\begin{equation}\label{eq:effpotentialV}
		V_{\mrm{eff}}(x)\defi V(x)-2\int_\R\log|x-y|\diff\mu _{\mrm{eq}}(y)-C_\mrm{eq}\begin{cases}
			\geq 0\hspace{1cm}& \text{for all }x\in\R
			\\= 0&\mu_{\mrm{eq}}-\text{almost-surely}
		\end{cases}
	\end{equation} We work under the following assumptions throughout the article:
		\begin{assu}
	\label{assumptions}
	The potential $V$ satisfies:
	\begin{enumerate}[label=(\roman*)]
		\item\label{assumption1} $V\in\mcal{C}^{0}(\R)$, is real-analytic in a complex neighborhood $\mc{U}$ of $[-1,1]$ and $\liminf\limits_{|x|\rightarrow\infty}\dfrac{V(x)}{\log|x|}=+\infty$ 
		\item \label{assumption2}(One-cut) The support of $\mu _{\mrm{eq}}$ is the connected interval $[-1,1]$ and is given by:
		$$\dfrac{\diff\mu _{\mrm{eq}}(x)}{\diff x}=\dfrac{S(x)\sqrt{1-x^2}}{\pi}$$
		where $S(x)>0$ on $[-1,1]$. Recall that $S$ can be extended analytically on $\mcal{U}$ by
		\begin{equation}\label{eq: S}
			S(z)\defi\dfrac{1}{2\pi}\int_{-1}^1\dfrac{V'(t)-V'(z)}{t-z}\dfrac{\diff t}{\sqrt{1-t^2}}.
		\end{equation}
		\item \label{assumption3}(Off-criticality) The function $x\mapsto V(x)-2\displaystyle\int_{-1}^1\log|x-y|\diff\mu_{\mrm{eq}}(y)$ achieves its minimum value on $[-1,1]$ only.
	\end{enumerate}
	Furthermore, we assume that $\alpha=N\beta\gg( \log N)^{k}$ for all $k>0$.
\end{assu}

The assumption \textit{ \ref{assumption1}} is not so restrictive since by using Lemma \ref{lem:replacement}, we can assume, up to a fixed cost of size $e^{-c\alpha}$ for $c>0$, that $V$ is a $\mcal{C}_{\mrm{loc}}^{\infty}(\R)$ potential equal to $x^{2}+C$ for a constant $C>0$ large enough (see Lemma \ref{lem:replacement}) outside of a compact $K$. This type of errors being irrelevant in this paper, \textbf{we thus assume from now on that $V\in\mcal{C}_{\mrm{loc}}^{\infty}(\R)$ equal to $x^{2}+C$ outside of a compact $K$.} The assumption on the growth of $V$ is necessary in order to show that the model is well-defined. The analyticity assumption is not necessary in our setup and could be removed by replacing $V(z)$ by the pseudo-analytic extension of degree $k$ of $V$ for $k$ large enough, as in \cite{dworaczekguera2025clt}. This last assumption simplifies the notation.

The assumption \textit{\ref{assumption2}} is necessary to use the inversion of the master operator $\Xi$ defined in Definition \ref{def:ope D}. This property no longer holds in the multi-cut situation and the CLT \ref{thm:main CLT global} is also known to be false in this situation. The invertibility of $\Xi$ is a major ingredient of the proof.


The assumption \textit{\ref{assumption3}} is necessary to ensure that the global CLT is true. Indeed, if $V_{\mrm{eff}}$ were to vanish at some point $x_0>1$, taking a test function $f$ supported around $x_0$, $\bm{\mu}_N(f)$ would not be asymptotically Gaussian  by \cite{fan2015convergence}.

The assumption on $\alpha$ comes from the fact that for $\alpha=\log N$, the large deviation principle for the largest particle $\lambda_{\max}$, Proposition \ref{appprop: LDP lambdamax} no longer holds. Indeed, with high probability there will be outliers (\textit{i.e.} particles lying outside of the asymptotic support) which weakens an approach using truncation arguments.

We recall that for all $1\leq p<\infty$, the definition of the \textit{$p$-Wasserstein} distance between two random variables $X,Y\in L^p$ not necessarily defined on the same probability space, is given by:
$$\mathbf{W}_p(X,Y)\defi\inf\left\{\|X'-Y'\|_p,X'\overset{\mrm{law}}{=}X,Y'\overset{\mrm{law}}{=}Y\right\}.$$
In Proposition \ref{prop:asympcorrelators}, we prove that for some $\beta\gg N^{-1+\kappa}$ for $\kappa>0$, there exist linear functionals $(\nu_i)_{i\geq1}$ such that for all $f\in\mcal{C}^{r}_{\mrm{loc}}(\R)$ such that $x\in\R\mapsto f(x)e^{-2\kappa V(x)}\in\mcal{C}^{r}(\R)$ for some $\kappa>0$ and $r$ large enough, there exists $j>0$ such that:
$$\Ec{\int_\R f(x)\diff (\mu_N-\mu_{\mrm{eq}})(x)}=\sum_{i=1}^{j}\dfrac{\nu_i(f)}{\alpha^{i}}+O\left(\dfrac{\|fe^{-\kappa V}\|_{\mcal{C}^{r}(\R)}}{N}\right).$$
Therefore, for all $k\geq0$, we define the following linear maps: $$\mu_k\coloneqq\mu_{\mrm{eq}}+\sum_{i=1}^{k}\frac{\nu_i}{\alpha^{i}},\quad\quad\quad L_N^{(k)}\coloneqq\bm{\mu}_N-\mu_k.$$

We start with the following global CLT with explicit rates of convergence. 
\begin{Thm}[Global CLT]\label{thm:main CLT global}Let $k\geq0$, $\beta\gg N^{-1+\kappa}$ for some $\kappa>0$. Let $f\in\mcal{C}^{r}(\R)$ for $r$ big enough (depending only on $\kappa$), then for all $q\geq1$:
	$$	\mathbf{W}_q\left(\sqrt{N\alpha}L_N^{(k)}(f),\mcal{N}\left(0,\sigma^{2}(f) \right)\right) \leq C(f) \Big(\dfrac{q}{\sqrt{N\alpha}}+\dfrac{\sqrt{q}}{\alpha}+\dfrac{\sqrt{N\alpha}}{\alpha^{k+1}}\Big)$$
	where
	\begin{equation}\label{eq: norm H 1/2 global}
		\sigma^2(f)\defi\dfrac{1}{2\pi^2}\iint_{[-1,1]^2}\left(\dfrac{f(x)-f(y)}{x-y}\right)^2 \cdot \dfrac{1-xy}{\sqrt{1-x^2}\sqrt{1-y^2}}\diff x\diff y.
	\end{equation}
\end{Thm}


The following mesoscopic CLT in the bulk is the main result of the present paper. This exhibits a transition that occurs at the mesoscopic scale $\eta\asymp1/\alpha$. 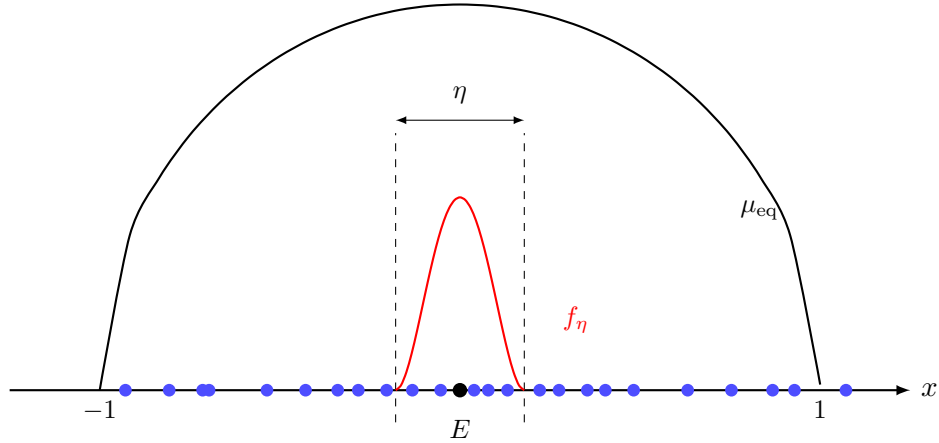
\begin{figure}[h]
	\centering
	\begin{tikzpicture}[>=latex,scale=1.7]
		
		\draw[->,thick] (-3.5,0) -- (3.5,0) node[right] {$x$};
		
		\node[below] at (-2.8,0) {\small $-1$};
		\node[below] at (2.8,0) {\small $1$};
		
		\draw[thick,smooth,domain=-2.8:2.8,variable=\x]
		plot ({\x},{3*(1-(\x/2.8)^2)^(0.5)});
		
		\node[right] at (2.1,1.4) {\small $\mu_{\mathrm{eq}}$};
		
		\foreach \x in {-2.6,-2.26,-2.0,-1.95,-1.5,-1.2,-0.95,-0.79,-0.57,-0.37,
			-0.15,0.11,0.22,0.37,0.62,0.77,0.99,1.13,1.35,1.77,2.11,2.43,2.6}
		{
			\fill[blue!70] (\x,0) circle (1.4pt);
		}
		
		\fill[blue!70] (3.0,0) circle (1.4pt);
		
		\fill (0,0) circle (1.6pt);
		\node[below] at (0,-0.15) {\small $E$};
		
		\def\phi{1}
		
		\draw[dashed] ({-1/2},-0.25) -- ({-1/2},2.0);
		\draw[dashed] ({1/2},-0.25) -- ({1/2},2.0);
		
		\draw[<->] ({-1/2},2.1) -- ({1/2},2.1);
		\node at (0,2.3) {\small $\eta$};
		
		\draw[red,thick,domain={-1/2}:{1/2},smooth,variable=\x]
		plot ({\x},{1.5*(1 - (2*\x)^2)^2});
		
		\node[red] at (0.9,0.55) {\small $f_\eta$};
		

	\end{tikzpicture}
	\caption{An illustration of a mesoscopic function $f_\eta$ centered at a point $E$ in the bulk.}
\end{figure}


\begin{Thm}[Mesoscopic CLT in the bulk]\label{thm:mesoc clt bulk}
	Let $N^{-1+\kappa}\ll\beta_N\ll \tfrac{1}{(\log N)^{1+\kappa}}$ for some $\kappa>0$. For some $k\geq1$, large enough, for all $E\in(-1+\varepsilon,1-\varepsilon)$ for some fixed $\varepsilon>0$ and all $q\geq1$, with $f_\eta\defi f\left(\eta^{-1}(\cdot-E)\right)$ and explicit rates $(\mfrak{R}_i)_{i\in\llbracket1,4\rrbracket}$ given in Theorem \ref{thm:meso sec5}, the following holds:
	\begin{enumerate}
		\item \textbf{\textit{(Random Matrix regime)}} if $\tfrac{1}{\alpha}\ll\eta\ll1$ then for all $f\in\mcal{C}_c^{5}(\R)$:
\begin{equation*}
	\mathbf{W}_q\left (\sqrt{N\alpha}L_N^{(k)}(f_\eta), \mcal{N}\left(0,\|f\|_{\msf{H}^{1/2}}^{2}\right)\right ) \lesssim \begin{cases}
		\mfrak{R}_1 & \text{if } \eta \geq \tfrac{N^\delta}{\alpha}, \\
	\mfrak{R}_2&\text{otherwise}.\end{cases}
\end{equation*}Furthermore, for all fixed $q\geq1$, $\mfrak{R}_1$ and $\mfrak{R}_2$ are a $o(1)$.
	\item \textbf{(\textit{Poisson regime})} if $\tfrac{1}{N}\ll\eta\ll\theta\ll\tfrac{1}{\alpha}$ then for all $f\in\mcal{C}_c^{3}(\R)$:
		$$		\mathbf{W}_q\left(\sqrt{\tfrac{N}{\eta \mu_{\mrm{eq}}(E)}}L_N^{(k)}(f_\eta), \mathcal{N}(0, \|f\|_{L^2(\R)}^2)\right) \lesssim \mfrak{R}_3.$$
		Furthermore, for all fixed $q\geq1$, $\mfrak{R}_3=o(1)$.
			\item \textbf{(\textit{Critical regime})} If $\alpha\eta\mu_\mrm{eq}(E)\tau= 1$ for some $\tau>0$ then for all $f\in\mcal{C}_c^{3}(\R)$:
				$$\mbf{W}_q\left (\sqrt{\tfrac{N}{\eta \mu_{\mrm{eq}}(E)}}L_N^{(k)}(f_\eta), \mathcal{N}(0, \Sigma_\tau^2(f))\right )\lesssim\mfrak{R}_4.$$
				Furthermore, for all fixed $q\geq1$, $\mfrak{R}_4=o(1)$.
	\end{enumerate}
	The semi-norm $\|\cdot\|_{\msf{H}^{1/2}}$ is given in \eqref{def: H12 norm} and
	\begin{equation*}\label{eq:Sigmac}
	\Sigma_\tau^2(\phi):=\frac{\tau}{2\pi}\int_{\R} \frac{|\xi|\cdot|\widehat{\phi}(\xi)|^2}{\pi+\tau|\xi|}\diff \xi.
	\end{equation*}
\end{Thm}

The main ingredient of the proof is a local law valid at all scales. We introduce the Stieltjes transform of $\bm{\mu}_N$ and $\mu_p$ for all $p\geq0$ and $z\in\C_+\defi\{z\in\C,\, \im{z}>0\}$:
$$s_N(z)\defi\dfrac{1}{N}\sum_{i=1}^{N}\dfrac{1}{\lambda_i-z},\quad\quad\quad\quad m_p(z)\coloneqq\mu_p\left(\dfrac{1}{\cdot-z}\right).$$ 
	\begin{Thm}[Local law]\label{thm:local law}
		Let $V$  satisfying Assumptions \ref{assumptions}, $\beta_N\gg N^{-1+\kappa}$ for some $\kappa>0$, and $\varepsilon>0$, there exists $\tilde{\eta},C,k>0$ (depending on $V$, $\kappa$ and $\varepsilon$) such that for any $q\geq1$, fixed $\delta>0$ and $z=x+\ii \eta$ with $-1+\varepsilon\leq x\leq1-\varepsilon$, we have for $N$ large enough:
	 $$\|s_N(z)-m_k(z)\|_q\lesssim\ind{C\frac{N^{\delta}}{\alpha}<\eta<\tilde{\eta}}\sqrt{\dfrac{q}{N\alpha\eta^{2}}}+\ind{\eta\leq C\frac{N^{\delta}}{\alpha}}\left(\sqrt{\dfrac{q}{N\eta}}+\dfrac{q}{N\eta}\right) .$$

	\end{Thm}
		
	
%
%
	
	Applying the loop equations method in the spirit of \cite{BoG1}, we can deduce the asymptotic expansion for the free energy $\log \mcal{Z}_N[V]$. Note that the expansion of $\mcal{Z}_N[2x^{2}]$ is known because of Mehta's formula \cite[17.6.7]{mehta2004random} or \cite[Section 7.2.2]{BoG2}, see Remark \ref{rem:asymp gaussian partition}.
\begin{Thm}[Expansion of the free energy]\label{thm:free energy}
	Let $V, \beta_N$ satisfying Assumptions \ref{assumptions}, there exists real numbers $\left(F_{i,j}\right)_{i,j}$ such that for all $K,L\geq1$:
$$\dfrac{2}{N\alpha}\log\dfrac{\mcal{Z}_N[V]}{\mcal{Z}_N[2x^{2}]}=\sum_{j=0}^{L}\sum_{i=0}^{K-j}\dfrac{F_{i,j}}{\alpha^iN^{j}}+o\left(\dfrac{1}{ \alpha^{K+1}}+\dfrac{1}{N^{L+1}}\right).$$
Furthermore, $\lim\limits_{N\rightarrow\infty}\dfrac{2}{N\alpha}\log\mcal{Z}_N[V]=-\mcal{E}(\mu_{\mrm{eq}})$ and $F_{1,0}=2\mrm{Ent}\left[\mu_{\mrm{eq}}\right]-2\log\pi+1$.
\end{Thm}
	
%
%
	\subsection{Connection with the literature}\label{subsec:connection}
\paragraph{The fixed temperature regime}
The global fluctuations of $\beta$-ensembles for $\beta>0$ fixed were obtained for fixed $\beta$ in the seminal paper \cite{johansson}. It was later generalized in the works \cite{Shc2,BoG1,BoG2,bekerman2018clt,LambertLedouxWebb,angst2024sharp,dworaczekguera2025clt} allowing for more general potentials, general domains of integrations, critical potentials, explicit rates of convergence, minimal regularity conditions or singular potentials. The mesoscopic fluctuations were also obtained for $\beta=2$ \cite{lambert1511clt} using the determinantal structure and for general $\beta>0$ \cite{bekerman2018mesoscopic,lambert2021mesoscopic,peilen2024local} using respectively, the loop equations, Johansson's method and the next-order energy. The analogue of Theorem \ref{thm:free energy}, \textit{i.e.} the asymptotics of the partition function were established in \cite{Shc1,BoG1,BoG2,guera2024asymptotics}. Optimal local laws were established in \cite{bourgade2022optimal} generalizing the results of \cite{BEY14,bourgade2012bulk,BouErdYau2014}.
\paragraph{High temperature regime and the thermal equilibrium measure}
The high-temperature regime, which has been extensively studied over the past fifteen years,  corresponds to the case $\beta=\alpha/N$ for $\alpha>0$ fixed. In this regime, the limiting measure is not $\mu_{\mrm{eq}}$ anymore but the so-called \textit{thermal equilibrium measure }$\mu_\alpha$, \textit{i.e.} the unique minimizer of the following functional:
$$\mcal{E}_\alpha(\mu)\defi\int_\R V(x)\diff\mu(x)-\iint_{\R^2}\log|x-y|\diff\mu(x)\diff\mu(y)+\dfrac{2}{\alpha}\int_\R\log\left(\dfrac{\diff\mu(x)}{\diff x}\right )\diff\mu(x).$$
The convergence holds in a large deviation sense \cite{Garcia}. A remarkable feature of $\mu_\alpha$, is that it is supported on the whole real line because of the entropic term in $\mcal{E}_\alpha$ (which is the only difference with $\mcal{E}$). This measure is known to interpolate between the measure $e^{-V(x)\diff x}/Z$ and $\mu_{\mrm{eq}}$ taking $\alpha\rightarrow0$ and $\alpha\rightarrow+\infty$, this is called the \textit{Gauss-Wigner} crossover \cite{AllezBouchaudGuionnet}. The global fluctuations of linear statistics were established to be asymptotically Gaussian of size $\sqrt{N\alpha}\sim\sqrt{N}$ in \cite{nakano2018gaussian,hardy2021clt,dworaczek2024clt}. The local statistics were also established to converge to a Poisson point process in this regime \cite{benaych2015poisson,nakano2018gaussian,nakano2020poisson,lambert2021poisson} and the asymptotic expansion of the partition function was proven in \cite{guera2024asymptotics}. The connection between the high temperature regime and integrable systems, like the Toda chain, was studied in \cite{Spohn1,GMToda,mazzuca2023large,mazzuca2022mean,mazzuca2023clt}.

At intermediary temperature, when $\alpha\gg1$, one can still use $\mu_\alpha$ (thus depending on $N$) to obtain a more accurate description of the $1$-point function $\mbb{E}[\bm{\mu}_N]$ as was argued in \cite{armstrong2022thermal}. In principle, one expects that by using $\mu_\alpha$ for any $\alpha>0$ (possibly diverging in $N$) to recenter $\bm{\mu}_N$, one could obtain the mesoscopic fluctuations as it is supposed to be a more precise recentring than $\mu_k$ for all $k\geq0$. However, the drawback is while the linear maps $\mu_k$ are quite explicit and satisfy nice recurrence relations, $\mu_\alpha$ is non-linear and highly non-explicit. While this measure is well-understood in dimension greater than 2, using the Coulomb-gas structure, see \cite{armstrong2022thermal}, a good understanding of this object is still lacking in one dimension. On the real line, the only example known where $\mu_\alpha$ is explicit is for $V(x)=x^{2}$ where $\mu_\alpha$ is known as the \textit{Askey-Wimp-Kerov} distribution \cite{AllezBouchaudGuionnet}.
\paragraph{Comparison with the Coulomb gas literature}
The probability density in \eqref{def:mesureparticules} can equivalently be described as $p_N (\bm{\lambda})= \dfrac{1}{\mcal{Z}_N[V]}e^{-\beta H(\bm{\lambda})}$, with
\begin{equation*}
H(\bm{\lambda})=\sum_{i<j} \log \frac{1}{|\lambda_i-\lambda_j|}+\frac{N}{2}\sum_{i=1}^N V(\lambda_i).
\end{equation*}
This kind of interaction is a special case of the more general Riesz gas, an interacting particle system in $\R^{\mathsf{d}}$ whose interaction kernel replaces $-\log|\cdot|$ with
\begin{equation*}
\mathsf{g}(\cdot)=\begin{cases}
\frac{1}{\mathsf{s}}|x|^{-\mathsf{s}} & \text{if }0< \mathsf{s}<\mathsf{d} \\
-\log|x| & \text{if }\mathsf{s}=0.
\end{cases}
\end{equation*}
In the special case $\mathsf{d}-2$, these are referred to as Coulomb gases. For $\mathsf{d}=2$ at inverse temperature $\beta=2$, the Coulomb gases correspond to the eigenvalue density of the Ginibre ensemble. A thorough discussion of the statistical mechanics of these systems, as well as their applications to diverse areas of mathematics and physics, can be found in \cite{serfaty26}.

Fluctuations for the $\mathsf{d}=2$ Coulomb gas have been well-studied at constant and intermediate temperature regimes, and the full CLT in what we call the ``random-matrix" regime above has been derived in \cite{serfaty2023gaussian} for fluctuations against the thermal equilibrium measure (see also \cite{RV07}, \cite{leble2018fluctuations} and \cite{BBNY19} for fluctuations at constant $\beta$). In addition to the log-gas results discussed above, CLTs for the Riesz gas in dimensions $\mathsf{d}=1,2$ have been derived at fixed temperature $\beta \simeq 1$ in \cite{boursier1} and \cite{peilenserfaty25}. Results at intermediate temperatures, as well as the corresponding phase transition diagram for general $\mathsf{s}$ and $\mathsf{d}$, are largely open and present interesting directions for future work.

\subsection{Strategy of the proof}
The main result of this paper is  Theorem \ref{thm:mesoc clt bulk}, \textit{i.e.} the mesoscopic CLT in the bulk  which we prove by using Stein's method in the spirit of \cite{LambertLedouxWebb,angst2024sharp}. Schematically, our proof can be summarized as follows:
$$\text{Optimal local law}\hspace{,5cm}\overset{\text{Helffer–Sjöstrand formula}}{\rightsquigarrow}\hspace{,5cm}\text{Stein's method},$$
From the optimal local law, we deduce, thanks to Helffer–Sjöstrand formula, sharp mesoscopic concentration estimates which is enough to apply Stein's method and deduce the mesoscopic CLT.
\paragraph{Mesoscopic concentration estimates}
It is known that a key input to run Stein's method is an optimal concentration estimate for mesoscopic linear statistics $L_N^{(k)}(f_\eta)$ where $f\in\mcal{C}_c^{r}(\R)$, for some $r>0$ large enough and $\eta\rightarrow0$ is the mesoscopic scale. Such optimal concentration estimates are control on $	\E_N\left[|L_N^{(k)}(f_\eta)|^{q}\right]$ in terms of $q$, $N$ and $\eta$ for $f_\eta$ a mesoscopic functions. Our strategy to obtain these estimates rely on an optimal local law in the bulk (see Theorem \ref{thm:local law}).

\paragraph{Local laws}
We now explain how to obtain this local law. The main ingredient of the proof are the so-called \textit{loop equations}, introduced in Section \ref{sec: local laws}, which are relations between different linear statistics obtained by integration by parts. The first one is given by: 
\begin{equation}\label{introeq:1stloopeq}
	\E_N\left[s_N^2(z) + V'(z) s_N(z) +\mu_{\mrm{eq}}(W_z) + \Delta_0 + \left( \dfrac{2}{\alpha}-\dfrac{1}{N}\right) \partial_z s_N(z)\right] 
	=0,
\end{equation}
where
$$\Delta_0(z)=L_N^{(0)}(W_z),\hspace{1cm}W_z(\lambda) \coloneqq\dfrac{V'(\lambda)-V'(z)}{\lambda-z}.$$
It was shown in \cite[Theorem 1.1]{bourgade2022optimal} that, for fixed $\beta>0$, these equations combined with an optimal bound on $\Delta_0$ allow to obtain an optimal local law. In the intermediary temperature regime,  such a strategy allows to obtain a first non-optimal local law in the RM regime. To improve the local law one must consider better recentring than the Stieltjes transform of $\mu_{\mrm{eq}}$. Indeed, a term $\tfrac{1}{\alpha}\partial_zs_N(z)$ which is negligible in the regime $\beta$ constant, plays in this context a crucial role and is the largest errror term appearing in the loop equations. To obtain a finer local law, one must try to recenter it as well.
 A key ingredient of the proof of \cite{bourgade2022optimal} is the fact that the Stieltjes transform of $\mu_{\mrm{eq}}$ is the solution of a quadratic equation, we build on this idea and create a sequence of quadratic polynomials $(P_k)$ (defined in Definition $\ref{def:Pk}$), one of whose solution is approximately equal to the Stieltjes transform of $\mu_k$. The use of the latter solution allows to iteratively improve the local law in the RM regime leading to Theorem \ref{prop:LL1}. For the Poisson regime, noticing that the entropic term ($\partial_zs_N(z)$) becomes bigger than the energy contribution (coming from the quadratic equation), a combination of differential inequalities (such as the Gronwall Lemma) and the loop equations allow to conclude about Theorem \ref{prop:LL2}.
%

	\subsection{Notations and conventions}
	\textbf{Spaces and norms:} We denote by $\mcal{M}_1(\R)$ the set of probability measures on $\R$. The upper half-plane is denoted by $\C_+=\left\{z\in\C,\,\im{z}>0\right\}$.  We denote by $\mcal{C}^{r}_{\mrm{loc}}(X)$ the space of functions defined on $X$ continuously differentiable $0\leq r\leq\infty$ times. We denote by $\mcal{C}^{r}(X)$, the subspace of functions for which the following norm is finite:
	$$\|f\|_{\mcal{C}^{r}(X)}\defi\max_{k\in\llbracket0,r\rrbracket}\sup_{x\in X}|f^{(k)}(x)|.$$
	The space of functions in $\mcal{C}^{k}(\R)$ supported on a compact $K$ is denoted $\mcal{C}_c^{k}(K)$. Let $X$ be a random variable, we define for all $p\geq1$, $\|X\|_p\defi\Ec{|X|^{p}}^{\frac{1}{p}}$. We define the distance $d$ on the space of measures by:
	$$d(\mu,\mu')\defi\sup_{\substack{\|f\|_{\msf{H}^{1/2}}\leq1\\\|f\|_{\mrm{Lip}}\leq1}}\Big|\int_\R f(x)\diff(\mu-\mu')(x)\Big|,$$
	where $\|\cdot\|_{\msf{H}^{1/2}}$ is defined in \eqref{def: H12 norm} and $\|f\|_{\mrm{Lip}}\defi\sup_{x,y\in\R}\Big|\dfrac{f(x)-f(y)}{x-y}\Big|$. We denote by $\mcal{W}^{k,p}(\R)$ the space $\{f\in L^{p}(\R),\forall j\in\llbracket0,k\rrbracket,\,f^{(j)}\in L^{p}(\R) \}$, we denote by $\|f\|_{\mcal{W}^{k,p}(\R)}=\sum_{i=0}^{k}\|f^{(k)}\|_{L^{p}(\R)}$. The Schwartz space is denoted by $\mcal{S}(\R)$.
	
	\textbf{Operators.} Let $\mcal{A}:\mcal{C}_{\mrm{loc}}^{j}(\R)\rightarrow\mcal{C}_{\mrm{loc}}^{j'}(\R^{k})$ for some $k\geq1$, we extend the definition of $\mcal{A}$ on $\mcal{C}_{\mrm{loc}}^{j}(\R^{\ell+1})$  by $\mc{A}[\psi](x_1,\dots,x_{k+\ell})\coloneqq\mc{A}[\psi(\cdot,x_{k+1},\dots,x_{k+\ell})](x_1,\dots,x_k)$ for all $\psi\in\mcal{C}_{\mrm{loc}}^{j}(\R^{\ell+1})$ for some $\ell\geq 1$. The Fourier transform of $\phi\in L^1(\R)\cap L^2(\R)$ is denoted by
	$$\hat{\phi}(x)\coloneqq\int_\R\phi(y)e^{-\ii xy}\diff y.$$
	The Hilbert transform of $\phi\in L^2(\R)$, is defined by density by the following formula:
	$$\mcal{H}[\phi](x)\coloneqq\fint_\R\dfrac{\phi(y)}{y-x}\dfrac{\diff y}{\pi}$$
	where $\fint$ denotes the Cauchy principal value integral.
	
	\textbf{General notations.} In the present article, the constants we used are allowed to change from line to line, can never depend on $N$ or any variable but can depend on fixed parameters such as $V$, $\tilde{\eta}$ etc...  Given a linear map $T$ (for example induced by a measure) on a space of functions, we define $\braket{T,\varphi}=T(\varphi)$. We denote by $z^{*}$ the complex conjugate of $z\in\C$. Let $\mcal{V}$ be a complex open set of $\C$, we denote by $\mcal{V}_+\coloneqq \mcal{V}\cap\mbb{C}_+$. For $x,y>0$, we say that $x\ll y$ (respectively $x\lesssim y$, $x\asymp y$) if $y/x\rightarrow+\infty$ (repectively. there exists a constant $C>0$ depending only on $V$ such that $x\leq C y$, resp. if there are two constants $c_1,c_2>0$ such that $c_1 x\leq y\leq c_2 x$). For a function of $f\in\mcal{C}^{1}_{\mrm{loc}}(\R^{n})$, we denote by $\partial_kf$ the partial derivative of $f$ with respect to the $k$-th variable. We denote by $\bar{\partial}$ the derivative $\frac{1}{2}(\partial_x+\ii\partial_y)$. For two integers $m>n$, we denote by $\llbracket n,m\rrbracket$ the set $\{n,\hdots,m\}$. $\sqrt{z^{2}-1}$ denotes the unique square root with cut on $[-1,1]$ and which behaves at infinity like $z+o(z)$.

		\subsection{Outline of the paper} In Section \ref{sec:bound Delta}, we define the upgraded recentring via Proposition \ref{prop:asympcorrelators}. We also prove the global CLT Theorem \ref{thm:main CLT global} and the asymptotics of the free energy Theorem \ref{thm:free energy}. In Section \ref{sec: local laws}, we prove the local law Theorem \ref{thm:local law}. In Section \ref{sec: concentration}, we use Helffer–Sjöstrand formula to translate the local law into controls on the mesoscopic linear statistics. In Section \ref{sec: proof CLT}, we finally prove the CLT using the concentration bounds. The Appendix \ref{app:apriori} gathers the proof of the \textit{a priori bound} Proposition \ref{a priori bound} necessary to derive Proposition \ref{prop:asympcorrelators}. In Appendix \ref{app LDP}, we prove the replacement lemma which allows us to consider very general potentials and the large deviation principle for the largest/smallest particles. Finally, in Appendix \ref{app: Gaussian}, we consider the very special case of a quadratic potential. We gather numerical simulations as well as the convergence of the upgraded recentring. In Appendix \ref{app:Stein}, we gather inequalities from Stein's method that we use to obtain concentration and explicit rates of convergence for the CLTs.
		
		\section*{Acknowledgements} This project was initiated during the program "\textit{Random Matrices and Scaling Limits}" at the Mittag-Leffler Institute. C.D.G and G.L. acknowledge the support of the starting grant 2022-04882 from the Swedish Research Council and of the starting grant from the Ragnar Söderbergs Foundation. 		
		\section{Upgraded recentring and global fluctuations} 
	\label{sec:bound Delta}
	The goal of this section  is to study the global fluctuations of multi-linear statistics, in particular to obtain an asymptotic expansion of their mean.
	In what follows, we rely on the following notations, for a function $\psi_n\in\mcal{C}_{\mrm{loc}}^{r}(\R^{n})$  with $r\ge 0$, we define for $\kappa>0$, 
	\begin{equation}\label{def:norm kappa}
		\|\psi_n\|_{\kappa,r}\defi \bigg\| \psi_n(x_1,\dots,x_n)\prod_{a=1}^{n}e^{-\kappa V(x_a)} \bigg\|_{\mcal{C}^r(\R^n)},\qquad\mcal{A}_{\kappa,r}\defi\left \{\psi_n\in\mcal{C}_{\mrm{loc}}^{r}(\R^{n}),\,\|\psi_n\|_{\kappa,r}<\infty\right \}
	\end{equation}
	and, for $\psi_n\in\mcal{A}_{\kappa,r}$ we set the following notation for $n$-multilinear statistics:
	\begin{equation}\label{eq:def bracket}
		\braket{\psi_n}\defi\mbb{E}_N\left[\int_{\R^n}\psi_n(x_1,\dots,x_n)\cdot\prod_{i=1}^{n}\diff(\bm{\mu}_N-\mu_{\mrm{eq}})(x_i)\right].
	\end{equation}
	
	The following asymptotic expansion of multi-linear statistics is proved in Section \ref{subsec1:loopeq} using loop equations based on the methods developed in \cite{BoG1}, see \cite{guera2024asymptotics} for an overview of the method. 			
	﻿\begin{Prop}\label{prop:asympcorrelators}
		Let $n\geq1$, for all $i\geq1,j\geq0$, there exist linear maps $b_{i,j}^{(n)}:\mcal{C}^{\infty}([-1,1]^{n})\rightarrow\R$ such that for any $K\geq L\geq0$, there exist $r=r(n,K)\geq1$ such that for all $\psi_n\in\mcal{A}_{\kappa,r}$ for some $\kappa>0$, one has the expansion as $N\to\infty$:
		\begin{equation}\label{braexp}
			\braket{\psi_n}=\sum_{j=0}^{L}\sum_{i=0}^{K-j}\mbf{1}_{(i,j)\neq(0,0)}\dfrac{b_{i,j}^{(n)}(\psi_n)}{\alpha^iN^{j}}+O\left(\dfrac{\|\psi_n\|_{\kappa,r}}{\alpha^{K+1}}\right)+O\left(\dfrac{\|\psi_n\|_{\kappa,r}}{N^{L+1}}\right).
		\end{equation}
		Moreover, one has  $|b_{i,j}^{(n)}(\psi_n)| \lesssim \|\psi_n\|_{\mcal{C}^{q}([-1,1]^n)}$ for some $q=q(n,i,j) \le r$ and $b_{i,j}^{(n)} = 0$ if $i+j<n$. 
	\end{Prop}

	We can use this expansion to define a suitable recentring to $\bm{\mu}_N$ to prove a CLT (Theorem \ref{thm:main CLT global}) and the optimal moment estimates (Proposition \ref{prop: moment estimates} and \ref{prop: moment estimates infinity}) for smooth linear statistics.  In particular, we define a sequence of mean functionals which become more and more accurate when $\alpha\gg 1$. Just as the equilibrium measure, these functional are supported on the support $[-1,1]$, but they are not measures, rather Schwartz distributions.
	
	\begin{Def}\label{def:mup}
		Let $k\geq0$, we define the linear map $\mu_k$ defined on $\mcal{C}^{r}([-1,1])$ for $r$ large enough by:
		$$\mu_{k}\defi\mu_{\mrm{eq}}+\sum_{j=1}^{k}\dfrac{\nu_j}{\alpha^{j}},$$
		where we defined for all $i\geq1$, $\nu_i\defi b_{i,0}^{(1)}$.
	\end{Def}
	
	In Subsection \ref{subsec1:loopeq}, we state the loop equations and derive Proposition \ref{prop:asympcorrelators}. We also study the Stieltjes transforms of the linear functionals $\nu_k, \mu_k$, the goal being to make a connection with the \emph{local laws} proved with Section \ref{sec: local laws}. Finally in Subsection \ref{subsec2:CLT}, we obtain optimal moment estimates in  Proposition \ref{prop: moment estimates} and from it, we prove the global CLT Theorem \ref{thm:main CLT global} using Stein's method in the spirit of \cite{LambertLedouxWebb,angst2024sharp}. We also state and prove the free energy expansion Theorem \ref{thm:free energy} in this subsection.
	
	\subsection{Analysis of loop equations}\label{subsec1:loopeq}
	To apply the loop equations method developed in \cite{BoG1}, we need some basic \textit{a priori} concentration results. The proof of the following classical result can be found in Appendix \ref{app:apriori}.
	\begin{Prop}[A priori bound]\label{a priori bound} Let $n\geq1$, $\kappa>0$, then for all $\psi_n$ depending on $n$ variables:
		$$\left| \Braket{\psi_n}\right|\lesssim\left(\dfrac{\log N}{\alpha}\right) ^{\frac{n}{2}} \left \|\psi_n\right\|_{\kappa,n}.$$
	\end{Prop}

	To state the loop equations, we will also need the definitions of several operators.
	\begin{Def}\label{def:ope D}
		We define the operator $\Xi$ called the master operator which is defined for all $\phi\in\mcal{C}^{1}_\mrm{loc}(\R)$ and $x\in\R$ by:
		$$\Xi[\phi](x)\defi \dfrac{V'(x)}{2}\phi(x)-\int_{\R}\dfrac{\phi(x)-\phi(y)}{x-y}\diff \mu _{\mrm{eq}}(y).$$
		We define the operator $\mcal{D}$ for all $\phi\in\mcal{C}^{1}_\mrm{loc}(\R)$ and $x,y\in\R$ by:
		$$\mcal{D}[\phi](x,y)\defi
		\begin{cases}
			\displaystyle
			\dfrac{\phi(x)-\phi(y)}{x-y},
			& \text{if } x\neq y, \medskip \\
			\phi'(x)
			& \text{if } x=y,
		\end{cases}$$
	\end{Def}
	
	Under the Assumptions~\ref{assumptions} $(i)-(iii)$, the \emph{master operator} is invertible and there is a classical expression for $\Xi^{-1}$ in terms of the density $S$ of the equilibrium measure. The following proposition is classical and its proof can be found in \cite[Lemma 4.1]{dworaczekguera2025clt} or \cite[Lemma 3.3]{bekerman2018clt}.
	
	\begin{Prop}[Invertibility of $\Xi$]\label{prop:inverse Xi}
		Let $f\in\mcal{C}^{1}(\R)$, there exists a unique $g\in\mcal{C}^{0}(\R)$ such that $\Xi[g]=f-a_f$ with $a_f\defi\int_{-1}^{1}f(t)\frac{\diff t}{\sqrt{1-t^{2}}}$.  Moreover, defining $\Xi^{-1}[f]\defi g$, for all $j\geq0$,
		$$\|\Xi^{-1}[f]\|_{\mcal{C}^{j}(\R)}\lesssim\|f\|_{\mcal{C}^{j+1}(\R)}.$$
		Furthermore $\Xi^{-1}[f]_{|[-1,1]}$ is given by:
		\begin{equation}\label{eq:inverse Xi formula}
			\Xi^{-1}[f]_{|[-1,1]}=\dfrac{1}{\pi S(x)}\int_{-1}^{1}\dfrac{f(x)-f(y)}{x-y}\dfrac{\diff y}{\sqrt{1-y^2}}.
		\end{equation}
		where $S$ is given in \eqref{eq: S} and for all $j\geq1$, $\|\Xi^{-1}[f]\|_{\mcal{C}^{j}([-1,1])}\lesssim\|f\|_{\mcal{C}^{j-1}([-1,1])}.$ Furthermore for $|x|\geq1$:
		%
					$$\Xi^{-1}[f](x)=
							\displaystyle 
							\frac{\displaystyle\int_{-1}^1 \frac{\Xi^{-1}[f](t)}{x-t} \diff \mu_{\mrm{eq}}(t)+f(x)-a}{ \dfrac{V'(x)}{2}-\displaystyle\int_{-1}^1 \frac{\diff \mu_{\mrm{eq}}(t)}{x-t}} ,
							$$
							and if $f(x)=0$ for $|x|\leq M$ for $M>1$ then, $\Xi^{-1}[f](x)=0$ for all $|x|\leq M$.
	\end{Prop}
	To state the loop equations for multilinear statistics, we need to define several other operators.
	\begin{Def}\label{def:nu1 theta}Let $n\geq1$ and $\psi_n\in\mcal{C}^{2}_{\mrm{loc}}(\R^{n})$, we define:
		$$\tau[\psi_n](x_2,\dots,x_n)\defi\int_{-1}^{1}\partial_1\Xi^{-1}[\psi_n](x_1,x_2,\dots,x_n)\diff \mu_{\mrm{eq}}(x_1).$$
		Note that when $n=1$, $\partial_1\Xi^{-1}[\psi_1]=\Xi^{-1}[\psi_1]'$ and $\nu_1(\psi_1)$ is just a constant.
		Furthermore if $n\geq2$, we define the operator $\Theta$ by:
		$$\Theta[\psi_n](x_1,\dots,x_{n-1})\defi \sum_{k=2}^{n}\Xi^{-1}[\partial_k\psi_{n}](x_1,\dots,x_{k-1},x_1,x_k,\dots,x_{n-1}).$$
		We also define:
		$$\theta[\psi_n](x_2,\dots,x_{n-1})\defi\int_{-1}^{1}\Theta[\psi_{n}](x_1,\dots,x_{n-1})\diff \mu_{\mrm{eq}}(x_1).$$
	\end{Def}
	
	We are now ready to deduce the loop equations. We recall that $\mcal{A}_{\kappa,1}$ was defined in \eqref{def:norm kappa} and $\braket{\cdot}$ in \eqref{eq:def bracket}.
	
	\begin{Prop}[Loop equations]\label{prop:DSequations}Let $\kappa>0$, the level $1$ loop equation holds for all function $\psi_1\in\mc{A}_{\kappa,1}$ and takes the form:
		\begin{equation}\label{DSlvl1}
			\Braket{\psi_1}=	\left (\dfrac{1}{\alpha }-\dfrac{1}{2N}\right)\tau[\psi_1]
			+\left (\dfrac{1}{\alpha }-\dfrac{1}{2N}\right)\Braket{\Xi^{-1}[\psi_1]'}+\dfrac{1}{2}\Braket{\mcal{D} \Xi^{-1}[\psi_1]}.
		\end{equation}
		For all $\psi_n\in\mc{A}_{\kappa,1}$, the level $n\geq2$ loop equations reads:
		\begin{equation}\label{DSlvlnthm}
			\begin{aligned}
				\Braket{\psi_{n}}
				&= \left (\dfrac{1}{\alpha }-\dfrac{1}{2N}\right )
				\Braket{\tau[\psi_{n}]}+\left (\dfrac{1}{\alpha }-\dfrac{1}{2N}\right )
				\Braket{\partial_1\Xi^{-1}[\psi_{n}]}+ \dfrac{1}{2}\Braket{\mcal{D} \Xi^{-1}[\psi_{n}]}\\
				&\quad + \dfrac{1}{N\alpha}	\Braket{\theta[\psi_n]}+\dfrac{1}{N\alpha}
				\Braket{\Theta[\psi_{n}]}.
			\end{aligned}
		\end{equation}
		%
		
	\end{Prop}
	
	The proof of the previous theorem is rather classical. We postpone it in Appendix \ref{app:apriori}.
	
	It is well known since \cite{BoG1}, that using the loop equations and the \textit{a priori} bounds from Proposition \ref{a priori bound}, one can obtain the large $N$-expansion of $n$-linear statistics $\braket{\psi_n}$ Proposition \ref{prop:asympcorrelators}. 
	
	\begin{Rem}\label{rem:continuity}	Every operator involved in the above loop equations are continuous. Indeed the following bounds hold for all $\kappa>0$, $k\geq0$:
		$$\|\tau[\psi_n]\|_{\kappa,k}\lesssim \|\psi_n\|_{\kappa,k+2},\quad\quad\|\Xi^{-1}[\psi_n]\|_{\kappa,k}\lesssim \|\psi_n\|_{\kappa,k+1},\quad\quad\quad\|\mcal{D}[\psi_n]\|_{\kappa,k}\lesssim \|\psi_n\|_{\kappa,k+1}$$
		and
		$$\|\Theta[\psi_n]\|_{\kappa,k}\lesssim \|\psi_n\|_{\kappa,k+2},\quad\quad\|\theta[\psi_n]\|_{\kappa,k}\lesssim \|\psi_n\|_{\kappa,k+2}.$$
		These bounds are crucial to apply Proposition \ref{a priori bound} for brakets of the form $\Braket{\mcal{W}[\psi_n]}$ where $\mcal{W}$ is any finite word in $\tau$, $\partial_1$, $\Xi^{-1}$, $\mcal{D}$, $\Theta$ and $\theta$. Because of these bounds and Proposition \ref{a priori bound}, we have for any $n\geq1$,
		$$\braket{\mcal{W}[\psi_n]}\lesssim\dfrac{\log^{k/2} N}{\alpha^{k/2}}\|\psi_n\|_{\kappa,r}$$
		for some $r\geq0$ where $k$ is the number of variables of the function $\mcal{B}[\psi_n]$.
	\end{Rem}
	
	The proof of Proposition \ref{prop:asympcorrelators} goes by induction. The next lemma covers the case $L=0$ in Proposition \ref{prop:asympcorrelators}. Let $\psi_n\in\mcal{A}_{\kappa,r}$ for some $\kappa>0$ and $r$ sufficiently large, we define the space of linear combination of brackets with bounded coefficients:
	$$\mcal{B}(\psi_n)\coloneq\left\{\sum_{i=1}^{p}c_i\braket{A_i[\psi_n]},\,p\geq0,\,c_i(N,\alpha)=O(1),\,A_i\text{ finite word in }\tau,\partial_1,\Xi^{-1}, \mcal{D}, \Theta \text{ and } \theta\right\} $$
	
	\begin{Lem}\label{lem:L=0}
			Let $n\geq1$, for all $i\geq1$ there exists $b_{i,0}^{(n)}:\mcal{C}^{\infty}([-1,1]^{n})\rightarrow\R$ such that for all $K\geq0$, there exists $r=r(n,K)\geq1$ such that for all $\psi_n\in\mcal{A}_{\kappa,r}$ for some $\kappa>0$, one has the expansion as $N\to\infty$:
		\begin{equation}\label{eq: proof prop21 eq intermed}
			\braket{\psi_n}=\sum_{i=n}^{K}\dfrac{b_{i,0}^{(n)}(\psi_n)}{\alpha^{i}}+o\left(\dfrac{\|\psi_n\|_{\kappa,r}}{\alpha^{K}}\right)+\dfrac{1}{N}\msf{B}_{0,K}^{(n)}(\psi_n),
		\end{equation}
		where $\msf{B}_{0,K}^{(n)}(\psi_n)\in\mcal{B}(\psi_n)$. Moreover, one has $|b_{i,0}^{(n)}(\psi_n)| \lesssim \|\psi_n\|_{\mcal{C}^{q}([-1,1]^n)}$ for some $q=q(n,i,j) \le r$.
	\end{Lem}
	
	\begin{Pro}
				We prove this statement by induction on $K\geq1$.
		
		By Proposition \ref{prop:DSequations}, the second loop equations gives:
		\begin{align}\label{eq:2nd loop eq}
			\Braket{\psi_{2}}
			&= \left (\dfrac{1}{\alpha }-\dfrac{1}{2N}\right )
			\Braket{\tau[\psi_{2}]} +\left (\dfrac{1}{\alpha }-\dfrac{1}{2N}\right )
			\Braket{\partial_1\Xi^{-1}[\psi_{2}]}+ \dfrac{1}{2}\Braket{\mcal{D} \Xi^{-1}[\psi_{2}]}+\dfrac{1}{N\alpha}\nonumber
			\Braket{\Theta[\psi_{2}]}\\
			&\quad  + \dfrac{1}{N\alpha}\theta[\psi_2]\nonumber
			\\&=\dfrac{1}{\alpha}\left(\Braket{\tau[\psi_{2}]}+\Braket{\partial_1\Xi^{-1}[\psi_{2}]}\right)+ \dfrac{1}{2}\Braket{\mcal{D} \Xi^{-1}[\psi_{2}]}+\dfrac{1}{N}\msf{B}_{0,1}^{(2)}(\psi_2),
		\end{align}
		where $\msf{B}_{0,1}^{(2)}(\psi_2)\defi \tfrac{1}{\alpha}(\Braket{\Theta[\psi_{2}]}+\theta[\psi_2])-\tfrac{1}{2}(\Braket{\tau[\psi_{2}]}+\Braket{\partial_1\Xi^{-1}[\psi_{2}]})\in\mcal{B}(\psi_2)$. Furthermore, because of Proposition \ref{a priori bound}, of Remark \ref{rem:continuity} and since $\alpha\gg\log ^{k} N$ for all $k\geq1$, we have the following bounds:
		$$\dfrac{1}{\alpha }|\Braket{\tau[\psi_{2}]}|\lesssim\dfrac{\log^{1/2} N}{\alpha^{3/2}}\|\psi_2\|_{\kappa,r}=o\left(\dfrac{\|\psi_2\|_{\kappa,r}}{\alpha}\right),$$
		$$\dfrac{1}{\alpha }|\Braket{\partial_1\Xi^{-1}[\psi_{2}]}|\lesssim\dfrac{\log N}{\alpha^{3/2}}\|\psi_2\|_{\kappa,r}=o\left(\dfrac{\|\psi_2\|_{\kappa,r}}{\alpha}\right),$$
		and
		$$ |\Braket{\mcal{D} \Xi^{-1}[\psi_{2}]}|\lesssim\dfrac{\log^{3/2} N}{\alpha^{3/2}}\|\psi_2\|_{\kappa,r}=o\left(\dfrac{\|\psi_2\|_{\kappa,r}}{\alpha}\right).$$
		Using \eqref{eq:2nd loop eq} and collecting all the previous bounds, we deduce that for general $\psi_2\in\mcal{A}_{\kappa,r}$ for some $r$ large enough:
		\begin{equation}\label{eq:improved bound psi2}
			|\braket{\psi_2}-\dfrac{1}{N}\msf{B}_{0,1}^{(2)}(\psi_2)|= o\left (\dfrac{\|\psi_2\|_{\kappa,r}}{\alpha}\right )
		\end{equation} for some $r\geq0$. Now the first loop equation, by Proposition \ref{prop:DSequations} gives for $\psi_1\in\mcal{A}_{\kappa,r}$:
		\begin{align*}
			\Braket{\psi_1}&=\dfrac{\tau[\psi_1]}{\alpha }-\dfrac{\tau[\psi_1]}{2N}
			+\left (\dfrac{1}{\alpha }-\dfrac{1}{2N}\right)\Braket{\Xi^{-1}[\psi_1]'}+\dfrac{1}{2}\Braket{\mcal{D} \Xi^{-1}[\psi_1]}
			\\&=\dfrac{\tau[\psi_1]}{\alpha }+\dfrac{1}{\alpha }\Braket{\Xi^{-1}[\psi_1]'}+\dfrac{1}{2}\Braket{\mcal{D} \Xi^{-1}[\psi_1]}+\dfrac{1}{N}\widetilde{\msf{B}}_{0,1}^{(1)}(\psi_1),
		\end{align*}
		where $\widetilde{\msf{B}}_{0,1}^{(1)}(\psi_1)=\tfrac{-1}{2}\tau[\psi_1]-\tfrac{1}{2}\Braket{\Xi^{-1}[\psi_1]'}$.
		Using Proposition \ref{a priori bound} as before, we have:
		$$\dfrac{1}{\alpha }|\Braket{\Xi^{-1}[\psi_1]'}|\lesssim\dfrac{1}{\alpha}\cdot\dfrac{\log ^{1/2} N}{\alpha^{1/2}}\|\psi_1\|_{\kappa,r},$$
		and by \eqref{eq:improved bound psi2} since $\mcal{D} \Xi^{-1}[\psi_1]\in\mcal{A}_{\kappa,r}$ is a function of 2 variables, we have:
		$$\Braket{\mcal{D} \Xi^{-1}[\psi_1]}= o\left(\dfrac{\|\psi_1\|_{\kappa,r}}{\alpha}\right)+\dfrac{1}{N}\msf{B}_{0,1}^{(2)}(\mcal{D} \Xi^{-1}[\psi_1]).$$
		Finally, collecting the previous bounds and noticing that $\widetilde{\msf{B}}_{0,1}^{(1)}(\psi_1)+\tfrac{1}{2}\msf{B}_{0,1}^{(2)}(\mcal{D} \Xi^{-1}[\psi_1])\in\mcal{B}(\psi_1)$, we obtain:
		$$\Braket{\psi_1}=\dfrac{\tau[\psi_1]}{\alpha }+o\left(\dfrac{\|\psi_1\|_{\kappa,r}}{\alpha}\right)+\dfrac{1}{N}\widetilde{\msf{B}}_{0,1}^{(1)}(\psi_1)+\dfrac{1}{2N}\msf{B}_{0,1}^{(2)}(\mcal{D} \Xi^{-1}[\psi_1]).$$
		Finally, for $n\geq3$:
		\begin{equation}\label{eq: n th loop eq}
			\begin{aligned}
				\Braket{\psi_{n}}
				&= \dfrac{1}{N\alpha}	\Braket{\theta[\psi_n]}+\dfrac{1}{N\alpha}
				\Braket{\Theta[\psi_{n}]}+\left (\dfrac{1}{\alpha }-\dfrac{1}{2N}\right )
				\Braket{\tau[\psi_{n}]} +\left (\dfrac{1}{\alpha }-\dfrac{1}{2N}\right )
				\Braket{\partial_1\Xi^{-1}[\psi_{n}]}
				\\	\quad&\quad+ \dfrac{1}{2}\Braket{\mcal{D} \Xi^{-1}[\psi_{n}]}
				,
				\\&=\dfrac{1}{\alpha }
				\Braket{\tau[\psi_{n}]}+\dfrac{1}{\alpha }
				\Braket{\partial_1\Xi^{-1}[\psi_{n}]}+ \dfrac{1}{2}\Braket{\mcal{D} \Xi^{-1}[\psi_{n}]}+\dfrac{1}{N}\msf{B}_{0,1}^{(n)}(\psi_n)
				\\&=o\left(\dfrac{\|\psi_1\|_{\kappa,r}}{\alpha}\right)+\dfrac{1}{N}\msf{B}_{0,1}^{(n)}(\psi_n),
			\end{aligned}
		\end{equation}
		where the last lign comes from Proposition \ref{a priori bound}. Noticing that $|\tau[\psi_1]|\lesssim \|\psi_1\|_{\mc{C}^{q}([-1,1])}$ for some $q\geq1$ because of Proposition \ref{prop:inverse Xi} concludes the case $K=1$.
		
		We now assume that for all $n\geq1$ and $\psi_n\in\mcal{A}_{\kappa,r}$ for $r$ large enough, \eqref{eq: proof prop21 eq intermed} holds for $K$, we show that the property holds for $K+1$.
		
		Let $n\geq1$ such that $\left(\tfrac{\log N}{\alpha}\right)^{n/2}\ll \tfrac{1}{\alpha^{K+1}}$, thus for all $m\geq n-1$:
		\begin{equation}\label{eq: n th loop eq'}
			\begin{aligned}
				\Braket{\psi_{m}}
				&= \dfrac{1}{N\alpha}	\Braket{\theta[\psi_m]}+\dfrac{1}{N\alpha}
				\Braket{\Theta[\psi_{m}]}+\left (\dfrac{1}{\alpha }-\dfrac{1}{2N}\right )
				\Braket{\tau[\psi_{m}]} +\left (\dfrac{1}{\alpha }-\dfrac{1}{2N}\right )
				\Braket{\partial_1\Xi^{-1}[\psi_{m}]}
				\\	\quad&\quad+ \dfrac{1}{2}\Braket{\mcal{D} \Xi^{-1}[\psi_{m}]}
				,
				\\&=\dfrac{1}{\alpha }
				\Braket{\tau[\psi_{m}]}+\dfrac{1}{\alpha }
				\Braket{\partial_1\Xi^{-1}[\psi_{m}]}+ \dfrac{1}{2}\Braket{\mcal{D} \Xi^{-1}[\psi_{m}]}+\dfrac{1}{N}\widetilde{\msf{B}}_{0,K+1}^{(m)}(\psi_m)
			\end{aligned}
		\end{equation}
		The third term is a $m+1$ braket so by Proposition \ref{a priori bound}, we have:
		$$|\Braket{\mcal{D} \Xi^{-1}[\psi_m]}|\lesssim\left(\dfrac{\log N}{\alpha}\right)^{n/2}\|\psi_m\|_{\kappa,r}=o\left(\dfrac{\|\psi_m\|_{\kappa,r}}{\alpha^{K+1}}\right).$$
		By induction hypothesis and since they both exhibit a prefactor $\tfrac{1}{\alpha}$, we have:
		\begin{align*}
			\dfrac{1}{\alpha }
			\left(	\Braket{\tau[\psi_{m}]}+
			\Braket{\partial_1\Xi^{-1}[\psi_{m}]}\right) &=\sum_{i=m-1}^{K}\dfrac{b_{i,0}^{(m-1)}(\tau[\psi_{m}])+b_{i,0}^{(m)}(\partial_1\Xi^{-1}[\psi_{m}])}{\alpha^{i+1}}+o\left(\dfrac{\|\psi_m\|_{\kappa,r}}{\alpha^{K+1}}\right)
			\\&\quad+\dfrac{1}{N}\msf{B}_{0,K+1}^{(m-1)}(\tau[\psi_{m}])+\dfrac{1}{N}\msf{B}_{0,K+1}^{(m)}(\partial_1\Xi^{-1}[\psi_{m}])
		\end{align*}
		where $\msf{B}_{0,K+1}^{(m-1)}(\tau[\psi_{m}])\in\mcal{B}(\tau[\psi_{m}])$ and $\msf{B}_{0,K+1}^{(m-1)}(\partial_1\Xi^{-1}[\psi_{m}])\in\mcal{B}(\partial_1\Xi^{-1}[\psi_{m}])$. Thus by \eqref{eq: n th loop eq'}, we obtain that for all $m\geq n-1$, 
		$$\begin{aligned}
			\Braket{\psi_{m}}&=\sum_{i=m}^{K+1}\dfrac{b_{i-1,0}^{(m-1)}(\tau[\psi_{m}])+b_{i-1,0}^{(m)}(\partial_1\Xi^{-1}[\psi_{m}])}{\alpha^{i}}+o\left(\dfrac{\|\psi_m\|_{\kappa,r}}{\alpha^{K+1}}\right)
			\\&\quad+\dfrac{1}{N}\left(\widetilde{\msf{B}}_{0,K+1}^{(m)}(\psi_m)+\dfrac{1}{\alpha}\msf{B}_{0,K+1}^{(m-1)}(\tau[\psi_{m}])+\dfrac{1}{\alpha}\msf{B}_{0,K+1}^{(m)}(\partial_1\Xi^{-1}[\psi_{m}])\right) 
		\end{aligned}$$
		and since $\widetilde{\msf{B}}_{0,K+1}^{(m)}(\psi_m)+\tfrac{1}{\alpha}\msf{B}_{0,K+1}^{(m-1)}(\tau[\psi_{m}])+\tfrac{1}{\alpha}\msf{B}_{0,K+1}^{(m)}(\partial_1\Xi^{-1}[\psi_{m}])\in\mcal{B}(\psi_m)$, it establishes the induction for $m\geq n-1$ . By the induction hypothesis, we know that both $b_{i-1,0}^{(m-1)}(\tau[\cdot])]$ and $b_{i-1,0}^{(m)}(\partial_1\Xi^{-1}[\cdot])$ are continuous for $\|\cdot\|_{\mc{C}^{q}([-1,1]^{m})}$ and linear for $i\in\llbracket m,K+1\rrbracket$. This implies that $b_{i,0}^{(m)}(\cdot)$ is continuous for the same norm. For $m=n-2$, taking the loop equation at level $m$, we have similarly as \eqref{eq: n th loop eq'}:
		\begin{equation*}
				\Braket{\psi_{m}}=\dfrac{1}{\alpha }
				\Braket{\tau[\psi_{m}]}+\dfrac{1}{\alpha }
				\Braket{\partial_1\Xi^{-1}[\psi_{m}]}+ \dfrac{1}{2}\Braket{\mcal{D} \Xi^{-1}[\psi_{m}]}+\dfrac{1}{N}\widetilde{\msf{B}}_{0,K+1}^{(m)}(\psi_m).
		\end{equation*}
		Since we have now proved the statement for the $m+1$ brakets, and by induction hypothesis using the prefactors $\tfrac{1}{\alpha}$ of the first two terms in the RHS above, we obtain:
			$$\begin{aligned}
				\Braket{\psi_{m}}&=\sum_{i=m}^{K+1}\dfrac{b_{i-1,0}^{(m-1)}(\tau[\psi_{m}])+b_{i-1,0}^{(m)}(\partial_1\Xi^{-1}[\psi_{m}])}{\alpha^{i}}+\sum_{i=m+1}^{K+1}\dfrac{b_{i,0}^{(m)}(\mcal{D} \Xi^{-1}[\psi_{m}])}{\alpha^{i}}+o\left(\dfrac{\|\psi_m\|_{\kappa,r}}{\alpha^{K+1}}\right)
				\\&\quad+\dfrac{1}{N}\left(\widetilde{\msf{B}}_{0,K+1}^{(m)}(\psi_m)+\dfrac{1}{\alpha}\msf{B}_{0,K+1}^{(m-1)}(\tau[\psi_{m}])+\dfrac{1}{\alpha}\msf{B}_{0,K+1}^{(m)}(\partial_1\Xi^{-1}[\psi_{m}])+\dfrac{1}{2}\msf{B}_{0,K+1}^{(m)}(\mcal{D} \Xi^{-1}[\psi_{m}])\right).
			\end{aligned}$$
		Repeating the exact same argument succesively for $m=n-3$, and then $n-4$ etc... proves the statement for all $m\geq1$. This concludes the induction. \end{Pro}
	
	\begin{Pro}[of Proposition \ref{prop:asympcorrelators}]
		\textbf{1. Asymptotic expansion:} To show the statement, we first prove by induction on $L$ that for all $n\geq1$ and $\psi_n\in\mcal{A}_{\kappa,r}$ for $r$ large enough, $K\geq1$:
		$$\braket{\psi_n}=\sum_{j=0}^{L}\sum_{i=0}^{K-j}\mbf{1}_{(i,j)\neq(0,0)}\dfrac{b_{i,j}^{(n)}(\psi_n)}{\alpha^{i}N^{j}}+o\left(\dfrac{\|\psi_n\|_{\kappa,r}}{\alpha^{K}}\right)+\dfrac{1}{N^{L+1}}\msf{B}_{L,K}^{(n)}(\psi_n)$$
		where $\msf{B}_{L,K}^{(n)}(\psi_n)\in\mcal{B}(\psi_n)$ and $b_{i,j}^{(n)}$ linear and continuous for $\|\cdot\|_{\mc{C}^{q}([-1,1]^{n})}$. The case $L=0$ is proved in Lemma \ref{lem:L=0}. We now assume that the above equation is true for $L\geq0$ and we show it for $L+1$. 
		
		Since $\msf{B}_{L,K}^{(n)}(\psi_n)\in\mcal{B}(\psi_n)$, there exists $p\geq0$, $(c_\ell)_{\ell=1}^{p}\in\R^{\ell}$ and $(A_\ell)_{\ell=1}^{p}$ finite words in $\tau,\partial_1,\Xi^{-1}, \mcal{D}$, $\Theta$ and $\theta$ with $A_\ell[\psi_n]$ possibly depending on 0 variables and $c_\ell=O(1)$, such that:
		\begin{align*}
			\msf{B}_{L,K}^{(n)}(\psi_n)&=\sum_{\ell=1}^{p}c_\ell\braket{A_\ell[\psi_n]}
			\\&=\sum_{i=0}^{K-L-1}\dfrac{1}{\alpha^{i}}\left(\sum_{\ell=1}^{p}c_\ell b_{i,0}(A_\ell[\psi_n])\right) +\dfrac{1}{N}\left( \sum_{\ell=1}^{p}c_\ell \msf{B}_{0,K-L-1}^{(n)}(A_\ell[\psi_n])\right)+o\left(\dfrac{\|\psi_n\|_{\kappa,r}}{\alpha^{K-L-1}}\right).
		\end{align*}
		Note that the second equality holds because of the induction hypothesis and that the sum starts at $i=0$, because the sum might contain 0-brakets. The linearity and continuity of $b_{i,j}^{(n)}(\psi_n)$ follows from the one of $b_{i,0}^{(k_\ell)}(A_\ell[\cdot])$ for $\ell\in\llbracket1,p\rrbracket$. It is straightforward to verify that the second term belongs to $\mcal{B}(\psi_n)$, this concludes the induction. Finally, by construction, any element of $\mcal{B}(\psi_n)$ is a $O(1)$ which proves that the expansion of Proposition \ref{prop:asympcorrelators} holds.

		\textbf{2. $b_{i,j}^{(n)}=0$ for $i+j<n$.} In this proof, we only work with smooth functions with apropriate growth at infinity such that the expansion  proved above holds at all order and don't specify the dependance in $\psi_n$ in our estimates. We show by induction that for $\alpha=\beta N$ with fixed $\beta>0$ and all $n\geq1$, $\braket{\psi_n}=O\left(N^{-n}\right)$. This is enough to show that $b_{i,j}^{(n)}(\psi_n)=0$ for $i+j<n$ because the obtained expansion gives: 
		$$\braket{\psi_{n}}=\sum_{k=1}^{n-1}\dfrac{1}{N^{k}}\sum_{i=0}^{k}\beta^{-i}b_{i,k-i}^{(n)}(\psi_n)+O\left(N^{-n}\right)=O\left(N^{-n}\right).$$
		Thus for every $k<n$ and every $\beta>0$:
		$$\sum_{i=0}^{k}\beta^{-i}b_{i,k-i}^{(n)}(\psi_n)=0.$$
		This implies that $b_{i,j}^{(n)}(\psi_n)=0$ for $i+j<n$.
		
		We now go back to the proof that for $\beta>0$ independent of $N$, $\braket{\psi_n}=O\left(N^{-n}\right)$. By Lemma \ref{lem:L=0} with $n=1$, $\braket{\psi_1}=O(1/N)$. Now assume that for all $k< n$, $\braket{\psi_k}=O(N^{-k})$. Using the $2n$-th loop equation:
		\begin{equation*}
			\begin{aligned}
				\Braket{\psi_{2n}}
				&= \dfrac{1}{N}\left (\dfrac{1}{\beta }-\dfrac{1}{2}\right )
				\Braket{\tau[\psi_{2n}]} +\dfrac{1}{N}\left (\dfrac{1}{\beta }-\dfrac{1}{2}\right )
				\Braket{\partial_1\Xi^{-1}[\psi_{2n}]}+ \dfrac{1}{2}\Braket{\mcal{D} \Xi^{-1}[\psi_{2n}]}\\
				&\quad + \dfrac{1}{N^2\beta}	\Braket{\theta[\psi_{2n}]}+\dfrac{1}{N^2\beta}
				\Braket{\Theta[\psi_{2n}]} .
			\end{aligned}
		\end{equation*}
		The \textit{a priori} bound Proposition \ref{a priori bound} gives the following estimates:
		$$\dfrac{1}{N}|
		\Braket{\tau[\psi_{2n}]}\lesssim\dfrac{\log ^{n-\frac{1}{2}}N}{N^{n+\frac12}},\quad\dfrac{1}{N} |
		\Braket{\partial_1\Xi^{-1}[\psi_{2n}]} |\lesssim\dfrac{\log ^{n}N}{N^{n+1}},\quad|\Braket{\mcal{D} \Xi^{-1}[\psi_{2n}]}|\lesssim\dfrac{\log ^{n+\frac12}N}{N^{n+\frac12}},$$
		and $$\dfrac{1}{N^2}|\Braket{\theta[\psi_{2n}]}|\lesssim\dfrac{\log ^{n-1}N}{N^{n+1}}, \quad\dfrac{1}{N^2}	|
		\Braket{\Theta[\psi_{2n}]}|\lesssim\dfrac{\log ^{n-\frac{1}{2}}N}{N^{n+\frac32}}.$$
		This leads to the improved bound $|\braket{\psi_{2n}}|\lesssim N^{-n-\frac12}\log ^{n+\frac12}N$. Applying to the $(2n-1)$-th loop equation the previous control for $\Braket{\mcal{D} \Xi^{-1}[\psi_{2n-1}]}$ and the \textit{a priori} bound for the remaining terms leads to the estimate $\braket{\psi_{2n-1}}=o(N^{-n})$. Doing the same for the $2n-2$, $2n-3$ ... $n$-th equation and applying the best bound between the improved bound and the bound from the induction hypothesis when available or the \textit{a priori} bound otherwise leads to $\braket{\psi_n}=O(N^{-n})$. This concludes the proof.
	\end{Pro}
	
	\begin{Rem}
		We proved Proposition \ref{prop:asympcorrelators}, in the case $\alpha\gg(\log N)^{k}$ for all $k\geq1$ but up to a slight modification of the proof, the result remains true for $\alpha\gg(\log N)^{1+\varepsilon}$ for some $\varepsilon>0$. 
	\end{Rem}
	
			The goal of the following proposition is to give a tensor formula for the linear maps $b_{i,0}^{(n)}$ appearing in the expansion. The latter is new to the knowledge of the authors and is crucial to show that the linear maps $\mu_k$, given in Definition \ref{def:mup} constitute the appropriate recentring that allows to obtain the optimal concentration bound and the global CLT. 	It also establishes a constraint on the measures $\nu_i$ similar to  \begin{equation}\label{eq:contrainst mu infinity}
				\int_{\R}V'(x)f(x)\diff\mu_{\mrm{eq}}-\iint_{\R^{2}}\dfrac{f(x)-f(y)}{x-y}\diff\mu_{\mrm{eq}}(x)\diff\mu_{\mrm{eq}}(y)=0
			\end{equation} and an analyticity preserving property.
	\begin{Prop}\label{prop: nu holomoprhic}
		For all $n\geq1$ and $f\in\mcal{C}^{r}([-1,1])$ for $r$ large enough:
		\begin{equation}\label{eq:constrait nu i}
			\nu_n(f)=\nu_{n-1}\left(\Xi^{-1}[f]'\right)+\dfrac{1}{2}\Bigg(\sum^{n-1}_{k=1}\nu_k\otimes\nu_{n-k}\Bigg) \left(\mcal{D} \Xi^{-1}[f]\right) 
		\end{equation}
		Furthermore, denoting the set of $n$-compositions of by $$\mscr{C}_{k}^{n}\defi\bigg\{\bm{i}=(i_1,\dots,i_n),\,\sum_{j=1}^{n}i_j=k,\,i_j\geq1\bigg\},$$
		we have for all $m,n\geq1$, $$b_{m,0}^{(n)}=\sum_{\bm{i}\in\mscr{C}_m^{n}}\bigotimes_{k=1}^{n}\nu_{i_k}.$$
		Finally, let $\mcal{O}\subset\C$ an open set, and for all $z\in\mcal{O}$, let $f_z\in\mcal{A}_{\kappa,r}$ such that for all $\lambda\in\R$, $z\in\mcal{O}\mapsto f_z(\lambda)$ is analytic. Thus $z\in\mcal{O}\mapsto\nu_n(f_z)$ is analytic for all $n\geq0$.
	\end{Prop}
	
	\begin{Pro} \textbf{1. Tensor formula for $b_{m,0}^{(n)}$ and constraints on $\nu_i$'s.}
		The scheme of the proof is represented in Figure \ref{fig: triangular induction}.
		\begin{figure}[h]
			\centering
			\begin{tikzpicture}[
				node distance=2cm,
				every node/.style={font=\small},
				box/.style={minimum width=1.6cm, minimum height=0.8cm},
				arr/.style={->, thick, green!60!black}
				]
				
				\node (p1) at (0,0) {$\psi_1$};
				\node (p2) at (2,0) {$\psi_2$};
				\node (p3) at (4,0) {$\psi_3$};
				\node (p4) at (6,0) {$\psi_4$};
				\node (p5) at (8,0) {$\psi_5$};
				\node (p6) at (10,0) {$\psi_6$};
				
				\node[left] at (-1,-1) {$b_1$};
				\node[left] at (-1,-2.5) {$b_2$};
				\node[left] at (-1,-4) {$b_3$};
				\node[left] at (-1,-5.5) {$b_4$};
				\node[left] at (-1,-7) {$b_5$};
				
				\node[box] (v1) at (0,-1) {$\nu_1$};
				\node[box] (v2) at (0,-2.5) {$\nu_2$};
				\node[box] (v3) at (0,-4) {$\nu_3$};
				\node[box] (v4) at (0,-5.5) {$\nu_4$};
				\node[box] (v5) at (0,-7) {$\nu_5$};
				
				\node[box] (e12) at (2,-1) {$\varnothing$};
				\node[box] (v1b) at (2,-2.5) {$\nu_1^{\otimes2}$};
				\node[box] (c22) at (2,-4) {$2\nu_1\nu_2$};
				\node[box] (c23) at (2,-5.5) {\scriptsize{$2\nu_1\nu_3+\nu_2^{\otimes2}$}};
				\node[box] (c24) at (2,-7) {$\substack{2\nu_1\nu_4
						\\+2\nu_2\nu_3}$};
				
				\node[box] (e13) at (4,-1) {$\varnothing$};
				\node[box] (e23) at (4,-2.5) {$\varnothing$};
				\node[box] (v1c) at (4,-4) {$\nu_1^{\otimes3}$};
				\node[box] (c32) at (4,-5.5) {$3\nu_1^{\otimes2}\nu_2$};
				\node[box] (c33) at (4,-7) {$\substack{3\nu_1^{\otimes2}\nu_3
						\\+3\nu_1\nu_2^{\otimes2}}$};
				
				\node[box] (e14) at (6,-1) {$\varnothing$};
				\node[box] (e24) at (6,-2.5) {$\varnothing$};
				\node[box] (e34) at (6,-4) {$\varnothing$};
				\node[box] (v1d) at (6,-5.5) {$\nu_1^{\otimes4}$};
				\node[box] (c42) at (6,-7) {$4\nu_1^{\otimes3}\nu_2$};
				
				\node[box] (e15) at (8,-1) {$\varnothing$};
				\node[box] (e25) at (8,-2.5) {$\varnothing$};
				\node[box] (e35) at (8,-4) {$\varnothing$};
				\node[box] (e45) at (8,-5.5) {$\varnothing$};
				\node[box] (v1e) at (8,-7) {$\nu_1^{\otimes 5}$};
				
				\node[box] (e16) at (10,-1) {$\varnothing$};
				\node[box] (e26) at (10,-2.5) {$\varnothing$};
				\node[box] (e36) at (10,-4) {$\varnothing$};
				\node[box] (e46) at (10,-5.5) {$\varnothing$};
				\node[box] (v1f) at (10,-7) {$\varnothing$};
				
				
				
				\draw[arr] (e12) -- (v1b);
				\draw[arr] (e23) -- (v1c);
				\draw[arr] (e34) -- (v1d);
				\draw[arr] (e45) -- (v1e);
				\draw[arr] (v1b) -- (c22);					
				\draw[arr] (c22) -- (c23);		
				\draw[arr] (c23) -- (c24);						
				\draw[arr] (v1c) -- (c32);
				\draw[arr] (c32) -- (c33);
				\draw[arr] (v1d) -- (c42);
				
				\draw[arr] (v1c) -- (c22);
				\draw[arr] (v1d) -- (c32);
				\draw[arr] (v1e) -- (c42);
				\draw[arr] (c32) -- (c23);
				\draw[arr] (c42) -- (c33);
				\draw[arr] (c33) -- (c24);

				\draw[arr] (v1) -- (v1b);
				\draw[arr] (v2) -- (c22);
				\draw[arr] (v3) -- (c23);
				\draw[arr] (c23) -- (c33);
				\draw[arr] (c32) -- (c42);
				\draw[arr] (c22) -- (c32);
				\draw[arr] (v4) -- (c24);
				
				\draw[arr] (v1b) -- (v1c);

				\draw[arr] (v1c) -- (v1d);

				\draw[arr] (v1d) -- (v1e);

			\end{tikzpicture}
			\caption{Coefficients $b_i(\psi_j)=b_{i,0}^{(j)}(\psi_j)$ where $\psi_j$ is a function of $j$ variables. We prove by induction that the result holds on the $k$-th diagonal \textit{i.e.} for $(b_{n+k}(\psi_n))_{n\geq1}$. For each $n,k$, the contributions to the $n$-th term on the $k$-th diagonal $b_{n+k}(\psi_n)$ are given by the previous term on the $k$-th diagonal $b_{n+k-1,0}^{(n-1)}$ (the top-left term), and two terms on the $(k-1)$-th diagonal, a $b_{n+k-1,0}^{(n)}$ (the top-middle term) and by a $b_{n+k,0}^{(n+1)}$-term (the right term) which is summarized in \eqref{eq:figure}.}
			\label{fig: triangular induction}
		\end{figure}
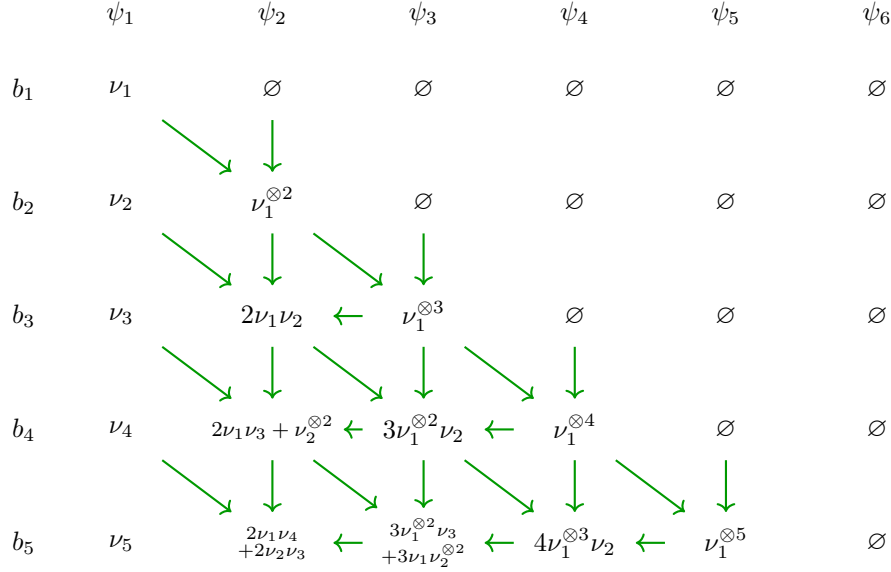We first show that, $b_{n,0}^{(n)}=\nu_1^{\otimes n}$. By Proposition \ref{prop:asympcorrelators} the $n$-th loop equation (see \eqref{DSlvlnthm}), we have for $\psi_{n}:(x_1,\dots,x_{n})\mapsto\prod_{i=1}^{n}\phi_i(x_i)$ with for all $i\in\llbracket1,n\rrbracket$, $\phi\in\mcal{A}_{\kappa,r}$ for $\kappa>0$ and $r$ large enough:
		$$
		\Braket{\psi_{n}}
		=  \dfrac{1}{\alpha }
		\braket{\tau\big[\bigotimes_{j=1}^{n}\phi_j\big]} +o\left(\dfrac{1}{\alpha^{n}}\right)=\dfrac{\tau[\phi_1]}{\alpha}\braket{\bigotimes_{j=2}^{n}\phi_j}+o\left(\dfrac{1}{\alpha^{n}}\right)=\dfrac{\displaystyle\prod_{j=1}^{n}\tau[\phi_i]}{\alpha^{n}}+o\left(\dfrac{1}{\alpha^{n}}\right).
		$$
		Thus since $\tau[\phi_1]=\nu_1(\phi_1)$, $b_{n,0}^{(n)}(\psi_{n})=\nu_1^{\otimes n}(\psi_{n}).$
		We conclude by density of $\mcal{C}^{r}([-1,1])^{\otimes n}$ in 
		$\allowbreak\mcal{C}^{r}([-1,1]^{n})$ and the continuity of $b_{n,0}^{(n)}$ for $\|\cdot\|_{\mc{C}^{r}([-1,1])}$. We now proceed by induction and assume that:
		$$b_{j,0}^{(n)}=\sum_{\bm{i}\in\mscr{C}_j^{n}}\bigotimes_{k=1}^{n}\nu_{i_k}$$
		holds true for the $k$ first diagonals in Figure \ref{fig: triangular induction} \textit{i.e.} for the linear maps $b_{n+j,0}^{(n)}$ for $j\in\llbracket0,k-1\rrbracket$ and all $n\geq1$. We now show that this decomposition holds for the diagonal $b_{n+k,0}^{(n)}$ for all $n\geq1$. We only treat the case of functions of the form $\psi_{n}=\bigotimes_{j=1}^{n}\phi_{j}$ where each $\phi_{j}\in\mcal{C}^{r}([-1,1])$ for $r$ large enough and conclude on the general case by density.
		
		By definition, $b_{1+k,0}^{(1)}=\nu_{k+1}$. We proceed by induction and assume that the result is true for $n-1$ and show it for $n$.  By Proposition \ref{prop:asympcorrelators}, identifying the coefficients in \eqref{DSlvl1}, we have for all $i\geq1$:
		\begin{equation*}
			b_{i,0}^{(1)}=\nu_i=\mbf{1}_{i=1}\nu_{1}+\mbf{1}_{i>1}\nu_{i-1}\left(\Xi^{-1}[\cdot]'\right) +\dfrac{1}{2}b_{i,0}^{(2)}\left(\mcal{D} \Xi^{-1}[\cdot]\right).
		\end{equation*}
		By the induction hypothesis (that the result holds on the $k$-first diagonals), for all $j\in\llbracket0,k-1\rrbracket$,
		$$b_{2+j,0}^{(2)}\left(\mcal{D} \Xi^{-1}[\cdot]\right)=\sum_{\bm{i}\in\mscr{C}_{j+2}^{2}}\bigotimes_{k=1}^{2}\nu_{i_k}\left(\mcal{D} \Xi^{-1}[\cdot]\right)=\sum_{k=1}^{j+1}\nu_i\otimes\nu_{j+2-i}\left(\mcal{D} \Xi^{-1}[\cdot]\right)$$
		from which we deduce that for all $i\in\llbracket1,k+1\rrbracket$:
		\begin{equation}\label{eq:constraint nui}
			\nu_i=\mbf{1}_{i=1}\nu_{1}+\mbf{1}_{i>1}\nu_{i-1}\left(\Xi^{-1}[\cdot]'\right) +\dfrac{1}{2}\sum_{\ell=1}^{i-1}\nu_\ell\otimes\nu_{i-\ell}\left(\mcal{D} \Xi^{-1}[\cdot]\right).
		\end{equation}
		Using again Proposition \ref{prop:asympcorrelators} and identifying the coefficients in the $n$-th loop equation \eqref{DSlvlnthm}, we obtain:
		\begin{align}\label{eq:figure}
			b_{n+k,0}^{(n)}(\psi_{n})&=\tau[\phi_1]b_{n+k-1,0}^{(n-1)}\Big(\bigotimes_{j=2}^{n}\phi_j\Big) +b_{n+k-1,0}^{(n)}\Big(\partial_1\Xi^{-1}[\phi_1]\bigotimes_{j=2}^{n}\phi_j\Big) 
			\\&\quad+ \dfrac{1}{2}b_{n+k,0}^{(n+1)}\Big(\mcal{D} \Xi^{-1}[\phi_{1}]\bigotimes_{j=2}^{n}\phi_j\Big) \nonumber.
		\end{align}
		Now, using the induction hypothesis (that the result holds true for the $(n-1)$-th coefficient on the $k$-th diagonal and on the whole $(k-1)$-th diagonal) and that $\tau[\phi_1]=\nu_1(\phi_1)$, leads to:
		\begin{align*}
			b_{n+k,0}^{(n)}(\psi_{n})
			&=\sum_{\bm{i}\in\mscr{C}^{n}_{n+k}}\ind{i_1=1}\prod_{j=1}^{n}\nu_{i_j}\left(\phi_{j}\right)+\sum_{\bm{i}\in\mscr{C}^{n}_{n+k-1}}\ind{i_1>0}\;\nu_{i_1}\left(\Xi^{-1}[\phi_{1}]'\right) \prod_{j=2}^{n}\nu_{i_j}\left(\phi_{j}\right)
			\\&\quad\quad+\dfrac{1}{2}\sum_{\bm{i}\in\mscr{C}^{n}_{n+k}}\sum_{\ell=1}^{i_1-1}\nu_\ell\otimes\nu_{i_1-\ell}\left(\mcal{D} \Xi^{-1}[\phi_{1}]\right) \prod_{j=2}^{n}\nu_{i_j}\left(\phi_{j}\right).
		\end{align*}
		We now use, use $i_1+1\rightarrow i_1$ in the second sum and \eqref{eq:constraint nui} to deduce that:
		\begin{align*}
			b_{n+k,0}^{(n)}(\psi_{n})&=\sum_{\bm{i}\in\mscr{C}^{n}_{n+k}}\Big(\ind{i_1>1}\;\nu_{i_1-1}\left(\Xi^{-1}[\phi_{1}]'\right)+\dfrac{1}{2}\sum_{\ell=1}^{i_1-1}\nu_\ell\otimes\nu_{i_1-\ell}\left(\mcal{D} \Xi^{-1}[\phi_{1}]\right)\Big)  \prod_{j=2}^{n}\nu_{i_j}\left(\phi_{j}\right)
			\\&\quad+\sum_{\bm{i}\in\mscr{C}^{n}_{n+k}}\ind{i_1=1}\prod_{j=1}^{n}\nu_{i_j}\left(\phi_{j}\right)
			\\&=\sum_{\bm{i}\in\mscr{C}^{n}_{n+k}}\bigotimes_{j=1}^{n}\nu_{i_j}\left(\psi_{n}\right).
		\end{align*}
		This establishes the result for $b_{n+k,0}^{(n)}$ and thus for the $(k+1)$-th diagonal by induction in Figure \ref{fig: triangular induction}, which concludes the induction on $k$. This also shows that \eqref{eq:constraint nui} holds for all $i\geq1$. 
		
		\textbf{2. Analyticity of $z\mapsto\nu_i(f_z)$.}  We first show that for all $x\in[-1,1]$, $z\mapsto\Xi^{-1}[f_z](x)$ is analytic. The map $g_{x,t}:z\in\mcal{O}\mapsto\dfrac{1}{\sqrt{1-t^2}}\dfrac{f_z(x)-f_z(t)}{x_1-t}$ is holomorphic for all $x\in[-1,1]$, $t\in(-1,1)$ and for all compact set $K\subset\mcal{O}$ and $z\in K$, $|g_{x,t}(z)|\leq \|\partial_xf_z\|_{\infty,K}$ which is $t$-integrable on $(-1,1)$. Thus the result follows by \eqref{eq:inverse Xi formula} and the analyticity under the integral sign Theorem. The analyticity of $z\mapsto\mcal{D} \Xi^{-1}[f_z](x_1,x_2)$ for all $x_1,x_2\in[-1,1]$, follows easily. Furthermore, the map $$h_{x,t}:z\in\mcal{O}\mapsto\dfrac{1}{\sqrt{1-t^2}}\dfrac{f_z(t)-f_z(x)-\partial_zf_z(x)(t-x)}{(x-t)^2}$$ is holomorphic for all $-1\leq t,x\leq1$ and for all compact set $K\subset\mcal{O}$ and $z\in K$, $|h_{x,t}(z)|\leq \|\partial_x^{2}f_z\|_{\infty,K}$, thus $z\mapsto\partial_{x}\Xi^{-1}[f_z](x)$ is analytic for all $x\in[-1,1]$.
		
		We now show the result by induction on $n$. For $n=1$, we have by Definition \ref{def:nu1 theta}:
		$$\nu_1(f_z)=\tau[f_z]=\int_{-1}^{1}\partial_x\Xi^{-1}[f_z](x)\diff\mu_{\mrm{eq}}(x).$$
		Thus since $z\mapsto\partial_x\Xi^{-1}[f_z](x)$ is holomorphic for all $x\in[-1,1]$ and $$|\partial_x\Xi^{-1}[f_z](x)|\leq\max_{(x,z)\in[-1,1]\times K}|\partial_x\Xi^{-1}[f_z](x)|$$ for all compact $K\subset\mcal{O}$ and $z\in K$, $z\in\mcal{O}\mapsto\nu_1(f_z)$ is holomorphic. The result follows by induction using \eqref{eq:constrait nu i}.
	\end{Pro}

	\begin{Rem}\begin{enumerate}
			\item The map $\nu_1$  can be made explicit in the Gaussian case, it is given by the mean in the Johansson's CLT \cite{johansson}. It was already shown in \cite{benaych2015poisson}, that:
			$$\bm{\mu}_N-\mu_{\mrm{eq}}\tend{N\rightarrow\infty}\nu_1=\dfrac{1}{2}\left(\delta_{-1}+\delta_{1}\right)-\dfrac{1}{\pi\sqrt{1-x^{2}}}\diff x $$
			where the convergence holds when tested against sufficiently smooth function. Furthermore $\nu_2$ is given by, using \eqref{eq:constraint nui}:
			$\nu_2(\psi)=\nu_1\left(\Xi^{-1}[\psi]'\right)+\dfrac{1}{2}\nu_{1}^{\otimes2}\left(\mcal{D} \Xi^{-1}[\psi]\right).$
			\item For example, for a symmetric function $\psi_3$ of 3 variables, $$b_{5,0}^{(3)}(\psi_3)=\left[3\nu_1^{\otimes2}\nu_3+3\nu_1\nu_2^{\otimes2}\right](\psi_3), $$
			see Figure \ref{fig: triangular induction} for the general pattern.
			\item From, $\Xi^{-1}[\widetilde{1}]=0$ on $[-1,1]$ by \eqref{eq:inverse Xi formula}, one can see by \eqref{eq:constrait nu i} that $\nu_k(\widetilde{1})=0$ for all $k\geq1$. This shows that for all $k\geq0$, $\mu_k(\widetilde{1})=1$.
		\end{enumerate}
		
	\end{Rem}
	The next proposition establishes an equation satisfied by $m_p$, the Stieltjes transforms of $\mu_p$. As similar equations are studied in Section \ref{sec: local laws} to establish local laws, the following establishes that the recentring used in the local laws are given by the $\mu_p$'s.
	
	\begin{Prop}\label{prop:mun stieltjes}Let $z\in\mc{U}_+$, denoting by $m_p(z)\defi\mu_p\left(\dfrac{1}{\cdot-z}\right) $, we have for all $p\geq0$ with the convention $m_{-1}=0$ and where $W_z:\lambda\mapsto\tfrac{V'(\lambda)-V'(z)}{\lambda-z}$:
		$$m_{p}(z)^{2}+V'(z)m_{2p}(z)+\mu_{2p}\left(W_z\right)+\dfrac{2}{\alpha}\partial_zm_{2p-1}(z)=0.$$
		Furthermore,  setting $\sigma_k(z)\defi\nu_k\left(\tfrac{1}{\cdot-z}\right)$, we have for all $k\geq1$:
		$$\sigma_k(z) =-\dfrac{\partial_z\sigma_{k-1}(z)+\dfrac{1}{2}\displaystyle\sum_{j=1}^{k-1}\sigma_j(z)\sigma_{k-j}(z)+\dfrac{\nu_k\left(W_z\right)}{2}}{S(z)\sqrt{z^{2}-1}}.$$
		Let $\delta>0$, then for all $z\in\mcal{U}_+$ such that $\re{z}\in(-1+\delta,1-\delta)$ and $p\geq0$, $\im m_p(z)\geq0$ for $\alpha$ large enough (depending on $\delta$, $p$).
	\end{Prop}
	
	\begin{Pro}
		We use \eqref{eq:constrait nu i} with $f=\Xi[g]$, Definition \ref{def:ope D} to get:
		\begin{equation}\label{eq:constraintnukproof}
			\dfrac{1}{2}\nu_k(V'g)-\nu_k\otimes\mu_{\mrm{eq}}\left(\mcal{D}[g]\right) =\nu_{k-1}\left(g'\right)+\dfrac{1}{2}\Bigg(\sum_{j=1}^{k-1}\nu_j\otimes\nu_{k-j}\Bigg) \left(\mcal{D}[g]\right).
		\end{equation}
		Combining the second term in the LHS with the sum in the RHS, multiplying by $\alpha^{-k}$ and summing for $k\in\llbracket1,2n\rrbracket$, we get:
		$$\mu_{2n}(V'g)-\mu_{\mrm{eq}}(V'g) =\dfrac{2}{\alpha}\mu_{2n-1} \left(g'\right)+\mu_n^{\otimes2}  \left(\mcal{D}[g]\right)-\mu_{\mrm{eq}}^{\otimes2} \left(\mcal{D}[g]\right).$$
		Finally using \eqref{eq:contrainst mu infinity}, taking $g(x)=(x-z)^{-1}$ and using:
		$$\mu_{2n} (V'g)=V'(z)m_{2n}(z)+\left(\sum_{k=0}^{2n}\dfrac{\nu_k\left(W_z\right)}{\alpha^{k}}\right),\quad\mu_{n}^{\otimes2}  \left(\mcal{D}[g]\right)=-m_n(z)^{2},\quad \partial_zm_{2n-1}(z)=-\mu_{2n-1} \left(g'\right)$$
		yields the first formula.
		
		For the second one, we use again \eqref{eq:constraintnukproof} to deduce that:
		$$	\dfrac{1}{2}\nu_k\left(W_z\right)+\left[\dfrac{V'(z)}{2}+\mu_{\mrm{eq}}\left(\dfrac{1}{\cdot-z}\right)\right] \sigma_k(z) =-\partial_z\sigma_{k-1}(z)-\dfrac{1}{2}\Bigg(\sum_{j=1}^{k-1}\sigma_j(z)\sigma_{k-j}(z)\Bigg).$$
		Thus  to conclude about the formula, we just use \eqref{roots}: $$V'(z)+2\mu_{\mrm{eq}}\left(\dfrac{1}{\cdot-z}\right)=S(z)\sqrt{z^{2}-1}$$ where the RHS is non zero close to the bulk. Finally, the positivity of $\im m_n(z)$ for all $n\geq0$ above the bulk, comes from the fact $m_0$ is the Stieltjes transform of $\mu_{\mrm{eq}}$ whose density is strictly positive and that $\nu_k \left(\dfrac{1}{\cdot-z}\right)$ is bounded above the bulk in $\|\cdot\|_{\mcal{C}^{k}}$-norm on $(-1+\delta,1-\delta)+\ii[0,+\infty)$. From this, it is easy to see by a simple induction and the fact $\nu_k\left(W_z\right)/[S(z)\sqrt{z^{2}-1}]$ is bounded at infinity, that $\sigma_{k}(z)$ is bounded above the bulk for all $k\geq0$. Thus for $\alpha$ large enough it yields the result.
	\end{Pro}

	\subsection{Concentration, global CLT and free energy}\label{subsec2:CLT}
	In this section, we derive all the main results of this section, namely the optimal moment estimates Proposition \ref{prop: moment estimates}, the global CLT Theorem \ref{thm:main CLT global} and the expansion of the free energy Theorem \ref{thm:free energy}.
	
	\begin{Def}[Recentring]\label{def:improved recentring}
		Let $k\geq1$ and $f\in\mcal{C}^{r}(\R)$ for $r$ large enough (depending on $k$), we define the linear map $L_N^{(k)}(f)$ by:
		$$L_N^{(k)}(f)\defi (\bm{\mu}_N-\mu_k)(f).$$
		We also define the anisotropy term by:
		$$A_N^{(k)}(f)\defi \left(\bm{\mu}_N-\mu_k\right)^{\otimes2}\left(\mcal{D}[f]\right) .$$
	\end{Def}
	
	Using \eqref{eq:constrait nu i}, we can rewrite the first loop equation by using the upgraded recentring $L_N^{(k)}$. The following recentered loop equation is needed to derive the optimal moment estimates of  Proposition \ref{prop: moment estimates} and the global CLT Theorem \ref{thm:main CLT global}.
	
	\begin{Prop}\label{prop:rewriting1stloopeq}
		Let $j\geq1$ and $g\in\mcal{C}^{r}(\R)$ for $r$ large enough, then the following holds true:
		\begin{multline*}
			-L_N^{(0)}(g)+\dfrac{1}{2}A_N^{(0)}(\Xi^{-1}[g])+ \left(\dfrac{1}{\alpha}-\dfrac{1}{2N}\right) \bm{\mu}_N(\Xi^{-1}[g]')
			\\=-L_N^{(j)}(g)+\dfrac{1}{2}A_N^{(\lfloor j/2\rfloor)}(\Xi^{-1}[g])+(\sum_{i=1}^{\lfloor j/2\rfloor}\frac{\nu_i}{\alpha^{i}}\otimes L_N^{(j-i)})\left(\mcal{D} \Xi^{-1}[g]\right) + \dfrac{1}{\alpha}  L_N^{(j-1)}(\Xi^{-1}[g]')-\dfrac{\bm{\mu}_N(\Xi^{-1}[g]')}{2N} .
		\end{multline*}
	\end{Prop}
	
	\begin{Pro}
		We prove the statement by induction. For $j=1$, we just use the fact that, by \eqref{eq:constrait nu i} with $n=1$: $$\nu_1(g)=\mu_{\mrm{eq}}(\Xi^{-1}[g]').$$
		Assume that the proposition is true for a $j\geq1$, defining:
		$$F_j\coloneqq\dfrac{1}{2}A_N^{(\lfloor j/2\rfloor)}(\Xi^{-1}[g])+(\sum_{i=1}^{(\lfloor j/2\rfloor)}\frac{\nu_i}{\alpha^{i}}\otimes L_N^{(j-i)})\left(\mcal{D} \Xi^{-1}[g]\right),$$
		we first notice that if $j=2\ell$ is even then:
		\begin{equation*}
			F_{2\ell}=F_{2\ell+1}+\dfrac{1}{\alpha^{2\ell+1}}(\sum_{i=1}^{\ell}\nu_i\otimes \nu_{2\ell+1-i})\left(\mcal{D} \Xi^{-1}[g]\right)
		\end{equation*}
		If $j=2\ell+1$ is odd, then using that for a symmetric function $g_2(x,y)$:
		$$L_N^{(k)\otimes2}(g_2)=L_N^{(k+1)\otimes2}(g_2)+\dfrac{2}{\alpha^{k+1}}\big(\nu_{k+1}\otimes L_N^{(k+1)}\big)(g_2)+\dfrac{1}{\alpha^{2(k+1)}}\nu_{k+1}^{\otimes2}(g_2),$$
		we obtain: 
		\begin{equation*}
			F_{2\ell+1}=F_{2\ell+2}+\dfrac{1}{\alpha^{2\ell+2}}(\sum_{i=1}^{\ell+1}\nu_i\otimes \nu_{2\ell+2-i})\left(\mcal{D} \Xi^{-1}[g]\right).
		\end{equation*}
		Thus for all $j\geq1$:
		\begin{equation}\label{eq:proof rewriting loopeq anis}
			F_j=F_{j+1}+\dfrac{1}{\alpha^{j+1}}\Big(\sum_{i=1}^{\lfloor \frac{j+1}{2}\rfloor}\nu_i\otimes \nu_{j+1-i}\Big)\left(\mcal{D} \Xi^{-1}[g]\right).
		\end{equation}
		Now, using \eqref{eq:proof rewriting loopeq anis}, we have:
		\begin{multline*}
			-L_N^{(j)}(g)+F_j+\dfrac{1}{\alpha} L_N^{(j-1)}(\Xi^{-1}[g]')=-L_N^{(j+1)}(g)+F_{j+1} + \dfrac{1}{\alpha}  L_N^{(j)}(\Xi^{-1}[g]')
			\\-\dfrac{1}{\alpha^{j+1}}\bigg(\nu_{j+1}(g)-\dfrac{1}{2}\sum_{i=1}^{\lfloor  \frac{j+1}{2}\rfloor}\nu_{j+1-i}\otimes\nu_i\left(\mcal{D} \Xi^{-1}[g]\right)-\nu_{j}(\Xi^{-1}[g]')\bigg).
		\end{multline*}
		Finally, using the induction hypothesis and \eqref{eq:constrait nu i} with $n=j+1$ conclude the proof.
	\end{Pro}
	
	We first state an \textit{a priori} concentration inequality appearing in \cite[Corollary 4.16]{guionnet2019asymptotics}, see \cite[Lemma 5.2]{dworaczekguera2025clt} for an alternative version more suited to our context. The adaptation to the regime $\alpha\gg(\log N)^{1+\varepsilon}$ is trivial.

	\begin{Lem}[Initial concentration]\label{lem:initial concentration}
	\textit{(i)}There exists $c,C>0$ such that for all $t>C\sqrt{\frac{\log N}{\alpha}}$:
		$$\P_N\left(d(\bm{\mu}_N,\mu_{\mrm{eq}})>t\right)\leq e^{-ct^{2}N\alpha}.$$
		\textit{(ii)} Furthermore, let $M>0$ such that $V(x)=x^{2}+C$ for all $|x|\geq M$ and $\chi\in\mcal{C}_c^{\infty}(\R)$ such that $\ind{[-M,M]}\leq\chi\leq\ind{[-2M,2M]}$, for all $g\in\mcal{C}^{0}(\R)$ such that $\|g\|_{\mcal{C}^{0}(\R)}\leq1$, there exists $c,C'>0$ such that for all $t\geq C'\sqrt{\frac{\log N}{\alpha}}$:
		$$\P_N\left(|\bm{\mu}_N(g(1-\chi))|>t\right)\leq e^{-ct^{2}N\alpha}.$$
	\end{Lem}
	
	\begin{Pro}\textit{(i)} We start with $\widetilde{\bm{\mu}}_{N,u}$ defined in Definition \ref{def:spacedconfig}, we have:
		$$d(\bm{\mu}_N,\mu_{\mrm{eq}})\leq d(\bm{\mu}_N,\widetilde{\bm{\mu}}_{N,u})+d(\widetilde{\bm{\mu}}_{N,u},\mu_{\mrm{eq}}).$$
		In the RHS, the first is term is a $O(1/N^{3})$ by Definition \ref{def:spacedconfig} and for a constant $C>0$:
		$$\Big|\int_\R f(x)\diff(\widetilde{\bm{\mu}}_{N,u}-\mu_{\mrm{eq}})(x)\Big|^{2}\leq C\Big|\int_\R|t|^{1/2}\widehat{f}(t)\dfrac{\overline{\widehat{(\mu_{N,u}-\mu_{\mrm{eq}}})}(t)}{|t|^{1/2}}\diff t\Big|\leq C\|f\|_{\msf{H}^{1/2}}^{2}\mfrak{D}^{2}[\widetilde{\bm{\mu}}_{N,u},\mu_{\mrm{eq}}].$$
		We thus deduce that by Proposition \ref{thm:concentration}, there exists $C,c>0$ such that for all $t>C\sqrt{\frac{\log N}{\alpha}}$:
		$$\P_N\left(d(\widetilde{\bm{\mu}}_{N,u},\mu_{\mrm{eq}})>t\right)\leq\P_N\Big(\mfrak{D}^{2}[\widetilde{\bm{\mu}}_{N,u},\mu_{\mrm{eq}}]>\dfrac{1}{C}(t-\dfrac{C}{N^{3}})^{2}\Big)\leq e^{-ct^{2}N\alpha}.$$
		This concludes the proof of \textit{(i)}. For the point \textit{(ii)}, we use again Markov inequality:
		$$\P_N\left(|\bm{\mu}_N(g(1-\chi))|>t\right)\leq e^{-t^{2}N\alpha}\E_N\left[\exp(N\alpha\bm{\mu}_N(g(1-\chi))^{2})\right].$$
		Furthermore, by Proposition \ref{thm:concentration}, we have for some $C'''>0$
		$$\E_N\left[\exp(N\alpha\bm{\mu}_N(g(1-\chi))^{2})\right]\leq e^{C'' N\log N}\left(\int_{\R}\exp\left(-\alpha \left(\frac{1}{2}V_{\mrm{eff}}(\lambda)-g(\lambda)\right)\right)\diff \lambda\right)^{N} $$
		Now since there exists $c>0$ (independent of $g$) such that for all $\lambda\geq 2M$, $\frac{1}{2}V_{\mrm{eff}}(\lambda)-g(\lambda)\geq c \lambda^{2}$, the integral converges and we have:
		$$\E_N\left[\exp(N\alpha\bm{\mu}_N(g(1-\chi))^{2})\right]\leq e^{C'' N\log N}\left(C+\int_{\R}\exp\left(-c\lambda^{2}\right)\diff \lambda\right)^{N}\leq e^{C''' N\log N},$$
		which allows to conclude that for constants $c'',C'''>0$:
			$$\P_N\left(|\bm{\mu}_N(f(1-\chi))|>t\right)\leq e^{-t^{2}N\alpha+C'''N\log N}\leq e^{-c''t^{2}N\alpha}.$$\qedsymbol{}
	\end{Pro}
	
	From the previous concentration inequalities, we deduce the following non-optimal moment estimates.
	
	\begin{Prop}[Initial moment estimates]\label{prop:inital moment estimate}
		\textit{(i)} Let $f\in\mcal{C}^{1}(\R)$ such that $\|f\|_{\msf{H}^{1/2}}<\infty$, then for all $q\geq1$:
		$$\|L_N^{(0)}(f)\|_{q}\lesssim \left(\|f\|_{\mcal{C}^{1}(\R)}+\|f\|_{\msf{H}^{1/2}}\right) \left(\sqrt{\dfrac{q}{N\alpha}}+\sqrt{\dfrac{\log N}{\alpha}}\right ).$$
		\textit{(ii)} Let $g\in\mcal{C}^{0}(\R)$, under the same assumptions as in \textit{(ii)} in Lemma \ref{lem:initial concentration}, we have for all $q\geq1$:
			$$\|\bm{\mu}_N(g(1-\chi))\|_{q}\lesssim \|g\|_{\mcal{C}^{0}(\R)}\left(\sqrt{\dfrac{q}{N\alpha}}+\sqrt{\dfrac{\log N}{\alpha}}\right ).$$
	\end{Prop}
	\begin{Pro}
Setting $\widetilde{f}=f/(\|f\|_{\msf{H}^{1/2}}+\|f\|_{\mcal{C}^{1}(\R)})$, we have for a constant $C>0$ independent of $f$:		\begin{align*}
		\E[|L_N^{(0)}(\widetilde{f})|^{q}]&=q\int_{0}^{+\infty}x^{q-1}\P\left(|L_N^{(0)}(\widetilde{f})|>x\right)\diff x
		\\&\lesssim C^{q}q\left(\dfrac{\log N}{\alpha}\right)^{q/2}+C^{q}q \int_{C\sqrt{\frac{\log N}{\alpha}}}^{+\infty}x^{q-1}e^{-cx^{2}N\alpha}\diff x
		\\&\lesssim C^{q}q\left(\dfrac{\log N}{\alpha}\right)^{q/2}+C^{q}q \left(\dfrac{1}{N\alpha}\right)^{q/2}\Gamma(q/2)
		\\&\lesssim C^{q}q\left(\dfrac{\log N}{\alpha}\right)^{q/2}+C^{q}q\left(\dfrac{q}{N\alpha}\right)^{q/2}.
	\end{align*}
	We thus, obtain:
	$$\|L_N^{(0)}(f)\|_q\lesssim (\|f\|_{\msf{H}^{1/2}}+\|f\|_{\mcal{C}^{1}(\R)})\left(\sqrt{\dfrac{q}{N\alpha}}+\sqrt{\dfrac{\log N}{\alpha}}\right).$$
	This concludes the proof of \textit{(i)}. Setting $\widetilde{g}=g/\|g\|_{\mcal{C}^{0}(\R)}$, the proof of \textit{(ii)} is identical.
	\end{Pro}

	We now use Stein's method and follow the strategy developed in \cite{LambertLedouxWebb,angst2024sharp}. It is based on the generator of the Dyson Brownian motion $	\mc{L} $ which acts on smooth functions on $\R^{N}$ by:
	\begin{equation}\label{eq:def L}
		\mc{L} [F]\defi-\sum_{i=1}^N\partial_{\lambda_i}^2F+\dfrac{\alpha}{2}\sum_{i=1}^NV'(\lambda_i)\partial_{\lambda_i}F-\dfrac{\beta}{2}\sum_{i\neq j}\dfrac{\partial_{\lambda_i}F-\partial_{\lambda_j}F}{\lambda_i-\lambda_j}.
	\end{equation}
	Note that $\mbb{E}_N\left[	\mc{L} [F]\right]=0$ for $F(\bm{\lambda})=\bm{\mu}_N(f)$ with $f\in\mcal{C}_c^{\infty}(\R)$ is equivalent to the first loop equation Proposition \ref{prop:DSequations}. We write $\mathrm{Dom}(\mc{L} )$ for the domain of the operator $\mc{L} $, the set of $F \in \mcal{C}^{2}(\R^{N})$ such that $\mc{L} [F] \in L^2(\P_N)$, it is very clear that $\bm{\lambda}\in\R^{N}\mapsto\bm{\mu}_N(f)\in \mathrm{Dom}(\mc{L} )$ for $f\in\mcal{C}^{2}(\R)$.
	
	\begin{Prop}[Improved moment estimates]\label{prop: moment estimates}
			Let $k\geq0$, there exists $\ell_k\geq1$ such that for all $g\in\mcal{C}_c^{\ell_k}(\R)$ and for all $q\geq1$:
		$$	\|L_N^{(k)}(g)\|_q \lesssim\|g\|_{\mcal{C}^{\ell_k}(\R)}\left(\sqrt{\dfrac{q}{N\alpha}}+\dfrac{1}{\alpha^{k+1}}\right).$$
	\end{Prop}
	
	\begin{Pro}
			\textbf{1. Case $k=0$.} Taking $F=N\bm{\mu}_N(f)$ in \eqref{eq:def L}, it is easy to see that by using \eqref{eq:contrainst mu infinity}:
		\begin{align}\label{eq:eq boucle L}
			\dfrac{\mc{L} [F]}{N\alpha}&=-\left(\dfrac{1}{\alpha}-\dfrac{1}{2N}\right) \bm{\mu}_N(f'')+\dfrac{1}{2}\bm{\mu}_N(V'f')-\dfrac{1}{2}\bm{\mu}_N^{\otimes2}(\mcal{D}[f'])\nonumber
			\\&=L_N^{(0)}(\Xi[f'])-\dfrac{1}{2}A_N^{(0)}(f')- \left(\dfrac{1}{\alpha}-\dfrac{1}{2N}\right) \bm{\mu}_N(f'').
		\end{align}
		One can then obtain the following decomposition:
		$$X_0=\dfrac{1}{N}\left(\mc{L} \left[F\right] +Z_0\right) ,$$
		where $X_0=\alpha L_N^{(0)}\left(\Xi[f']\right)$ and $$Z_0\defi N\alpha\left[\left(\dfrac{1}{\alpha}-\dfrac{1}{2N}\right) \bm{\mu}_N(f'')+\dfrac{1}{2}A_N^{(0)}(f')\right] .$$
		We now set $g=\Xi[f']$ and suppose that $g\in\mcal{C}_c^{\ell}(\R)$ for $\ell$ large enough, thus by Proposition \ref{prop:inverse Xi} and \ref{prop:bound_normq} with $\lambda=N$:
		\begin{equation}\label{eq:stein inequality}
					\|X_0\|_q \le \sqrt{\dfrac{q \|\Gamma_0 \|_{q/2}}{N}} +\dfrac{ \|Z_0\|_q }{N}\lesssim\alpha\left(\sqrt{\dfrac{q\|g\Xi^{-1}[g]'\|_{\infty}}{N\alpha}}+\|A_N^{(0)}(\Xi^{-1}[g])\|_q+\dfrac{\|\Xi^{-1}[g]'\|_{\infty}}{\alpha}\right)  , 
		\end{equation}
		where we noticed that
		$$\Gamma_0\defi \nabla F \cdot\nabla X_0=\dfrac{\alpha}{N}\sum_{i=1}^{N}g(\lambda_i)\Xi^{-1}[g]'(\lambda_i)=\alpha\bm{\mu}_N(g\Xi^{-1}[g]').$$
		To control $|A_N^{(0)}(g)|$, we use Helffer–Sjöstrand formula from Proposition
		\ref{prop:Helffer–Sjöstrand} with $\delta=1$ and a pseudo-analytic extension of degree $4$. We start by taking $\chi\in\mcal{C}^{\infty}_c(\R)$ such that $\ind{[-1/2,1/2]}\leq\chi\leq\ind{[-1,1]}$ and $h\in\mcal{C}_c^{5}(\R)$, this yields:
		$$\dfrac{h(\lambda)-h(\lambda')}{\lambda-\lambda'}=\dfrac{-1}{\pi}\iint_{\C}\dfrac{\bar{\partial}\big(h_4(z)\chi(y)\big)}{(\lambda-z)(\lambda'-z)}\diff^{2}z$$
		and thus by Fubini:
		$$|A_N^{(0)}(h)|\lesssim\iint_{\C}|\bar{\partial}\big(h_4(z)\chi(y)\big)|\cdot|(s_N-m_\mrm{eq})(z)|^{2}\diff^{2}z.$$
		We thus obtain:
		\begin{align*}
				\Ec{|A_N^{(0)}(\Xi^{-1}[g])|^{q}}&\lesssim\Ec{\iint\hdots\iint\prod_{j=1}^{q}|\bar{\partial}\big(h_4(z_j)\chi(y_j)\big)|\cdot|(s_N-m_\mrm{eq})(z_j)|^{2} \cdot\diff^{2}z_j}
				\\&\lesssim\iint\hdots\iint\Ec{\prod_{j=1}^{q}|(s_N-m_\mrm{eq})(z_j)|^{2} }\prod_{j=1}^{q}|\bar{\partial}\big(h_4(z_j)\chi(y_j)\big)|\cdot\diff^{2}z_j.
		\end{align*}
		Now for all $y_1,\dots,y_q\in(0,1)$ with $z_j=x_j+\ii y_j$, we have by Hölder inequality and Proposition \ref{prop:inital moment estimate}:
		\begin{align*}
			\Ec{\prod_{j=1}^{q}|(s_N-m_\mrm{eq})(z_j)|^{2} }&\leq\prod_{j=1}^{q}\|(s_N-m_\mrm{eq})(z_j)\|_{2q}^{2}
			\\&=\prod_{j=1}^{q}\|L_N^{(0)}(\dfrac{1}{\cdot-z_j})\|_{2q}^{2}
			\\&\lesssim\left(\dfrac{q}{N\alpha}+\dfrac{\log N}{\alpha}\right)^{q} \prod_{j=1}^{q}\left(\|\dfrac{1}{\cdot-z_j}\|_{\mcal{C}^{1}(\R)}+\|\dfrac{1}{\cdot-z_j}\|_{\msf{H}^{1/2}}\right)^{2}
			\\&\lesssim\left(\dfrac{q}{N\alpha}+\dfrac{\log N}{\alpha}\right)^{q} \prod_{j=1}^{q}\left(\dfrac{1}{y_j^{2}}+\dfrac{1}{y_j}\right)^{2}
			\\&\lesssim\left(\dfrac{q}{N\alpha}+\dfrac{\log N}{\alpha}\right)^{q} \prod_{j=1}^{q}\dfrac{1}{y_j^{4}}.
		\end{align*}
		A simple computation gives: $$\bar{\partial}\big(h_4(z)\chi(y)\big)=\bar{\partial}\big(\sum_{j=0}^4 \frac{1}{j!} (\ii y)^j h^{(j)}(x)\chi(y)\big)=\left( \dfrac{\chi(y)}{2}\frac{(\ii y)^4 }{4!}h^{(5)}(x)+\dfrac{\ii\chi'(y)}{2}h_4(z) \right)\underset{y\rightarrow0}{=}O(y^{4}).$$
		We thus obtain:
		\begin{align*}
			\|A_N^{(0)}(h)\|_q&\lesssim\left(\dfrac{q}{N\alpha}+\dfrac{\log N}{\alpha}\right)\iint\dfrac{1}{y^{4} }|\bar{\partial}\big(h_4(z)\chi(y)\big)|\diff^{2}z
			\\&\lesssim\left(\dfrac{q}{N\alpha}+\dfrac{\log N}{\alpha}\right)\int_\R\diff x\left(|h^{(5)}(x)|+\int_{1/2}^{1}\dfrac{\chi'(y)}{2}|h_4(z)|\diff y\right) 
			\\&\lesssim\left(\dfrac{q}{N\alpha}+\dfrac{\log N}{\alpha}\right)\|h \|_{\mcal{C}^{5}(\R)}.
		\end{align*}
		This leads to the following bound (noticing that $\Xi^{-1}[\mcal{C}_c^{6}(\R)]\subset\mcal{C}_c^{5}(\R)$):
		$$\|A_N^{(0)}(\Xi^{-1}[g])\|_{q}\lesssim\|g\|_{\mcal{C}^{6}(\R)} \left(\dfrac{q}{N\alpha}+\dfrac{\log N}{\alpha} \right )$$
		and thus for all $1\leq q\leq N\log N$:
				\begin{equation*}
			\| L_N^{(0)}(g)\|_q \lesssim\|g\|_{\mcal{C}^{\ell}(\R)}\left(\sqrt{\dfrac{q}{N\alpha}}+\dfrac{q}{N\alpha}+\dfrac{1}{\alpha}+\dfrac{\log N}{\alpha}\right)\lesssim\|g\|_{\mcal{C}^{\ell}(\R)}\left(\sqrt{\dfrac{q}{N\alpha}}+\dfrac{\log N}{\alpha}\right) .
		\end{equation*}
		We can now bootstrap this estimate, using that:
		$\|\dfrac{1}{\cdot-z}\|_{\mcal{C}^{\ell}(\R)}\lesssim\dfrac{1}{y^{\ell+1}},$
		and taking the pseudo-analytic extension of degree $\ell+1$, we obtain in the exact same way, for some $\ell'\geq1$:
		$$\|A_N^{(0)}(\Xi^{-1}[g])\|_{q}\lesssim\|g\|_{\mcal{C}^{\ell'}(\R)} \left(\dfrac{q}{N\alpha}+\dfrac{(\log N)^{2}}{\alpha^{2}} \right ).$$
		It leads to for all $q\leq N\log N$, using that $\alpha\gg(\log N)^{k}$ for all $k\geq1$:
			\begin{equation*}
			\| L_N^{(0)}(g)\|_q \lesssim\|g\|_{\mcal{C}^{\ell'}(\R)}\left(\sqrt{\dfrac{q}{N\alpha}}+\dfrac{1}{\alpha}\right) .
		\end{equation*}
		By Proposition \ref{prop:inital moment estimate}, when $q\geq N\log N$, it is easy to see that:
		$$\|L_N^{(0)}(g)\|_{q}\lesssim\|g\|_{\mcal{C}^{\ell'}(\R)}\left(\sqrt{\dfrac{q}{N\alpha}}+\sqrt{\dfrac{\log N}{\alpha}}\right)\lesssim\|g\|_{\mcal{C}^{\ell'}(\R)}\sqrt{\dfrac{q}{N\alpha}}.$$
		This concludes the proof for the case $k=0$.
		
		\textbf{2. Induction step.} Suppose for all $j\leq k-1$, there exists some $\ell_j\geq1$ such that:
		$$\|L_N^{(j)}(g)\|_{q}\lesssim \|g\|_{\mcal{C}^{\ell_j}(\R)}\left(\sqrt{\dfrac{q}{N\alpha}}+\dfrac{1}{\alpha^{j+1}}\right).$$We now use again \eqref{eq:eq boucle L} together with Proposition \ref{prop:rewriting1stloopeq}, to obtain with $g=\Xi[f']$:
		\begin{align}\label{eq:X_k+Z_k}
			\begin{aligned}
				\alpha^{k+1}\dfrac{\mc{L} [F]}{N\alpha}&=\alpha^{k+1}L_N^{(k)}(g)-\alpha^{k+1}\dfrac{1}{2}A_N^{(\lfloor\frac{k}{2} \rfloor)}(\Xi^{-1}[g])-\alpha^{k+1}(\sum_{i=1}^{\lfloor \frac{k}{2}\rfloor}\frac{\nu_i}{\alpha^{i}}\otimes L_N^{(k-i)})\left(\mcal{D} \Xi^{-1}[g]\right)
				\\&\quad- \alpha^{k} L_N^{(k-1)}(\Xi^{-1}[g]')+\alpha^{k+1}\dfrac{\bm{\mu}_N(\Xi^{-1}[g]')}{2N}
				\\&=X_k-\dfrac{\alpha^{k}}{N}Z_k,
			\end{aligned}
		\end{align}
		with $X_k=\alpha^{k+1}L_N^{(k)}(g)$ and
		$$Z_k=N\alpha\left(\dfrac{1}{2}A_N^{(\lfloor\frac{k}{2} \rfloor)}(\Xi^{-1}[g])+(\sum_{i=1}^{\lfloor \frac{k}{2}\rfloor}\frac{\nu_i}{\alpha^{i}}\otimes L_N^{(k-i)})\left(\mcal{D} \Xi^{-1}[g]\right)+ \dfrac{1}{\alpha} L_N^{(k-1)}(\Xi^{-1}[g]')	-\dfrac{\bm{\mu}_N(\Xi^{-1}[g]')}{2N}\right).$$
		We now use again
		Proposition \ref{prop:bound_normq} with $\lambda=N/\alpha^{k}$, to deduce that:
		\begin{equation}\label{eq:stein inequality k}
			\begin{aligned}
				\|X_k\|_q &\le \sqrt{\dfrac{q\alpha^{k} \|\Gamma_k \|_{q/2}}{N}} +\dfrac{ \alpha^{k}\|Z_k\|_q }{N}
				\\&\lesssim\alpha^{k+1}\Bigg(\sqrt{\dfrac{q\|g\Xi^{-1}[g]'\|_{\infty}}{N\alpha}}+\|A_N^{(\lfloor\frac{k}{2} \rfloor)}(\Xi^{-1}[g])\|_q+\sum_{i=1}^{\lfloor \frac{k}{2}\rfloor}\|\frac{\nu_i}{\alpha^{i}}\otimes L_N^{(k-i)}\left(\mcal{D} \Xi^{-1}[g]\right)\|_q
				\\&\quad+ \dfrac{1}{\alpha}\|L_N^{(k-1)}(\Xi^{-1}[g]')\|_q+\dfrac{1}{N}\|\Xi^{-1}[g]'\|_\infty\Bigg),
			\end{aligned}
		\end{equation}
		where we noticed that
		$$\Gamma_k=\alpha^{k}\Gamma_0= \alpha^{k}\nabla F \cdot\nabla X_0=\alpha^{k+1}\bm{\mu}_N(g\Xi^{-1}[g]').$$
		By induction hypothesis, we have the following estimates for some $\ell\geq1$:
		$$\sum_{i=1}^{\lfloor \frac{k}{2}\rfloor}\|\frac{\nu_i}{\alpha^{i}}\otimes L_N^{(k-i)}\left(\mcal{D} \Xi^{-1}[g]\right)\|_q\lesssim\|g\|_{\mcal{C}^{\ell}(\R)} \left(\dfrac{1}{\alpha}\sqrt{\dfrac{q}{N\alpha}}+\dfrac{1}{\alpha^{k+1}}\right) $$
		and
		$$\dfrac{1}{\alpha}\|L_N^{(k-1)}(\Xi^{-1}[g]')\|_q\lesssim\|g\|_{\mcal{C}^{\ell}(\R)}\left(\dfrac{1}{\alpha}\sqrt{\dfrac{q}{N\alpha}}+\dfrac{1}{\alpha^{k+1}}\right).$$
		Finally, to control the anisotropy, we use the same Helffer–Sjöstrand formula as before, choosing a pseudo analytic extension of degree $m$ sufficiently large (depending only on $k$), we obtain just as before with $h=\Xi^{-1}[g]$:
			\begin{align*}
			\|A_N^{(\lfloor\frac{k}{2}\rfloor)}(h)\|_q&\lesssim\left(\dfrac{q}{N\alpha}+\dfrac{1}{\alpha^{2\lfloor\frac{k}{2}\rfloor+2}}\right)\int_\R\diff x\left(|h^{(m+1)}(x)|+\int_{1/2}^{1}\dfrac{\chi'(y)}{2}|h_m(z)|\diff y\right) 
			\\&\lesssim\left(\dfrac{q}{N\alpha}+\dfrac{1}{\alpha^{k+1}}\right)\|h\|_{\mcal{C}^{m+1}(\R)}\\&\lesssim\left(\dfrac{q}{N\alpha}+\dfrac{1}{\alpha^{k+1}}\right)\|g\|_{\mcal{C}^{\ell}(\R)}
		\end{align*}
		where we used that $2\lfloor\frac{k}{2}\rfloor+2\geq k+1$ and Proposition \ref{prop:inverse Xi}.

		Using \eqref{eq:stein inequality k}, and collecting all the bounds, it leads to, when $q\leq N\log N$:
		\begin{equation*}
			\|L_N^{(k)}(g)\|_q \lesssim\|g\|_{\mcal{C}^{\ell}(\R)}\left(\sqrt{\dfrac{q}{N\alpha}}+\dfrac{1}{\alpha^{k+1}}\right).
		\end{equation*}
		Furthermore, when $q\geq N\log N$, from Proposition \ref{prop:inital moment estimate}, we have:
	\begin{align*}
		\|L_N^{(k)}(g)\|_q\lesssim\|L_N^{(0)}(g)\|_q+\dfrac{\|g\|_{\mcal{C}^{\ell}(\R)}}{\alpha}&\lesssim(\|g\|_{\msf{H}^{1/2}}+\|g\|_{\mcal{C}^{\ell}(\R)})\left(\sqrt{\dfrac{q}{N\alpha}}+\dfrac{1}{\alpha}\right)
		\\&\lesssim\|g\|_{\mcal{C}^{\ell}(\R)}\sqrt{\dfrac{q}{N\alpha}}.
	\end{align*}
		This concludes the induction and the proof.
	\end{Pro}

	We now give optimal moment estimates for linear statistics with test-functions supported near infinity.

		\begin{Prop}[Improved moment estimates at infinity]\label{prop: moment estimates infinity}
			Let $\beta\gg N^{-1+\varepsilon}$ for some $\varepsilon>0$, there exists $\ell>0$ (depending only on $\varepsilon$) such that for all $g\in\mcal{C}^\ell(\R)$, under the same assumptions as in \textit{(ii)} in Lemma \ref{lem:initial concentration}, we have for all $q\geq1$:
		$$\|\bm{\mu}_N(g(1-\chi))\|_{q}\lesssim\|g\|_{\mcal{C}^{\ell}(\R)} \sqrt{\dfrac{q}{N\alpha}}.$$
	\end{Prop}
	
	\begin{Pro}\textbf{1. Preliminaries.} We proceed in a similar way as in the proof of Proposition \ref{prop: moment estimates}. We take $f(x)=O(1/x)$ at infinity and $\chi$ as in Lemma \ref{lem:initial concentration} \textit{(ii)} and set $g=\Xi[(f(1-\chi))']$ and $F=N\bm{\mu}_N(f(1-\chi))$ in \eqref{eq:def L}. We obtain for all $k\geq0$:
\begin{align*}
	\dfrac{\mc{L} [F]}{N\alpha}&=-\left(\dfrac{1}{\alpha}-\dfrac{1}{2N}\right) \bm{\mu}_N((f(1-\chi))'')+\dfrac{1}{2}\bm{\mu}_N(V'(f(1-\chi))')-\dfrac{1}{2}\bm{\mu}_N^{\otimes2}(\mcal{D}[(f(1-\chi))'])
	\\&=L_N^{(2k)}(\Xi[(f(1-\chi))'])-\dfrac{1}{2}A_N^{(k)}((f(1-\chi))')-(\sum_{i=1}^{k}\frac{\nu_i}{\alpha^{i}}\otimes L_N^{(2k-i)})\left(\mcal{D} [(f(1-\chi))']\right)
	\\&\quad- \dfrac{1}{\alpha} L_N^{(2k-1)}((f(1-\chi))'')+\dfrac{\bm{\mu}_N((f(1-\chi))'')}{2N}
	\\&=X-\dfrac{Z}{N\alpha}
\end{align*}
where by Proposition \ref{prop:inverse Xi}, $X\defi L_N^{(2k)}(g)=L_N^{(2k)}(g(1-\chi))=\bm{\mu}_N(g(1-\chi))$ and:
\begin{align*}
	Z&\defi \dfrac{N\alpha}{2}A_N^{(k)}(\Xi^{-1}[g](1-\chi))+N\alpha(\sum_{i=1}^{k}\frac{\nu_i}{\alpha^{i}}\otimes L_N^{(2k-i)})\left(\mcal{D}[\Xi^{-1}[g](1-\chi)]\right)
	\\&\quad+NL_N^{(k-1)}((\Xi^{-1}[g](1-\chi))')-\dfrac{\alpha}{2}\bm{\mu}_N((\Xi^{-1}[g](1-\chi))').
\end{align*}
\textbf{2. Estimation of $\|Z\|_q$.} Since $f(x)=O(1/x)$ at infinity, by Proposition \ref{prop:inverse Xi}, $g(x)\sim V'(x)f(x)\sim xf(x)=O(1)$ at infinity, thus $$\Xi^{-1}[g(1-\chi)](x)=\dfrac{g(x)}{V(x)}(1+o(1))=O\left (\dfrac{1}{x}\right ),$$
where the above estimate is differentiable.
Since $ L_N^{(2k-1)}((\Xi^{-1}[g](1-\chi))')= \bm{\mu}_N((\Xi^{-1}[g](1-\chi))')$ and $g=O(1)$, we have for all $q\geq1$ by Proposition \ref{prop:inital moment estimate}:
$$\Big\|\left(\dfrac{1}{\alpha}-\dfrac{1}{2N}\right) \bm{\mu}_N((\Xi^{-1}[g](1-\chi))')\Big\|_q\lesssim\|g\|_{\mcal{C}^{2}(\R)}\dfrac{1}{\alpha}\left(\sqrt{\dfrac{q}{N\alpha}}+\sqrt{\dfrac{\log N}{\alpha}}\right ).$$
Furthermore,
$$\sum_{i=1}^{k}\dfrac{1}{\alpha^{i}}\|(\nu_i\otimes L_N^{(2k-i)})\left(\mcal{D}[\Xi^{-1}[g](1-\chi)]\right)\|_q=\sum_{i=1}^{k}\dfrac{1}{\alpha^{i}}\Big\|\int_{[-M,M]^{c}}\nu_i\left(\mcal{D}[\Xi^{-1}[g](1-\chi)](\cdot,y)\right)\diff\bm{\mu}_N(y)\Big\|_q,$$
and since for all $|y|>M$,  $\mcal{D}[\Xi^{-1}[g](1-\chi)](\cdot,y)$ is smooth and bounded so is the map $y\allowbreak\mapsto\nu_i\left(\mcal{D}[\Xi^{-1}[g](1-\chi)](\cdot,y)\right)$, thus by Proposition \ref{prop:inital moment estimate} the continuity of $\nu_i$, $\mcal{D}$ and $\Xi^{-1}$, we have for some $r\geq0$:
\begin{multline*}
	\sum_{i=1}^{k}\dfrac{1}{\alpha^{i}}\|(\nu_i\otimes L_N^{(2k-i)})\left(\mcal{D}[\Xi^{-1}[g](1-\chi)]\right)\|_q\\\lesssim\dfrac{1}{\alpha}\left(\sqrt{\dfrac{q}{N\alpha}}+\sqrt{\dfrac{\log N}{\alpha}}\right )\max_{i, |y|>M}\nu_i\left(\mcal{D}[\Xi^{-1}[g](1-\chi)](\cdot,y)\right)
	\\\lesssim\dfrac{1}{\alpha}\left(\sqrt{\dfrac{q}{N\alpha}}+\sqrt{\dfrac{\log N}{\alpha}}\right )\|g\|_{\mcal{C}^{r}(\R)}.
\end{multline*}
Finally, for the anisotropy, noticing that $\mcal{D}[h(1-\chi)](x,y)=0$ for $|x|,|y|\leq M$, we have:
\begin{align*}
	A_N^{(k)}(\Xi^{-1}[g](1-\chi))&=\int_{([-M,M]^{c})^{2}}\mc{D}[\Xi^{-1}[g](1-\chi)](x,y)\diff\bm{\mu}_N(x)\diff\bm{\mu}_N(y)
	\\&\quad+\int_{[-M,M]^{c}}L_N^{(k)}\left(\mc{D}[\Xi^{-1}[g](1-\chi)](\cdot,y)\right)\diff\bm{\mu}_N(y).
\end{align*}
For the first term, we simply use that $\mc{D}[\Xi^{-1}[g](1-\chi)](\cdot,\cdot)$ is bounded by Proposition \ref{prop:inverse Xi}, thus by Proposition \ref{prop:inital moment estimate}:
\begin{align*}
	\Big\|\int_{([-M,M]^{c})^{2}}\mc{D}[\Xi^{-1}[g](1-\chi)](x,y)\diff\bm{\mu}_N(x)\diff\bm{\mu}_N(y)\Big\|_q&\lesssim\|g\|_{\mcal{C}^{2}(\R)}\cdot\|\bm{\mu}_N(1-\chi)\|_{2q}^{2}
	\\&\lesssim\|g\|_{\mcal{C}^{2}(\R)}\left(\dfrac{q}{N\alpha}+\dfrac{\log N}{\alpha}\right).
\end{align*}
For the last term, we use:
\begin{multline*}
	\Big\|\int_{[-M,M]^{c}}L_N^{(k)}\left(\mc{D}[\Xi^{-1}[g](1-\chi)](\cdot,y)\right)\diff\bm{\mu}_N(y)\Big\|_q
	\\\lesssim\Big\|\max_{|z|>M}|L_N^{(k)}\left(\mc{D}[\Xi^{-1}[g](1-\chi)](\cdot,z)\right)|\cdot\bm{\mu}_N(1-\chi)\Big\|_q
	\\\lesssim\Big\|\max_{|z|>M}|L_N^{(k)}\left(\mc{D}[\Xi^{-1}[g](1-\chi)](\cdot,z)\right)\Big\|_{2q}\Big\|\bm{\mu}_N(1-\chi)\Big\|_{2q}.
\end{multline*}
We now use the fact that for all $|z|>M$, $x\mapsto\mc{D}[\Xi^{-1}[g](1-\chi)](x,z)\in\mcal{W}^{1,2}(\R)$. Indeed, it is $\mcal{C}^{1}_{\mrm{loc}}(\R)$ and both it and its derivative behave like $O\left(1/x^{2}\right)$ at infinity. Thus, we use the fact that, setting $Y(z)\defi L_N^{(k)}\left(\mc{D}[\Xi^{-1}[g](1-\chi)](\cdot,z)\right)$ and using that $Y\in\mcal{W}^{1,2}(\R)$ and $Y(z)\rightarrow0$ as $z\rightarrow\infty$:
$$|Y(z)|^{2}\leq\int_{-\infty}^{z}2Y(t)Y'(t)\diff t\leq\|Y\|_{L^2(\R)}\cdot\|Y'\|_{L^2(\R)},$$
thus the following inequality holds:
$$\max_{|z|>M}|Y(z)|\leq\sqrt{\|Y\|_{L^2(\R)}\cdot\|Y'\|_{L^2(\R)}}.$$
We then obtain by Cauchy-Schwarz, Fubini, Hölder inequality and Proposition \ref{prop: moment estimates}:
\begin{align*}
	\Big\|\max_{|z|>M}|L_N^{(k)}\left(\mc{D}[\Xi^{-1}[g](1-\chi)](\cdot,z)\right)\Big\|_{2q}&\lesssim\Big\|\|Y\|_{L^2(\R)}\Big\|_{2q}^{1/2}\cdot\Big\|\|Y'\|_{L^2(\R)}\Big\|_{2q}^{1/2}
	\\&\lesssim\Big\|\|Y(\cdot)\|_{2q}\Big\|_{L^2(\R)}^{1/2}\cdot\Big\|\|Y'(\cdot)\|_{2q}\Big\|_{L^2(\R)}^{1/2}.
\end{align*}
Now, decomposing the function into a compactly supported part and a part supported near infinity and using Proposition \ref{prop: moment estimates} and Proposition \ref{prop:inital moment estimate} for the part near infinity, we get for $i\in\{0,1\}$:
\begin{align*}
	\int_{\R}\|\partial^{i}_tY(t)\|_{2q}^{2}\diff t&\lesssim\int_{\R}\|L_N^{(k)}\left(\chi(\cdot)\partial^{i}_t\mc{D}[\Xi^{-1}[g](1-\chi)](\cdot,t)\right)\|_{2q}^{2}\diff t
	\\&\quad+\int_{\R}\|\bm{\mu}_N\left((1-\chi)(\cdot)\partial^{i}_t\mc{D}[\Xi^{-1}[g](1-\chi)](\cdot,t)\right)\|_{2q}^{2}\diff t
		\\&\lesssim\int_\R\|\partial^{i}_t\mc{D}[\Xi^{-1}[g](1-\chi)](\cdot,t)\|_{\mcal{C}^{r}(\R)}^{2}\diff t\cdot\left(\sqrt{\dfrac{q}{N\alpha}}+\dfrac{1}{\alpha^{k+1}}+\sqrt{\dfrac{\log N}{\alpha}}\right) 
\end{align*}
For all $i\in\{0,1\}$, the map $t\mapsto\|\partial_z^{i}\mc{D}[\Xi^{-1}[g](1-\chi)](\cdot,t)\|_{\mcal{C}^{r}(\R)}$ is bounded on compact sets and at infinity behaves at least like $O(1/t)$ since $\Xi^{-1}[g]$ and its derivatives vanish like $O(1/t)$ at infinity, thus the integral is finite and we have:
$$\int_\R\|\partial^{i}_t\mc{D}[\Xi^{-1}[g](1-\chi)](\cdot,t)\|_{\mcal{C}^{r}(\R)}^{2}\diff t\lesssim\|g\|_{\mcal{C}^{r'}(\R)}\int_{\R}\dfrac{\diff x}{x^{2}+1}\lesssim\|g\|_{\mcal{C}^{r'}(\R)}$$
and finally we obtain:
$$\Big\|\int_{[-M,M]^{c}}L_N^{(k)}\left(\mc{D}[\Xi^{-1}[g](1-\chi)](\cdot,y)\right)\diff\bm{\mu}_N(y)\Big\|_q\lesssim\|g\|_{\mcal{C}^{r'}(\R)}\left(\dfrac{q}{N\alpha}+\dfrac{\log N}{\alpha}\right) .$$
We thus, obtain that for all $q\leq N\log N$:
$$\|Z\|_q\lesssim N\alpha\|g\|_{\mcal{C}^{r'}(\R)}\left(\sqrt{\dfrac{q}{N\alpha}}+\dfrac{\log N}{\alpha}\right).$$
Finally, there exists $h\in\mcal{C}^{0}(\R)$ with $\|h\|_{\mcal{C}^{0}(\R)}\lesssim\|g\|_{\mcal{C}^{r}(\R)}^{2}$ such that:
$$\Gamma= \nabla F \cdot\nabla X=\dfrac{1}{N}\sum_{i=1}^{N}\left(f(1-\chi)\right)' (\lambda_i)\left(g(1-\chi)\right) '(\lambda_i)=\bm{\mu}_N(h(1-\chi)),$$
thus $\|\Gamma\|_{q/2}\lesssim\|g\|_{\mcal{C}^{r}(\R)}^{2}\left(\sqrt{\tfrac{q}{N\alpha}}+\sqrt{\tfrac{\log N}{\alpha}}\right).$

\textbf{3. Conclusion.} We can now use Proposition \ref{prop:bound_normq} to obtain that there exists $\ell\geq1$, such that for all $q\leq N\log N$:
\begin{equation*}
	\|X\|_q \lesssim \sqrt{\dfrac{q \|\Gamma \|_{q/2}}{N\alpha}}+\dfrac{\|Z\|_q}{N\alpha}\lesssim\|g\|_{\mcal{C}^{\ell}(\R)}\left(\sqrt{\dfrac{q}{N\alpha}}+\dfrac{\log N}{\alpha}\right) .
\end{equation*}
We can now bootstrap this inequality to obtain for all $k\geq0$ and $q\leq N\log N$:
\begin{equation*}
	\|X\|_q \lesssim \sqrt{\dfrac{q \|\Gamma \|_{q/2}}{N\alpha}}+\dfrac{\|Z\|_q}{N\alpha}\lesssim\|g\|_{\mcal{C}^{\ell'}(\R)}\left(\sqrt{\dfrac{q}{N\alpha}}+\dfrac{1 }{\alpha^{k+1}}\right) ,
\end{equation*}
and thus remove the second term in the RHS, taking $k$ large enough. Finally by Proposition \ref{prop:inital moment estimate} for all $q\geq N\log N$:
$$\|X\|_q\lesssim\|g\|_{\mcal{C}^{0}(\R)}\sqrt{\dfrac{q}{N\alpha}}\lesssim\|g\|_{\mcal{C}^{\ell'}(\R)}\sqrt{\dfrac{q}{N\alpha}}.$$
This concludes the proof.
	\end{Pro}

	A consequence of the optimal moment estimates obtained in Proposition \ref{prop: moment estimates} and Proposition \ref{prop: moment estimates infinity} is the global CLT. 
The following lemma is used to recast the asymptotic variance in Theorem \ref{thm:main CLT global} and is proved in \cite[Eq (4.15)]{LambertLedouxWebb}.

\begin{Lem}[Variance formula]\label{lem:varformula}
	The following identity holds for all $g\in\mcal{C}^{1}(\R)$:
	$$\int_{-1}^1\Xi^{-1}[g]'(x)g(x)\diff\mu_{\mrm{eq}}(x)=
\sigma^2(g)$$
where $\sigma^2$ was defined in \eqref{eq: norm H 1/2 global} and $\Xi^{-1}[g]$ given in \eqref{eq:inverse Xi formula}. 
\end{Lem}
	We can now prove Theorem \ref{thm:main CLT global}.
	\begin{Pro}[of Theorem \ref{thm:main CLT global}] Let $f$ such that $f'=\Xi^{-1}[g]$, we take again the decomposition of \eqref{eq:X_k+Z_k} with $F=N\bm{\mu}_N(f)$:
		$$	\widetilde{X}_k=\sqrt{N\alpha}L_N^{(k)}(g)=\dfrac{\sqrt{N\alpha}}{\alpha^{k+1}}X_k=\dfrac{\mc{L} [F]+Z_k}{\sqrt{N\alpha}},$$
		Using Proposition \ref{prop:bound Wp}, we obtain with $\lambda=\sqrt{N\alpha}$:
			$$\mathbf{W}_q\left(\widetilde{X}_k,\mcal{N}\left(0,\sigma^2(g)\right)\right)\leq \dfrac{1}{\sqrt{N\alpha}}\big(\dfrac{\sqrt{q}}{\sigma(g)}\| W_k\|_{q}+\|Z_k\|_{q} \big), \qquad  W_k = \nabla \widetilde{X}_k \cdot \nabla F - \sqrt{N\alpha}\sigma^{2}(g). $$
			We obtain easily using Lemma \ref{lem:varformula}:
			$$W_k =\sqrt{N\alpha}\left(\dfrac{1}{N}\sum_{i=1}^{N}g'(\lambda_i)\Xi^{-1}[g](\lambda_i)-\int_{-1}^1\Xi^{-1}[g]'(x)g(x)\diff\mu_{\mrm{eq}}(x)\right)=\sqrt{N\alpha}L_N^{(0)}(\Xi^{-1}[g]'g),$$
			and from Proposition \ref{prop: moment estimates}, Proposition \ref{prop:inverse Xi} and \ref{prop: moment estimates infinity} and using a cutoff $\chi$:, 
			\begin{align*}
				\dfrac{	\|W_k\|_q}{\sqrt{N\alpha}}=\|L_N^{(0)}(\Xi^{-1}[g]'g)\|_q &\lesssim\|\Xi^{-1}[g]'g\|_{\mcal{C}^{\ell}(\R)}\left(\sqrt{\dfrac{q}{N\alpha}}+\dfrac{1}{\alpha}\right) \\&\lesssim\|g\|_{\mcal{C}^{\ell}(\R)}\|\Xi^{-1}[g]'\|_{\mcal{C}^{\ell}(\R)}\left(\sqrt{\dfrac{q}{N\alpha}}+\dfrac{1}{\alpha}\right).
				\\&\lesssim\|g\|_{\mcal{C}^{\ell'}(\R)}^{2}\left(\sqrt{\dfrac{q}{N\alpha}}+\dfrac{1}{\alpha}\right).
			\end{align*}
			Using the exact same bounds as in the proof of Proposition \ref{prop: moment estimates}, we obtain for all $q\geq1$ and for some $r\geq1$:
			\begin{align*}
				\dfrac{1}{\sqrt{N\alpha}}\|Z_k\|_q&\lesssim\|g\|_{\mcal{C}^{\ell}(\R)}\sqrt{N\alpha} \left(\dfrac{1}{\alpha}\sqrt{\dfrac{q}{N\alpha}}+\dfrac{q}{N\alpha}+\dfrac{1}{\alpha^{k+1}}+\dfrac{1}{N}\right)\\&\lesssim\|g\|_{\mcal{C}^{\ell}(\R)} \left(\dfrac{\sqrt{q}}{\alpha}+\dfrac{q}{\sqrt{N\alpha}}+\dfrac{\sqrt{N\alpha}}{\alpha^{k+1}}\right).
			\end{align*}
			We thus obtain for some $r'\geq1$ and all $q\geq1$:
			$$\mathbf{W}_q\left(\sqrt{N\alpha}L_N^{(k)}(g),\mcal{N}\left(0,\sigma^2(g)\right)\right)\lesssim C(g) \left (\dfrac{q}{\sqrt{N\alpha}}+\dfrac{\sqrt{q}}{\alpha}+\dfrac{\sqrt{N\alpha}}{\alpha^{k+1}}\right )$$
			which concludes the proof.
	\end{Pro}
	
	From Proposition \ref{prop:asympcorrelators}, we deduce the expansion of the free energy Theorem \ref{thm:free energy}. The first proof of these expansions are due to \cite{BoG1} in the $\beta$-fixed case.

	\begin{Pro}[of Theorem \ref{thm:free energy}] First notice that the equilibrium measure associated with the potential $V(x)=2x^{2}$ is supported on $[-1,1]$, see \cite[Section 7.2.2]{BoG2}. Using differentation under the integral sign and setting $V_t(x)=tV(x)+(1-t)2x^{2}$, we obtain:
		\begin{align*}
			-\dfrac{2}{N\alpha}\log\dfrac{\mcal{Z}_N[V]}{\mcal{Z}_N[2x^{2}]}&=\int_0^1\E_N^{V_t}\left[\int_\R (V(x)-2x^2)\diff \bm{\mu}_N(x)\right]\diff t
			\\&=\int_0^1\braket{V-2x^2}^{V_t}\diff t+\int_0^1\int_{-1}^{1}(V(x)-2x^2) \diff\mu_{\mrm{eq}}^{V_t}(x)\diff t.
		\end{align*}
		Above, we have set $\mu_{\mrm{eq}}^{V_t}=t\mu_{\mrm{eq}}+(1-t)\mu_{\mrm{eq}}^{V_0}$ for the equilibrium associated with the potential $V_t$ and \textit{idem} for $\E_N^{V_t}$ and $\braket{\cdot}^{V_t}$. The integrand in the first term of the RHS admits an expansion at all order. Showing that the each coefficients is $t$-integrable relies on \cite[Lemma 5.1]{BoG1} and the fact that the remainder is uniform in $t$ comes from the constant in the control of $\Xi^{-1}$. The fact that $\lim_{N\rightarrow\infty}\dfrac{1}{N\alpha}\log\mcal{Z}_N\left[V\right]=-\dfrac{\mcal{E}(\mu_{\mrm{eq}})}{2}$ comes from a large deviation principle argument, see \cite[Proof of Theorem 2.6.1]{anderson2010introduction}. Finally, the term proportional to $\dfrac{1}{\alpha}$ of $-\dfrac{2}{N\alpha}\log\dfrac{\mcal{Z}_N[V]}{\mcal{Z}_N[2x^{2}]}$ is, by Proposition \ref{prop:asympcorrelators}: $$\int_0^1\nu_1(V-2x^2)^{V_t}\diff t=\int_0^{1}\diff t\int_{-1}^{1}\partial_1\Xi_t^{-1}[V-2\cdot^{2}](x)\diff\mu_{\mrm{eq}}(x)$$
		with $\nu_1(\cdot)^{V_t}$ and $\Xi_{t}$ denotes respectively the linear functional $\nu_1$ and the master operator under the choice of potential $V_t$. By \cite[Proof of Corollary 1.3]{dworaczekguera2025clt},
		\begin{align*}
			\int_0^{1}\diff t\int_{-1}^{1}\partial_1\Xi_t^{-1}[V-2\cdot^{2}](x)\diff\mu_{\mrm{eq}}(x)=2\int_0^{1}\partial_t\mrm{Ent}\left[\mu_{\mrm{eq}}^{V_t}\right]\diff t
			&= 2\left(\mrm{Ent}\left[\mu_{\mrm{eq}}\right]-\mrm{Ent}\left[\mu_{\mrm{eq}}^{V_0}\right]\right)
			\\&=2\mrm{Ent}\left[\mu_{\mrm{eq}}\right]-2\log\pi+1.
		\end{align*}
		Above we have used that $\mrm{Ent}\left[\mu_{\mrm{eq}}^{V_0}\right]=\log\pi-\frac{1}{2}$ by \cite{dworaczekguera2025clt}. This concludes the proof. \qedsymbol{}
	\end{Pro}
	
	\begin{Rem}\label{rem:asymp gaussian partition}We recall that the Gaussian partition is given, see \cite[Section 7.2.2]{BoG2}, by:
		$$\mcal{Z}_N[2x^{2}]=\dfrac{N!(2\pi)^{N}\left(\frac{N\alpha}{2}\right)^{\frac{\alpha}{2}-N} }{\Gamma\left(\dfrac{\alpha}{2N}\right)^{N}\Gamma_2\left(N+1,\dfrac{2N}{\alpha},1\right)  }\cdot(4N)^{-\frac{N\alpha}{4}-\frac{N}{2}+\frac{\alpha}{2}},$$
		where $\Gamma_2$ is the so-called Barnes double Gamma function, see \cite{alexanian2023barnes}.
	\end{Rem}

	\section{Local laws}\label{sec: local laws}

In this section, we derive local laws Theorem \ref{thm:local law} which is crucial to obtain the mesoscopic concentration estimates in Section \ref{sec: concentration} and thus conclude on the mesoscopic CLT in the bulk Theorem \ref{thm:mesoc clt bulk} via Stein's method in Section \ref{sec: proof CLT}.

In Subsection \ref{subsec3:loopeq}, we define a sequence of holomorphic function $(\xi_k(z))_{k\geq0}$ which is very close to $(m_k)_{k\geq0}$, the Stieltjes transforms of $\mu_k$'s defined in Defintion \ref{def:mup}. In Subsection \ref{subsec3: a priori BMP}, we adapt the proof of the local law  proved in \cite{bourgade2022optimal} to the intermediary temperature regime to obtain an \textit{a priori estimate} for $\eta\gtrsim\alpha^{-1}$, \textit{i.e.} in the so-called random matrix regime. In Subsection \ref{subsec3:RMT regime}, we improve this local law inductively and obtain the optimal one. In Subsection \ref{subsec3:poisson}, we obtain the optimal local law in the so-called Poisson regime, namely for $\eta\lesssim\alpha^{-1}$.

\subsection{Loop equations \& properties of the equilibrium measure}\label{subsec3:loopeq}
The goal of this subsection is to provide some background on loop equations for $\beta$-ensembles and the associated properties of the equilibrium measure.

We first start with the definition of the random variable $\Delta_k$ which is an important object in this version of the loop equations.

\begin{Def}\label{def: delta}
	We define for all $z\in\mcal{U}$ and $k\geq0$
	$$\Delta_k(z)\defi \left(\bm{\mu}_N-\mu_k\right) (W_z),\hspace{1.5cm}\text{with}\quad W_z\defi\dfrac{V'(\cdot)-V'(z)}{\cdot-z}.$$
\end{Def}

In the Gaussian case $V(x)=x^{2}$, $W_z=1$ so that  $\Delta_k=0$. In general, controlling the moments of $\Delta_k(z)$ for $z\in\mcal{U}$ is a crucial step to deduce a local law by the scheme developed in \cite{bourgade2022optimal}.

	From Proposition \ref{prop: moment estimates}, we deduce the following estimate on $L_N^{(k)}(W_z)$.

\begin{Cor}\label{cor: glocon}
	There exists $c>0$ such that for all $k\geq1$ and $z\in\mcal{U}$, for all $q\geq1$ and $M$ large enough, we have:
	$$\|L_N^{(k)}(W_z)\|_{q}\lesssim\left(\dfrac{1}{\alpha^{k+1}}+\sqrt{\dfrac{q}{N\alpha}}\right). $$
\end{Cor}

\begin{Pro} Let $\varepsilon>0$ and $\phi\in\mcal{C}_c^{\infty}(\R)$ such that $\ind{[-M,M]}\leq\phi\leq\ind{[-M-1,M+1]}$ decomposing $W_z=W_z\phi+(1-\phi)W_z$ and using that supp $\mu_{\mrm{eq}}=[-1,1]$, we have by Proposition \ref{prop: moment estimates} and \ref{prop: moment estimates infinity} using that $W_z\in\mcal{C}^{\infty}(\R)$ and hence bounded, we have:
\begin{align*}
	\|L_N^{(k)}(W_z)\|_{q}&\leq\|L_N^{(k)}(\phi W_z)\|_{q}+\|L_N^{(k)}\left((1-\phi)W_z\right) \|_{q}
	\\&\lesssim \|W_z\phi\|_{\mcal{C}^{\ell}(\R)}\left(\dfrac{1}{\alpha^{k+1}}+\sqrt{\dfrac{q}{N\alpha}}\right)+\|\bm{\mu}_N\left((1-\phi)W_z\right) \|_{q}
	\\&\lesssim \left(\dfrac{1}{\alpha^{k+1}}+\sqrt{\dfrac{q}{N\alpha}}\right).
\end{align*}
\end{Pro}

We recall the following identity, see \cite[Eq. (2.9)]{dworaczekguera2025clt} for a simple proof.
\begin{Lem}[Loop equation]\label{lem:Leq}
	For a random variable $G(z)\in\mathcal{C}^1(\R^N)$, one has for  $z\in \mathcal{U}_+$,
	\[
	\E_N\left[\left( s_N^2(z) + V'(z) s_N(z)  + \bm{\mu}_N(W_z) + \left( \dfrac{2}{\alpha}-\dfrac{1}{N}\right) \partial_z s_N(z) \right)G(z)\right] 
	= \frac1{\alpha N} 	\E_N\bigg[ \sum_{j=1}^{N} \frac{\partial_{\lambda_j}G(z)}{\lambda_j-z}  \bigg] . 
	\]
\end{Lem}

We recall the following classical result of Stieltjes transform $m_0(z)$ of the equilibrium measure $\mu_\mrm{eq}$. We define the following quadratic polynomial in $s(z)$ for which $m_0(z)$ is a solution, for all $z\in\mcal{U}$, 
\[
P_0(s) \defi s^2 + V'(z) s + \mu_{\mrm{eq}}(W_z).
\]

We recall that the holomorphic map $S$ on $\mcal{U}$ is defined in Assumptions \ref{assumptions}

\begin{Lem} \label{lem:eq}
	
	There exists $\eta_0>0$ such that for $\Omega \defi B+\ii(0,\eta_0)$ where $\bar{B} \subset (-1,1)$, there exists $c>0$ such that the discriminant $ \triangle_0(z)$ of the polynomial $P_0(\cdot)$ satisfies for all $z\in \Omega$:
	$$| \triangle_0(z)| \ge c.$$
	Moreover, for all $z\in\mcal{U}$ the Stieltjes transform $\xi_0(z)$ of the equilibrium measure is the root  of $P_0(\cdot)$ and satisfies:
	$$2m_0(z)=2\xi_0(z) = -V'(z)+2S(z)\sqrt{z^{2}-1}$$ where $\sqrt{\cdot}$ is the principal branch of the square root.
	In particular, there exists $C>0$ such that $$\|m_0\|_{\mathcal{C}^1(\Omega)} \le C.$$
\end{Lem}

\begin{Pro}
	By the law of large numbers, for a fixed $z\in\Omega$,  $s_N(z) \to m_0(z)$ and $\partial_z s_N(z) \to \partial_z m_0(z)$ with overwhelming probability as $N\rightarrow\infty$. Thus taking the limit as $N\to\infty$ in \eqref{Leq} with $G=1$, we obtain $P_0(m_0) =0$ for all $z\in\Omega$. 
	
	Recall that for all $E\in\R$, the density of the equilibrium measure satisfies:
	\begin{equation} \label{PS}
		\mu_{\mrm{eq}}(E) = \dfrac1\pi\lim_{\eta\to0}\im{ m_0(E+\ii \eta) }.
	\end{equation}
	The two roots of the polynomial $P_0$ are:
	\begin{equation}\label{roots}
		2m_0(z)=2\xi_0(z)=-V'(z)+2S(z)\sqrt{z^{2}-1}, \qquad  2\widetilde{\xi}_0 =  -V'(z)-2S(z)\sqrt{z^{2}-1}.
	\end{equation}
Also by \cite[Eq (B.7)]{bourgade2022optimal}, there exists $\eta_0>0$ and $c>0$ such that for all $z\in\Omega$: $$|\triangle(z)|^{1/2}=2|S(z)|\cdot|z^{2}-1|^{1/2}\geq c.$$
	Finally, the bound $\|m_0\|_{\mathcal{C}^1(\Omega)} \le C$ for $z\in \Omega$ follows directly from the expression  \eqref{roots}.
\end{Pro}


\begin{Def}\label{def:Pk}
	We define for $k\geq1$ and $z\in\Omega$:
	\[
	P_k : s \in \C\mapsto s^2 + V'(z) s + \mu_k(W_z)  + \tfrac2\alpha \partial_z \xi_{k-1}(z), 
	\]
	where $\xi_k(z)$  is the root of $P_k$  such that for all $z\in\Omega$, $\xi_k(z) = m_0(z) +O(\alpha^{-1})$.
\end{Def}

We can deduce the following properties of $(\xi_k)_{k\geq0}$.

\begin{Prop} \label{prop:m}
	For all $k\geq1$, $m_k$ is analytic in $\Omega$ and satisfies:
	\begin{equation}\label{mind}
		\xi_k  = \xi_{k-1} +  O(\alpha^{-k})
	\end{equation}
	where the error term is controlled in $\|\cdot\|_{\mathcal{C}^1(\Omega)}$. Moreover, for all $z\in\Omega$:
	\begin{equation} \label{Pfact}
		P_k(s) = (s- \xi_k(z))(s-\widetilde{\xi}_k(z)),
	\end{equation}
	where there exists $c>0$, such that for all $z\in \Omega$, $-\im\widetilde{ \xi}_k(z)  >c/2$. Finally, there exists $C_k>0$ such that for any $k\geq1$:
	\[
	\| \xi_k \|_{\mathcal{C}^1(\Omega)} \le C_k ,\qquad\qquad \| \widetilde{\xi}_k \|_{\mathcal{C}^1(\Omega)} \le C_k . 
	\]
	
\end{Prop}

\begin{Pro} 
	We proceed by induction on $k\geq0$ using that, according to Lemma~\ref{lem:eq}, 
	$$\|m_0\|_{\mathcal{C}^1(\Omega)}=\|\xi_0\|_{\mathcal{C}^1(\Omega)}\le C_0$$  and $| \triangle_0(z)| \ge c>0$ for $z\in \Omega$. 
	Moreover, one also has that for all $k\geq0$, $z\in\Omega\mapsto\mu_k(W_z)$ is holomorphic by Proposition \ref{prop: nu holomoprhic}, thus $\|z\mapsto\mu_k(W_z)\|_{\mathcal{C}^1(\Omega)} \le C_{k-1} $.
	
	We assume that there exists $k\ge 0$ such that $\xi_k  = \xi_{k-1}+  O(\alpha^{-k})$ holds in $\Omega$ where the error is analytic in $\Omega$ and controlled in $\|\cdot\|_{\mathcal{C}^1(\Omega)}$ and the discriminant of $P_k(\cdot)$ satisfies $| \triangle_k(z)| \ge 2^{-k}c $ and $\| \triangle_k \|_{\mathcal{C}^1(\Omega)} \le  C_k$.

	By construction, we can write $P_{k+1}(s) = P_k(s) + O(\alpha^{-(k+1)})$ in $\Omega$ where the perturbation is analytic in $\Omega$ and controlled in $\|\cdot\|_{\mathcal{C}^1(\Omega)}$ (by hypothesis, one has $\|\mu_{k+1}(W_z)\|_{\mathcal{C}^1(\Omega)}  , \|\xi_k\|_{\mathcal{C}^1(\Omega)} \le C_k$). 
	Then, the discriminant of the polynomial $P_{k+1}(\cdot)$ satisfies $\triangle_{k+1}=\triangle_k - O(\alpha^{-(k+1)}) $ so that since $\alpha\gg 1$, 
	\begin{equation} \label{deltacond}
		|\triangle_{k+1}(z)| \ge c/2^{k+1} \text{ for $z\in\Omega$ \quad\quad\quad  and }\quad\quad\quad
		\| \triangle_{k+1} \|_{\mathcal{C}^1(\Omega)} \le  C_{k+1}. 
	\end{equation}
	This implies that there exists a  branch $\sqrt{\cdot}$ such that the roots of $P_{k+1}$ (see \eqref{Pfact}) are analytic in $\Omega$ satisfy the following expansion and are given by:
	\[\begin{cases}
		2\xi_{k+1}(z)&\hspace{-0,3cm}=  -V'(z)+\sqrt{\triangle_{k+1}(z)} = 2\xi_k(z) +O(\alpha^{-(k+1)}) \\
		2	\widetilde{\xi}_{k+1}(z)&\hspace{-0,3cm}=-V'(z)- \sqrt{\triangle_{k+1}(z)}= 2\widetilde{\xi}_{k}(z) +O(\alpha^{-(k+1)}) 
	\end{cases}.\]
	Then, using \eqref{deltacond}, we verify that 
	$\| \xi_{k+1} \|_{\mathcal{C}^1(\Omega)} , \| \widetilde{\xi}_{k+1} \|_{\mathcal{C}^1(\Omega)} \le  C_{k+1}$. 
	
	Moreover,  according to \eqref{roots}, since there exists $c>0$ such that $|\im V'(z)| \le c$  ($V$ is real-analytic), taking $\eta_0$ sufficiently small, we deduce that there exists $c_1>0$ such that $-\im  \xi_0(z)\ge c_1$ for $z\in\Omega$
	.
	Hence, for any $k\in\N$ if $\alpha\gg1$ (depending on $k$), one also has for $z\in\Omega$, 
	\[ 
	-\im \xi_k(z) = -\im \xi_0(z) + O(\alpha^{-1}) \ge c_1/2 . 
	\]
	This completes the proof. 
\end{Pro}

\begin{Rem}\begin{enumerate}
		\item  In the one-cut regular case, choosing $B=[-1+\delta,1-\delta]$ for  $\delta>0$,  one has $c= c_0 \delta$ where the constant $c_0>0$ depends on $V$ and $\eta_0>0$.
		\item The same arguments show that for any $q\geq1$, there are constants so that for any $k\geq1$, $\| \xi_k \|_{\mathcal{C}^q(\Omega)} \le C_{k,q}$. 
	\end{enumerate}
\end{Rem}

The sequence of functions $(\xi_k(z))_{k\geq1}$ will be instrumental to re-center the Stieltjes transform $s_N$ of the empirical measure and obtain optimal fluctuation bounds. These functions are conceivably
related to the linear functionals $(\mu_k)_{k\geq1}$ as explained in the next proposition. We recall that $m_k$ stands for the Stieltjes transform of $\mu_k$.

\begin{Prop}\label{prop:approx zeta m}
	Let $k\geq0$ and $z\in\Omega$, one has $$\| \xi_k - m_k\|_{\mathcal{C}^1(\Omega)} =O(\alpha^{-(k+1)}).$$ 
\end{Prop}

\begin{Pro}
	We proceed by induction on $k\in\N$. In particular, the claim is true for $k=0$ since $\xi_0$ is exactly the Stieltjes transform of the equilibrium measure $\mu_\mrm{eq}$. 
	According to Proposition~\ref{prop:mun stieltjes}, the Stieltjes transform of $\mu_k$ satisfies for $z\in\Omega$:
	\[
	m_{k}^{2}(z)+V'(z)m_{k}(z)+\mu_k(W_z)+\dfrac{2}{\alpha}\partial_zm_{k-1}(z) + O(\alpha^{-(k+1)}) = 0
	\]
	using that $m_{2k} = m_{k} +  O(\alpha^{-(k+1)})$ and $\mu_{2k}(W_z) = \mu_{k}(W_z)+  O(\alpha^{-(k+1)})$ where the errors are controlled in $\|\cdot\|_{\mathcal{C}^1(\Omega)}$.
	Using the induction hypothesis, $\xi_{k-1} = m_{k-1} +O(\alpha^{-k})$, we obtain:
	\[
	P_k(m_k) = O(\alpha^{-(k+1)}).
	\]
	Then, using the factorisation \eqref{Pfact} from Proposition~\ref{prop:m}, one has  $|m_k-\widetilde{\xi}_k| \ge 
	-\im \widetilde{\xi}_k(z)  >c/2$ since $\im m_k \ge0$ by Proposition \ref{prop:mun stieltjes} so that: 
	\[
	(m_k-\xi_k) =  O(\alpha^{-(k+1)}) .
	\]
\end{Pro}

The goal of the next Proposition is to rewrite the main term in the \emph{loop equation} from Lemma~\ref{lem:Leq} by using the polynomail $P_k$. We record the relevant results. In the sequel, we denote for  $k\geq1$:
\begin{equation} \label{Delta}
	s_{N,k}(z)  \defi  s_N(z) -\xi_{k}(z) , \qquad z\in\Omega. 
\end{equation}

\begin{Prop} \label{prop:Leq}
	Let $k\geq1$, $G(z)\in\mathcal{C}^1(\R^N)$ and  $z\in\Omega$, then:
	\begin{equation} \label{Leq}
		\E_N\left[\left( P_k(s_N) +  \Delta_k(z)+  \dfrac{2}{\alpha'}   \partial_zs_{N,k-1}(z)  +O(N^{-1})\right)G(z)\right] 
		= \frac1{\alpha N} 	\E_N\bigg[ \sum_{j=1}^{N} \frac{\partial_{\lambda_j}G(z)}{\lambda_j-z}  \bigg] 
	\end{equation}
	where $\alpha' \defi (\frac1\alpha-\frac1{2N})^{-1}$ and where the error term is deterministic and controlled uniformly in $\Omega$ (the implied constants depend on $k$). 
	Moreover, there is a constant $C_k\ge 0$ such that for any $k\in\N$ and $z\in \Omega$:
	\begin{equation} \label{stab}
		|s_{N,k}(z)| \le C|P_{k}(s_N(z))|\wedge|P_{k}(s_N)|^{1/2} . 
	\end{equation}
\end{Prop}

\begin{Pro}
	Using the previous definitions, one has for all $k\in\geq1$ and $z\in\Omega$, 
	\[
	s_N^2(z) + V'(z) s_N(z)  + \bm{\mu}_N(W_z) + \left( \tfrac{2}{\alpha}-\tfrac{1}{N}\right) \partial_z s_N(z)
	= P_k(s_N(z)) +  \Delta_k(z)+  \tfrac{2}{\alpha'} \partial_z s_{N,k-1}(z)    -\tfrac{1}{N} \partial_z \xi_{k-1}(z) . 
	\]
	Then, by Lemma~\ref{lem:Leq}, this gives~\eqref{Leq} with $-\tfrac{1}{N} \partial_z \xi_{k-1}  = O(N^{-1})$
	(by Proposition~\ref{prop:m}).  Finally, by \eqref{Pfact}, one has $P_{k}(s_N) = s_{N,k}(s_N-\widetilde{\xi}_k)$ and since $\im(s_N-\widetilde{\xi}_k) >c/2$ for all $z\in \Omega$ if $\alpha\gg1$. This yields the first bound in \eqref{stab}. For the second one, we use the fact that for some $c>0$:
	$$|s_N-\widetilde{\xi}_{k}|\geq|\im{(s_N-\widetilde{\xi}_{k})}|\geq-c\im{\widetilde{\xi}_{0}}\geq c'/2>0.$$
	Now, using that:
	$$|P_k(s_N)|^{1/2}\geq|s_{N,k}|\wedge|s_N-\widetilde{\xi}_{k}|,$$
	we obtain the conclusion.
\end{Pro}

\subsection{A priori estimates}\label{subsec3: a priori BMP}
In this section, we review the local law from \cite{bourgade2022optimal} and present a simplified version of their proof in the bulk. This local law (Proposition~\ref{prop:BMP}) is optimal only in the regime where $\beta$ is fixed but it will only be an input to derive an optimal local law at intermediate temperature  ($\alpha\gg(\log N)^{k}$ for all $k\geq1$) in the \emph{random matrix regime}  (on scales $\eta \gg \alpha^{-1}$). 
The key idea from \cite{bourgade2022optimal} is to apply  \eqref{Leq}  with $G(z)  = |P_0(s_N)|^{2(q-1)} P_0(s_N)^* $ for $q\geq1$ in order to obtain some inequality relating $|P_0(s_N)|^{2q}$ and $|P_0(s_N)|^{2(q-1)}$ to estimates the moments of $P_0(s_N)$. 
We will use variations of these arguments. In particular, we record the following \emph{loop inequality}. 

\begin{Lem} \label{lem:Lineq}
	There exists a constant $C>0$ such that for any  $k\geq1$, $q\geq1$ and $z\in\Omega$:
		\begin{equation*}
	\E_N\left[ |P_k(s_N)|^{2q}  \right] 
	\leq C^{q}	\Ec{\dfrac{q^{q}}{(N\alpha\eta^{2})^{q}}\big(|s_{N,k}(z)| + 1 \big)^{q}+|\Delta_k(z)|^{2q}+  \dfrac{1}{\alpha^{2q}}   |\partial_zs_{N,k-1}(z)|^{2q} +N^{-2q} }
\end{equation*}
\end{Lem}

\begin{Pro}
	Let $k\geq1$, $z\in\Omega$ and $q\geq1$, we apply Proposition~\ref{prop:Leq} with $G(z)  = |P_k(s_N)|^{2(q-1)} P_k(s_N)^* $.
	First, observe that if $s\mapsto G(s)$ is a complex polynomial and $G(z) = G(s_N(z))$, one has:
	\begin{equation}\label{Gbound}
		\begin{aligned}
			\left | {\displaystyle  \sum_{j=1}^{N} \dfrac{\partial_{\lambda_j}G(z)}{\lambda_j-z}} \right | 
			&\le \big( |\partial_s G(s_N(z))| +  |\overline{\partial}_{s} G(s_N(z))| \big) \left | {\displaystyle  \sum_{j=1}^{N} \dfrac{\partial_{\lambda_j}{s_N}(z)}{\lambda_j-z}} \right | \\
			&\le 2\eta^{-2} | s_N(z)| \big(  |\partial_s G(s_N(z))| \vee | \overline{\partial}_{s}  G(s_N(z))| \big)
		\end{aligned},
	\end{equation}
	using the complex chain rule and that $$ \sum_{j=1}^{N} \dfrac{\partial_{\lambda_j}s_N(z)}{\lambda_j-z} =- \partial^{2}_zs_N(z)$$ and we use the bound  $|\partial^{2}_zs_N(z)| \le \eta^{-2} \im s_N(z) $ for all $z\in\Omega$.
	
	Then, applying \eqref{Gbound} to $G = P_k(s)^{q-1}P_k^*(s)^{q}$, using that $\partial_s P_k(s) = 2(s-\xi_{k}) +O(1)$ according to Proposition~\ref{prop:m} (uniformly in $z\in\Omega$), there exists a constant $C_k>0$ such that: 
	\[
	|\partial_s G(s_N)| \vee | \overline{\partial}_{s}  G(s_N)| \le q \big(|s_{N,k}| +C_k \big) |P_k(s_N)|^{2(q-1)}.
	\]
	Similarly, one can bound $| s_N|  \le |s_{N,k}| + |\xi_{k}| $ so that, combined with \eqref{Gbound}, we obtain :
	\[
	\left | {\displaystyle  \sum_{j=1}^{N} \frac{\partial_{\lambda_j}G(z)}{\lambda_j-z}} \right |  \le 2q \eta^{-2} 
	\big(|s_{N,k}(z)| + C_k \big)^2 |P_k(s_N)|^{2(q-1)}.
	\]
	This leads to for all $\lambda>0$:
	\begin{multline*}
			\E_N\left[ |P_k(s_N)|^{2q}  \right] 
		\leq	\E_N\bigg[\lambda| P_k(s_N)|^{2(q-1)}  \cdot\frac{q}{\alpha N\eta^{2}\lambda} \big(|s_{N,k}(z)| + C_k \big)^2  \bigg] \\+\Ec{\lambda|P_k(z)|^{2q-1}\cdot\dfrac{1}{\lambda}(|\Delta_k(z)|+  \dfrac{2}{\alpha'}   |\partial_zs_{N,k-1}(z) | +O(N^{-1}))}.
	\end{multline*}
	We then apply Young's inequality ($xy\leq \frac{x^{a}}{a}+\frac{y^{b}}{b}$ if $a^{-1}+b^{-1}=1$) with $a=\frac{2q}{2(q-1)}$ and $b=q$ $a=\frac{2q}{2q-1}$ and $b=2q$ and obtain:
		\begin{multline*}
		\E_N\left[ |P_k(s_N)|^{2q}  \right] 
		\leq	(\lambda^{\frac{q}{q-1}}+\lambda^{\frac{2q}{2q-1}})\E_N\left[ |P_k(s_N)|^{2q}  \right] + \dfrac{q^{q}}{\lambda^{q}(N\alpha\eta^{2})^{q}}	\E_N\bigg[\big(|s_{N,k}(z)| + C_k \big)^{q}  \bigg] \\+\Ec{\dfrac{1}{2q\lambda}(|\Delta_k(z)|^{2q}+  \dfrac{2}{\alpha'^{2q}}   |\partial_zs_{N,k-1}(z) |^{2q} +O(N^{-2q}))}.
	\end{multline*}
	Taking $\lambda=1/10$ allows to absorb the first in the RHS in the LHS.
	Using that $P_k(s_N) G  = |P_k(s_N)|^{2q}$ and $|s_{N,k}| \le C|P_{k}(s_N)|$, by Proposition~\ref{prop:Leq}, we conclude that  there exists constants $C>0$ such that:
		\begin{equation*}
	\E_N\left[ |P_k(s_N)|^{2q}  \right] 
	\leq C^{q}	\Ec{\dfrac{q^{q}}{(N\alpha\eta^{2})^{q}}\big(|s_{N,k}(z)| + C_k \big)^{q}+|\Delta_k(z)|^{2q}+  \dfrac{1}{\alpha^{2q}}   |\partial_zs_{N,k-1}(z)|^{2q} +N^{-2q} }
	\end{equation*}
	This completes the proof. 
\end{Pro}

Then, we exploit the bounds from Lemma~\ref{lem:Lineq} in the case $k=0$ to prove a first local law analogous to \cite{bourgade2022optimal}. We note that this local law is sharp only in the regime where $\beta$ is fixed. The proof is an adaptation of  \cite{bourgade2022optimal} valid at arbitrary temperatures $\alpha\gg(\log N)^{1+\kappa}$.

\begin{Prop} \label{prop:BMP}
	There exists $C> 0$ such that for any $q\geq1$ and $z\in \Omega$ with $\alpha\eta \ge C$: 
	\[
	\|s_N(z) -m_{0}(z)\|_q \lesssim \sqrt{\frac{ q}{\alpha N\eta^2}}
	+  \frac{1}{\alpha\eta} .
	\]
\end{Prop}

\begin{Pro}
	Let $z\in\Omega$ and set $Z \defi   |P_0(s_N)|^2$.
	We first apply Lemma~\ref{lem:Lineq} with $k=0$ $(s_{N,0}=s_N)$ using the trivial bound, to get:
	\[
	|\partial_z s_{N}| \le \frac1N \sum_{j=1}^{N} \frac{1}{|z-\lambda_j|^2} \le \eta^{-1}  \im s_N . 
	\]
	Then, by re-centering and using the bounds $\| \xi_0 \|_{\mathcal{C}^1(\Omega)} \le C$ and  \eqref{stab}  (with $s_{N,0}=s_N-\xi_0$)with Proposition  \ref{prop:Leq}, we obtain:
	\[
	|  \partial_zs_{N}| \lesssim \eta^{-1}  \big (1+ |P_0(s_N)| \big) .
	\]
	By  Lemma~\ref{lem:Lineq} and Proposition \ref{prop:Leq}, this yields that there is a constant~$C>0$ such that for any $q\geq1$: 
	\begin{equation} \label{Leq0}
		\E_N\left[|P_0(s_N)|^{2q}\right]   \le \frac{C^{q} q^{q}}{(\alpha N\eta^2)^{q}} \E\left[|P_0(s_N)|^{q}+1 \right] 
		+ C^{q}\E_N\left[  |\Delta_0|^{2q}+  \frac{1}{(\alpha\eta)^{2q}}   + \frac{1}{N^{2q}} \right]. 
	\end{equation}
	Now using Corollary \ref{cor: glocon} and that $x\leq a\sqrt{x}+b$ implies $x\leq a^{2}+b$  when $x,a,b\geq0$, implies that:
		\begin{equation*} 
		\E_N\left[|P_0(s_N)|^{2q}\right]   \le \frac{C^{q} q^{2q}}{(\alpha N\eta^2)^{2q}} +\frac{C^{q} q^{q}}{(\alpha N\eta^2)^{q}}
		+ C^{q} \left(\dfrac{q^{q}}{(\sqrt{N\alpha})^{q}}+\dfrac{1}{\alpha^{4q}}+  \frac{1}{(\alpha\eta)^{2q}}   + \frac{1}{N^{2q}}\right) . 
	\end{equation*}
	When $q\leq N\eta$, we obtain:
	\begin{align*}
		\|s_N(z)-\xi_0(z)\|_{2q}\lesssim	\|P_0(s_N)\|_{2q}&\lesssim  \sqrt{\frac{ q}{\alpha N\eta^2}}
	+  \frac{1}{\alpha\eta}.
	\end{align*}
	When $q\geq N\eta$, we use:
		\begin{align*}
		\|s_N(z)-\xi_0(z)\|_{2q}\lesssim	\|P_0(s_N)\|_{q}^{1/2}&\lesssim  \sqrt{\frac{ q}{\alpha N\eta^2}}
+  \frac{1}{\alpha\eta}.
	\end{align*}
\end{Pro}

\subsection{Random matrix regime}\label{subsec3:RMT regime}

We now use the bounds from Lemma~\ref{lem:Lineq} for $k\ge 1$ to improve the local law from Proposition~\ref{prop:BMP} by upgrading the re-centring to improve the concentration. The proof proceeds by induction on $k\ge 1$ starting from Proposition~\ref{prop:BMP} and eventually leads to  optimal bounds on the fluctuations in the \emph{random matrix regime} away from the critical window. 

\begin{Thm}[Local law -- RMT regime] \label{prop:LL1}
	Let $k\geq0$, for any $q\ge 1$ and $z\in \Omega$ with $\alpha\eta \ge C$, 
	\[
		\|s_{N}(z)-\xi_k(z) \|_q \lesssim  \sqrt{\frac{q}{\alpha N\eta^2} } +\frac{1}{(\alpha\eta)^{k+1}} . 
	\]
	In particular, if $\alpha \eta \geq N^\delta$ for some $\delta>0$, then for $k\geq0$ large enough we obtain the sharp bounds for any $q\ge 1$:
	\[
	\|s_{N}(z)-\xi_k(z) \|_q \lesssim   \sqrt{\frac{q}{\alpha N\eta^{2}}} .
	\]
\end{Thm}
\begin{Pro}
	Let $z\in\Omega$, t he starting point is again the inequality from Lemma~\ref{lem:Lineq}  for $k\ge 0$.
	Observe that to control the RHS, using the induction hypothesis (the bound holds for $k=0$ by Proposition~\ref{prop:BMP}), we can control (using Cauchy's formula) for any $q\geq1$: 
	\begin{equation}
		\begin{aligned} \label{sest}
			\|\partial_z s_{N,k}(z)\|_q &\le \eta^{-1} \sup_{w\in\overline{B}(z,\eta/2)}\| s_{N,k}(w)\|_{q} \\
			\alpha^{-1} \|\partial_z s_{N,k}(z)\|_q, &\lesssim \frac{1}{(\alpha\eta)^{k+1}}+\dfrac{1}{\alpha\eta}\sqrt{\frac{ q}{\alpha N\eta^2}}.
		\end{aligned}
	\end{equation}
%
%
By Lemma~\ref{lem:Lineq}  we have:
	\[
	\|  P_k(s_N)\|_{2q} \lesssim \sqrt{\frac{q}{\alpha N\eta^2}}\big(\|s_{N,k}(z)\|_q^{1/2} + 1 \big) + \|\Delta_k\|_{2q}+  \frac{\|\partial_zs_{N,k}\|_{2q}}{\alpha} + \frac{1}{N} .
	\]
	We thus obtain, using again, $|s_{N,k}(z)|\lesssim |P_k(s_N)|$ by Proposition \ref{prop:Leq} and the fact that $x\leq a\sqrt{x}+b$ implies that $x\leq a^{2}+b$:
	\[
	\|  P_k(s_N)\|_{2q} \lesssim \sqrt{\frac{q}{\alpha N\eta^2}}+\frac{q}{\alpha N\eta^2}+ \|\Delta_k\|_{2q}+  \frac{\|\partial_zs_{N,k}\|_{2q}}{\alpha} + \frac{1}{N} .
	\]
	The $\Delta_k$-term is controlled by Corollary \ref{cor: glocon} and the $\partial_zs_{N,k}$-term is controlled by \eqref{sest}.  Adjusting the constants, we obtain  that if $\alpha\eta \ge C$:
\begin{align*}
		\|  P_k(s_N)\|_{2q} &\lesssim\sqrt{\frac{q}{\alpha N\eta^2} }+  \frac{q}{\alpha N\eta^2} +\dfrac{1}{\alpha^{k+1}}+\sqrt{\dfrac{q}{N\alpha}} +\frac{1}{(\alpha\eta)^{k+1}}+\dfrac{1}{\alpha\eta}\sqrt{\frac{ q}{\alpha N\eta^2}}
+\dfrac{1}{N}
	\\&\lesssim\sqrt{\frac{q}{\alpha N\eta^2} }+  \frac{q}{\alpha N\eta^2}   +\frac{1}{(\alpha\eta)^{k+1}}
\end{align*}
	the terms $\frac{1}{\alpha\eta}\sqrt{\frac{ q}{\alpha N\eta^2}}$, $\sqrt{\frac{q}{N\alpha}}$ and $\frac{1}{N}$ being negligible.
	Now using for $q\leq N\eta$, $|s_{N,k}|\lesssim|P_k(s_N)|$, we obtain:
	\[
	\|  s_{N,k}(z)\|_{2q} \lesssim  \sqrt{\frac{q}{\alpha N\eta^2} } +\frac{1}{(\alpha\eta)^{k+1}}.
	\]
	For $q\geq N\eta$, we use $|s_{N,k}|\lesssim|P_k(s_N)|^{1/2}$, to obtain:
		\[
	\|  s_{N,k}(z)\|_{2q} \lesssim 	\|  P_k(s_N)\|_{q}^{1/2}\lesssim \sqrt{\frac{q}{\alpha N\eta^2} } +\frac{1}{(\alpha\eta)^{k+1}}.
	\]
	Finally, replacing $ Z \defi   |P_k(s_N)|^2$ and using the bound \eqref{stab}, this completes the proof. 
\end{Pro}

\subsection{Poisson regime}\label{subsec3:poisson}
The optimal local law Theorem \ref{prop:LL1} only holds above the transition \textit{i.e.} for $\eta\geq C\alpha^{-1}$. While the bound gets worse as $\eta$ decreases, an interesting phenomenon occurs. The bound saturates and changes from $\frac{1}{\sqrt{N\alpha}\eta}$ above the transition to $\frac{1}{\sqrt{N\eta}}$. In this regime, the largest term in the loop equation Proposition \ref{prop:Leq} is not the term in $P_k(s_N)$ like in the RM regime but the term $\alpha^{-1}\partial_z s_N(z)$. As the latter term is the contribution of the entropy and $P_k$ encodes the main contribution of the interaction and the confining potential, we can interprent the Poisson regime as the regime where entropy dominates energy.

\begin{Lem}\label{lem:gronwall square}
	Let $f\in\mcal{C}^{1}([a,b])$ be a positive function such that there exists positive functions $\phi,\psi\in\mcal{C}^{0}([a,b])$ such that for all $x\in[a,b]$
	$$-f'(x)\leq\dfrac{\phi(x)}{f(x)}+\psi(x),$$
	then
	$$f(a)\leq\int_{a}^{b}\psi(x)\diff x+\left(f(b)^{2}+2\int_{a}^{b}\phi(x)\diff x\right)^{1/2}.$$
\end{Lem}

\begin{Pro}Let $g(x)\defi\int_{x}^{b}\psi(x)\diff x+\left(f(b)^{2}+2\int_{x}^{b}\phi(x)\diff x\right)^{1/2}$, notice that:
	$$-g'(x)=\psi(x)+\dfrac{\phi(x)}{\left(f(b)^{2}+2\int_{x}^{b}\phi(x)\diff x\right)^{1/2}}\geq \psi(x)+\dfrac{\phi(x)}{g(x)}$$suppose there exists a maximal interval $[u,v]$ in $\{f>g\}$ with $f(v)=g(v)$, thus for all $x\in[u,v]$
	$$-f'(x)\leq\dfrac{\phi(x)}{f(x)}+\psi(x)< \dfrac{\phi(x)}{g(x)}+\psi(x)\leq-g'(x).$$
	Thus $g-f$ is non-increasing while equal to 0 at $v$, thus $g\geq f$ which is a contradiction. This proves the result.
\end{Pro}

\begin{Thm}[Local law -- Poisson regime] \label{prop:LL2}Let $\delta>0$, assume that $\alpha \ge N^{\delta}$, then for all $k\geq0$ large enough (depending only on $\delta>0$) and $q\ge 1$, $z\in \Omega$ with $\eta  \le \alpha^{-1/2}$:
	\[
	\| s_N(z) - \xi_k(z) \|_q \lesssim  \sqrt{\frac{q}{N\eta}}+\dfrac{q}{N\eta}+\dfrac{1}{\alpha^{k+1}}. 
	\]
	Thus for $k\geq0$ large enough, for any $q\geq1$:
		\[
	\| s_N(z) - \xi_k(z) \|_q \lesssim  \sqrt{\frac{q}{N\eta}}+\dfrac{q}{N\eta}. 
	\]
\end{Thm}

\begin{Pro}
	Here $k\geq0$ is fixed such that  $ \frac{1}{\eta^2 N}  \gg \frac{1}{\alpha^k} $ and we set $G = q|s_{N,k}|^{2(q-1)} s_{N,k}^*$ in the loop equation from Proposition~\ref{prop:Leq}. 
	Like in the proof of Lemma~\ref{lem:Lineq}, $G$ is a complex polynomial in $s_N$, so using \eqref{Gbound}, we obtain:
	\begin{equation}\label{Gbound2}
		\begin{aligned}
			\left | {\displaystyle  \sum_{j=1}^{N} \dfrac{\partial_{\lambda_j}G(z)}{\lambda_j-z}} \right | 
			& \le 2 \eta^{-2} q^2 | s_N(z)|\cdot |s_{N,k}(z)|^{2(q-1)} \\
			& \le  2 \eta^{-2} q^2 C_k |s_{N,k}(z)|^{2(q-1)}  + 2 \eta^{-2} q^2 |s_{N,k}(z)|^{2q-1}
		\end{aligned}
	\end{equation}
	using that $|m_k(z)| \le C_k$ for $z\in\Omega$.
	
	Then, using that $\xi_{k-1}= \xi_{k} + O(\alpha^{-k})$ according to \eqref{mind}, taking the imaginary part on the RHS of \eqref{Leq} and using the factorization \eqref{Pfact} for $P_k$, one has for all $z\in\Omega$:
	\begin{equation} \label{Leq2}
		\im\left\{\left( P_k(s_N) +   \dfrac{2}{\alpha'}  \partial_zs_{N,k}(z) \right)G(z)\right\} 
		= - \dfrac1{\alpha'}\partial_\eta [|s_{N,k}|^{2q}](z)  + q|s_{N,k}(z)|^{2q}  \im(s_N(z)-\widetilde{\xi}_k(z))    
	\end{equation}
	where we used that 
	\[
	-\partial_\eta |s_{N,k}|^{2q}  = 2  \im\big(s_{N,k}^{*q} \partial_z \overline{s_N}^q) \big)
	= 2  \im\big(G  \overline{s_N}' \big) .
	\]
	The key observation is that,  by Proposition~\ref{prop:m},  $ \im(s_N-\widetilde{\xi}_k) > 0$ for $z\in \Omega$ so that the second term on the RHS of \eqref{Leq2} is non-negative. 
	Thus, combining \eqref{Leq} with \eqref{Gbound2} and \eqref{Leq2}, there is a constant $C_k>0$ such that for any $q\geq1$ and $z\in\Omega$:
	\[
	\E_N\big[ -\partial_\eta |s_{N,k}(z)|^{2q} \big]  \le \frac{\alpha'}{\alpha} C_k \E_N\left[ \frac{q^2}{ N\eta^2} |s_{N,k}(z)|^{2q-2}+\bigg( \frac{\alpha |\Delta_k(z)|}2 + \frac{q^2}{\eta^2 N} + \frac{ \alpha}{N}  + \frac{1}{\alpha^k} \bigg)   |s_{N,k}(z)|^{2q-1} \right]. 
	\]
	At this stage, we observe that $\alpha' \le 2\alpha$ and using the conditions  $ \frac{1}{\eta^2 N}  \gg \frac{1}{\alpha^k} $ and $\eta \le \alpha^{-1/2}$, the last two terms can be neglected.
	Hence, there is a constant $C>0$ such that for any $q\geq1$ and $z\in\Omega$ and $\eta\leq\alpha^{-1/2}$:
	\begin{equation} \label{Leq4}
		\E_N\big[ -\partial_\eta |s_{N,k}(z)|^{2q} \big]  \le \frac{C q^2}{ N\eta^2}\E_N \left[|s_{N,k}(z)|^{2(q-1)} \right] 
		+  C\E_N\left[\bigg(  \alpha |\Delta_k(z)| + \frac{ q^2}{\eta^2 N} \bigg)   |s_{N,k}(z)|^{2q-1} \right]  . 
	\end{equation}
	\textbf{1. Case $q\leq N\eta$.} Now, let   $u_q : \eta \in\R_+ \mapsto   \|s_{N,k}(z)^2\|_{q}  $ for $q\ge 0$.
	By H\"older's inequality, we have the following bound:
	\[
	\E_N\left[   |\Delta_k| \cdot |s_{N,k}|^{2q-1}  \right] \le \| \Delta_k \|_{2q} \cdot\|s_{N,k} \|_{2q}^{{2q-1}} \lesssim\sqrt{\dfrac{q}{N\alpha}}  u_{q}^{q-\frac12} 
	\]
	where we used the bound Corollary \ref{cor: glocon} at last (the term proportional to $\tfrac{1}{\alpha^{k+1}}$ is again negligible).
	
	Observe that $q\mapsto u_q$ is non-decreasing, so we deduce from \eqref{Leq4} that:
	\begin{equation} \label{Leq5}
		- \partial_\eta(u_q^{q})(\eta) = \E_N\big[ -\partial_\eta |s_{N,k}|^{2q}(z) \big]  
		\le \frac{C q^2}{N\eta^2}   u_q^{q-1}(\eta) +  C\bigg( \sqrt{\dfrac{q\alpha}{N}}   \frac{ q^2}{\eta^2 N} \bigg) u_q^{q-\frac12}(\eta) . 
	\end{equation}
	
	This inequality is almost suitable to apply Gr\"onwall's estimate. 
	Multiplying \eqref{Leq5} by $\tfrac1{q}u_{q}^{1-q}$, we obtain  for any $\theta\ge 1$:
	\begin{equation} \label{Leq6}
		\begin{aligned}
			-  \partial_\eta(u_q) (\eta)
			&\le   \frac{C q}{N\eta^2} +  C\bigg(\frac{1}{\sqrt{N /\alpha}}+ \frac{q}{\eta^2 N} \bigg)\sqrt{u_q(\eta)}  \\
			&\le   \frac{C q}{ N\eta^2} + \frac{C}{2\theta\sqrt{N /\alpha}}  + C\theta   \bigg(\frac{1}{\sqrt{N /\alpha}}+ \frac{q}{\eta^2 N} \bigg)u_q (\eta)
		\end{aligned}
	\end{equation}
	by Young's inequality ($2\sqrt{u} \le \theta^{-1}+ \theta u$ for $u\ge 0$).

	Then, integrating the inequality \eqref{Leq6} between $[\eta, \epsilon]$ with $\eta\le \epsilon$ and $\epsilon \le \alpha^{-1/2}$,  we can apply Gr\"onwall Lemma.
	Choosing $\theta^{-1} = \frac{\epsilon}{\sqrt{N /\alpha}}$, we obtain if $N \eta \ge q$:
	\[
	\theta\int_\eta^{\epsilon} \bigg(\frac{1}{\sqrt{N /\alpha}}+ \frac{Cq}{t^2 N} \bigg) \diff t \le C +\frac{Cq}{N\eta} \le C. 
	\]	Similarly, with $\epsilon = \alpha^{-1/2}$ ($\theta=\sqrt{N}$), we also have:
\begin{align*}
		\int_\eta^{\epsilon}   \bigg(\frac{q}{t^2 N}+\frac{1}{\theta\sqrt{N /\alpha}}\bigg) \diff t
	\le  \frac{C q}{N\eta}  +  \dfrac{\epsilon}{\theta\sqrt{N/\alpha}} \lesssim\frac{C q}{N\eta}.
\end{align*}
	Applying Gr\"onwall's lemma to \eqref{Leq6}, we conclude that for $\eta \le \epsilon = \alpha^{-\frac12}$:
	\[
	u_q(\eta) \le u_q(\epsilon) + C\frac{ q}{N\eta} .
	\]
	Finally, by Theorem~\ref{prop:LL1} ($\alpha\epsilon = \epsilon^{-1}\ge N^{\delta/2}$ so that the first term dominates), 
	$\sqrt{u_q(\epsilon)} \le \sqrt{\frac{C q}{N}}$ so that for any $q\geq1$ and $\eta \in \left [\dfrac qN, \dfrac{1}{\sqrt{\alpha}}\right ]$:
	\[
	\| s_N(z) - \xi_k(z) \|_{2q}=\sqrt{u_q(\eta)}\lesssim   \sqrt{\frac{q}{N\eta}}. 
	\]
	\textbf{2. Case $q\geq N\eta$.} We start with \eqref{Leq4} and set $v_q(\eta)=\sqrt{u_q(\eta)}=\| s_{N,k}(z)\|_{2q}$, we thus obtain:
	$$- \partial_\eta(v_q^{2q})(\eta) =-2q \partial_\eta v_q(\eta)v_q^{2q-1}(\eta)
	\le \frac{C q^2}{N\eta^2}   v_q^{2q-2}(\eta) +  C\bigg( \sqrt{\dfrac{q\alpha}{N}} +    \frac{ q^2}{\eta^2 N} \bigg) v_q^{2q-1}(\eta).$$
	Dividing by $2qv_q^{2q-1}$ (because $|s_{N,k}(z)|$ is not equal to 0 almost surely), we obtain:
		$$ -\partial_\eta v_q(\eta)
	\le \frac{C q}{N\eta^2}   \dfrac{1}{v_q(\eta)} +  C\bigg( \sqrt{\dfrac{\alpha}{N}} +    \frac{ q}{\eta^2 N} \bigg).$$
	Thus by Lemma \ref{lem:gronwall square} and Theorem \ref{prop:LL1}, we obtain for $\eta_0=C/\alpha$:
	$$v_q(\eta)\leq C\int_{\eta}^{\eta_0}\left(\sqrt{\dfrac{\alpha}{N}} +    \frac{ q}{t^2 N}\right) \diff t+\left(v_q(\eta_0)^{2}+2\int_{\eta}^{\eta_0}\frac{C q}{Nt^2}  \diff t\right)^{1/2}\lesssim\dfrac{1}{\sqrt{N\alpha}}+\dfrac{q}{N\eta}+\sqrt{\dfrac{q}{N\eta}}\lesssim\dfrac{q}{N\eta}.$$
	Finally, we use Proposition \ref{prop:m} to extend the result to any $k\geq0$. This completes the proof. 
\end{Pro}

\begin{Rem}
	Combining Theorems \ref{prop:LL1} and \ref{prop:LL2} with Proposition \ref{prop:approx zeta m} gives Theorem \ref{thm:local law}.
\end{Rem}

\section{Mesosocopic concentration estimates}\label{sec: concentration}
The main goal of this section is to obtain explicit bounds on the error terms appearing in Stein's method (or equivalently loop equations) to obtain the mesoscopic CLT in the bulk Theorem \ref{thm:mesoc clt bulk}. These error terms are typically given in terms of the moments of the anisotropy $\|A_N^{(k)}(\phi_\eta)\|_q$ and of other linear statistics $\|L_N^{(k)}(\phi_\eta)\|_q$ (see Definition \ref{def:improved recentring}) for all $q\geq1$ where $\eta=\eta_N\tend{N\rightarrow\infty}0$ and $\phi\in\mcal{C}_c^{3}(\R)$. These bounds will be directly implied by using the Helffer–Sjöstrand formula and the sharp local laws obtained in Section \ref{sec: local laws}, more precisely Theorems \ref{prop:LL1} and \ref{prop:LL2}. An argument for the version we use can be found in \cite[Proposition C.1]{benaych2016lectures}.


\begin{Prop}[Helffer–Sjöstrand formula] \label{prop:Helffer–Sjöstrand}

	Let $r \in \N$ and $\phi \in \mcal{C}_c^{r+1}(\R)$. Let $\delta>0$. We define the pseudo-analytic extension of $\phi$ of degree $r$ by
	\[
	\Phi(z)=\Phi (x+\ii y) 
	\defi  \chi_\delta(y)\sum_{j=0}^r \frac{1}{j!} (\ii y)^j \phi^{(j)}(x),
	\]
	where $\chi_\delta=\chi(\delta^{-1}\cdot)$ and $\chi\in \mcal{C}_c^\infty(\R)$ is even and such that $\ind{[-1,1]}\leq\chi_\delta\leq \ind{[-2,2]} $.	Then, for all $\lambda \in \R$, and $r\geq1$:
	\[
	\phi(\lambda) = \re\left(\int_{\C_+}\frac{\overline{\partial}\Phi(z)}{\lambda-z}\frac{\diff^2 z}{\pi}\right).
	\]
\end{Prop}
Of course the above formula holds for generic cutoff $\chi;$ we introduce the scale dependent behavior of $\chi$ in this formula to make the notation cleaner in what follows.

We first apply the Helffer–Sjöstrand formula to show that a linear statistic $L_N^{(k)}(\phi)$ can be controlled by the imaginary part Stieltjes transform of the recentered empirical measure, $\im\widetilde{s_N}(z)$.


\begin{Lem}\label{eq:HS_L_N(f)}
	Let $k\geq0$, $\phi \in \mcal{C}_c^3(\R)$ supported on $[-1+\varepsilon,1-\varepsilon]$ and let $\Omega$ denote the neighborhood in Lemma \ref{lem:eq}. Let $\delta>0$ be arbitrary, and denote the Stieltjes transform of $\bm{\mu}_N-\mu_k$ by 
	\begin{equation*}
	\widetilde{s_N}(z)\defi s_N(z)-m_k(z).
	\end{equation*}
	Then, with $z=x+\ii y$ we have
	\begin{equation*} 
L_N^{(k)}(\phi)=-\frac{1}{2\pi}\int_{\Omega_+}\bigg(\phi''(x)y\chi_\delta(y)+\phi(x)\Big(2(\chi_\delta)'(y)+y(\chi_\delta)''(y)\Big)\bigg) \im \widetilde{s_N}(z)\diff^2z.
	\end{equation*}
\end{Lem}

\begin{Pro}
	We first apply Proposition \ref{prop:Helffer–Sjöstrand} with a degree $2$ pseudo-analytic extension and expand to obtain:
	\begin{align*}
	L_N^{(k)}(\phi)&=-\frac{1}{2\pi}\int_{\Omega_+}\phi''(x)y\chi_\delta(y)\im \widetilde{s_N}(z)\diff^2 z-\frac{1}{2\pi}\int_{\Omega_+}(\chi_\delta)'(y)\phi(x)\im \widetilde{s_N}(z)\diff^2z \\&\quad-\frac{1}{2\pi}\int_{\Omega_+}\phi'(x)(\chi_\delta)'(y)y\re \widetilde{s_N}(z)\diff^2z-\frac{1}{4\pi}\int_{\Omega_+}y\chi_\delta(y)\phi^{(3)}(x)\re \widetilde{s_N}(z)~\diff^2z \\&\quad+\frac{1}{4\pi}\int_{\Omega_+}(y\chi_\delta(y))'\phi''(x)\im\widetilde{s_N}(z)\diff^2z.
	\end{align*}
	The last two terms in fact cancel. Integrating by parts in $x$ and using Cauchy-Riemann,
	\begin{align*}
	-\frac{1}{4\pi}\int_{\Omega_+}y\chi_\delta(y)\phi^{(3)}(x)\re \widetilde{s_N}(z)~\diff^2z&=\frac{1}{4\pi}\int_{\Omega_+} y\chi_\delta(y)\phi''(x)\partial_x \re \widetilde{s_N}(z)\diff^2z \\
	&=\frac{1}{4\pi}\int_{\Omega_+} y\chi_\delta(y)\phi''(x)\partial_y \im \widetilde{s_N}(z)\diff^2z
	\end{align*}
	using the compact support of $\phi$ to guarantee the boundary term vanishes. Integrating by parts in $y$ now using that $y\chi_\delta(y)\equiv 0$ on the real axis, we find
	\begin{equation*}
	-\frac{1}{4\pi}\int_{\Omega_+}y\chi_\delta(y)\phi^{(3)}(x)\re \widetilde{s_N}(z)~\diff^2z=-\frac{1}{4\pi}\int_{\Omega_+} (y\chi_\delta(y))'\phi''(x)\im \widetilde{s_N}(z)\diff^2z.
	\end{equation*}
	This leaves us with 
	\begin{multline}\label{eq: firstfluctexp}
	L_N^{(k)}(\phi)=-\frac{1}{2\pi}\int_{\Omega_+}\phi''(x)y\chi_\delta(y)\im \widetilde{s_N}(z)\diff^2 z-\frac{1}{2\pi}\int_{\Omega_+}(\chi_\delta)'(y)\phi(x)\im \widetilde{s_N}(z)\diff^2z \\-\frac{1}{2\pi}\int_{\Omega_+}\phi'(x)(\chi_\delta)'(y)y\re \widetilde{s_N}(z)\diff^2z.
	\end{multline}
	
	We can also remove the $\phi'$ in the last term for easier estimation later on. Integrating by parts in $x$, using that $\phi$ is compactly supported 
	\begin{equation*}
	-\frac{1}{2\pi}\int_{\Omega_+}\phi'(x)(\chi_\delta)'(y)y\re \widetilde{s_N}(z)\diff^2z=\frac{1}{2\pi}\int_{\Omega_+}\phi(x)(\chi_\delta)'(y)y\partial_x(\re \widetilde{s_N}(z))\diff^2z
	\end{equation*}
	and using Cauchy-Riemann we may write 
	\begin{equation*}
	\frac{1}{2\pi}\int_{\Omega_+}\phi(x)(\chi_\delta)'(y)y\partial_x(\re \widetilde{s_N}(z))\diff^2z=\frac{1}{2\pi}\int_{\Omega_+}\phi(x)(\chi_\delta)'(y)y\partial_y(\im \widetilde{s_N}(z))\diff^2z.
	\end{equation*}
	Integrating by parts in $y$ using that $(\chi_\delta)'(y)y\im \widetilde{s_N}(z)$ vanishes on the real axis, we obtain 
	\begin{equation*}
	\frac{1}{2\pi}\int_{\Omega_+}\phi(x)(\chi_\delta)'(y)y\partial_y(\im \widetilde{s_N}(z))\diff^2z=-\frac{1}{2\pi}\int_{\Omega_+}\phi(x)\left((\chi_\delta)''(y)y+(\chi_\delta)'(y)\right)\im \widetilde{s_N}(z)\diff^2z.
	\end{equation*}
	Substituting into \eqref{eq: firstfluctexp} and simplifying yields the result.
\end{Pro}

We next use Theorems \ref{prop:LL1} and \ref{prop:LL2} and the above result to obtain the following control on fluctuations of linear statistics.

\begin{Prop}[Bound on linear statistics]\label{lem:HS_L_N}
	Let $\varepsilon>0$, and suppose $B \subset [-1+\varepsilon, 1-\varepsilon]$. Suppose $\alpha \geq N^\delta$,  for some $\delta\in (0,1]$, $q=q(N)\geq1$. Define for all $\phi\in \mcal{C}_c^3(B)$ and all $\eta\in(0,\eta_0]$, a constant $\mathtt{c}(\phi)>0$ such that
	for all $m\in\llbracket0,3\rrbracket$:
	\begin{equation}\label{eq: Cnorm}
		\|\phi^{(m)}\|_{L^1(\R)}\leq \mathtt{c}(\phi)\eta^{1-m}.
	\end{equation}Then for $k\geq1$ large enough (depending only on $\delta$) and all $\phi\in \mcal{C}_c^3(B)$,
we have: 
	\begin{equation*}
		\|L_N^{(k)}(\phi)\|_{q}
\lesssim	\mathtt{c}(\phi)\begin{cases}
\dfrac{qN^{\frac{2\delta}{3}}}{ N\alpha\eta}+\sqrt{\dfrac{q}{N\alpha}}  & \text{if }\alpha \eta \geq N^\delta, \\
\sqrt{\dfrac{q\eta}{N}}+\dfrac{q}{N}&\text{otherwise.}
\end{cases}
	\end{equation*}

\end{Prop}


\begin{Pro} 
We first choose the cutoff $\chi_\delta$ with $\delta=\eta$. We then write $L_N^{(k)}(\phi)$ using Lemma \ref{eq:HS_L_N(f)}; using the triangle inequality along with $y|(\chi_\eta)''(y)|\leq -(\chi_\eta)'(y)=\eta^{-1}|\chi'(\eta^{-1}y)|$, we obtain:
\begin{align*}\label{eq: secondfluctexp}
\left\|L_N^{(k)}(\phi)\right\|_q &\lesssim  \|\phi''\|_{L^1(\R)} \underbrace{\int_0^{2\eta} \hspace{-.2cm}  y\chi_\eta( y) \sup_{E\in B} \| \im \widetilde{s_N}(E+\ii y)\|_q   \diff y}_{\defi \mcal{I}_1}\nonumber
\\&\quad+O\bigg(
\|\phi\|_{L^1(\R)}  \underbrace{\sup_{\substack{E\in B\\y\in [\eta,2\eta]}} \| \im \widetilde{s_N}(z) \|_q}_{\defi \mcal{I}_2}  \bigg) 
\end{align*}
We consider the Random matrix and Poisson regimes separately. 

\textbf{1. Case $ \eta \geq \tfrac{N^{\delta}}{\alpha}$.}
First, suppose $ \eta \geq \tfrac{N^{\delta}}{\alpha}$. Then, applying Theorems \ref{prop:LL1} and \ref{prop:LL2}, we obtain with $k$ sufficiently large:
\begin{align*}
	\mcal{I}_1 \lesssim \int_0^{\frac{N^{\frac{2\delta}{3}}}{\alpha}}\left(\sqrt{\frac{q}{Ny}}+\frac{q}{Ny}\right)y\diff y+\int_{\frac{N^{\frac{2\delta}{3}}}{\alpha}}^{2\eta} \sqrt{\frac{q}{N\alpha}}\diff y \lesssim\sqrt{\frac{qN^{2\delta}}{\alpha^{3} N}}+\frac{qN^{\frac{2\delta}{3}}}{ N\alpha}+\sqrt{\frac{q\eta^{2}}{N\alpha}}\lesssim \frac{qN^{\frac{2\delta}{3}}}{ N\alpha}+\sqrt{\frac{q\eta^{2}}{N\alpha}}
\end{align*}
where we neglegted the term $\sqrt{\frac{qN^{2\delta}}{\alpha^{3} N}}$ because $\alpha \eta \geq N^{\delta}$.
Thus, we have the bound
\begin{equation*}
 \|\phi''\|_{L^1(\R)}\mcal{I}_1 \lesssim  \mathtt{c}(\phi)\eta^{-1}\left(\frac{qN^{\delta}}{ N\alpha}+\sqrt{\frac{q\eta^{2}}{N\alpha}}\right)= \mathtt{c}(\phi)\left(\frac{qN^{\delta}}{ N\alpha\eta}+\sqrt{\frac{q}{N\alpha}}\right).
\end{equation*}We also have 
\begin{equation*}
\|\phi\|_{L^1(\R)}\mcal{I}_2 \lesssim\mathtt{c}(\phi)\eta \sqrt{\frac{q}{N\alpha}}
\end{equation*}
from Theorem \ref{prop:LL1}; we thus obtain:
\begin{equation*}
\left\|L_N^{(k)}(\phi)\right\|_q \lesssim\mathtt{c}(\phi)\left(\frac{qN^{\delta}}{ N\alpha\eta}+\sqrt{\frac{q}{N\alpha}}\right).
\end{equation*}
\textbf{2. Case $ \eta \leq \tfrac{N^{\delta}}{\alpha}$.} If $\eta \leq\tfrac{N^{\delta}}{\alpha}$, the cutoff is equal to $1$ on $[0,\eta]$ and then decays outside. Using Theorem \ref{prop:LL2} yields 
\begin{equation*}
 \|\phi''\|_{L^1(\R)} 	\mcal{I}_1 \lesssim \mathtt{c}(\phi)\eta^{-1} \int_0^{2\eta}\left(\sqrt{\dfrac{qy}{N}}+\dfrac{q}{N}\right) \diff y 
	\lesssim\mathtt{c}(\phi)\left(\sqrt{\dfrac{q\eta}{N}}+\dfrac{q}{N}\right) 
\end{equation*}
and 
\begin{equation*}
\|\phi\|_{L^1(\R)}\mcal{I}_2 \lesssim \eta\mathtt{c}(\phi)\left(\sqrt{\frac{q\eta}{N}}+\frac{q\eta}{N}\right) .
\end{equation*}
Again using \eqref{eq: Cnorm}, we obtain 
\begin{equation*}
\left\|L_N^{(k)}(\phi)\right\|_q \lesssim \mathtt{c}(\phi)\left(\sqrt{\dfrac{q\eta}{N}}+\dfrac{q}{N}\right) ,
\end{equation*}
as desired.
\end{Pro}

Before we proceed to analyzing the anisotropy, we record the following extension of Proposition \ref{lem:HS_L_N} that will be useful in controlling some of the errors in our application of Stein's method in the following section.
\begin{Cor}\label{lem: genHS_L_N}
Under the same assumptions as in Proposition \ref{lem:HS_L_N}, let $m:\Omega_+\rightarrow \C$ satisfying $\|m\|_{\mcal{C}^0(\Omega_+)}\lesssim 1$. We have:

	\begin{equation*}
		\left\|\re\left(\int_{\Omega_+}\overline{\partial}\Phi(z)m(z)\widetilde{s_N}(z)\frac{\diff^2 z}{\pi}\right)\right\|_{q}
\lesssim\mathtt{c}(\phi)\begin{cases}
	\dfrac{qN^{\frac{2\delta}{3}}}{ N\alpha\eta}+\sqrt{\dfrac{q}{N\alpha}}  & \text{if }\alpha \eta \geq N^\delta, \\
	\sqrt{\dfrac{q\eta}{N}}+\dfrac{q}{N}&\text{otherwise.}
\end{cases}
	\end{equation*}
	
\end{Cor}
\begin{Pro}
The proof is analogous to that of Proposition \ref{lem:HS_L_N}.
\end{Pro}
The key to obtaining estimates on the anistotropy $A_N^{(k)}(\phi)$ is the following analogue of Lemma \ref{eq:HS_L_N(f)}, which follows from applying the Helffer-Sj\"ostrand formula to $
\frac{\phi(\lambda)-\phi(\lambda')}{\lambda-\lambda'}$ which leads to:
\[
\frac{\phi(\lambda)-\phi(\lambda')}{\lambda-\lambda'} 
=  \frac{-1}{2\pi} \iint_{\Omega_+}\re\left( \frac{\overline{\partial}\Phi(z)}{(\lambda - z)(\lambda' - z)} \right)\diff^{2}z .
\]

\begin{Lem}\label{eq:HS_A_N(f)}
	Suppose $k\geq0$, and $B \subset [-1+\varepsilon,1-\varepsilon]$. Let $\Omega_+$ be as in Lemma \ref{lem:eq}. Let $\phi\in \mcal{C}_c^3(B)$ and let $\Phi$ denote the degree $2$ pseudo-analytic extension. Then, we have
\[
A_N^{(k)}(\phi)   =  -\re\bigg(\int_{\Omega_+}  \overline{\partial}\Phi(z) \, \widetilde{s_N}(z)^2  \frac{\diff^2 z}{2\pi}  \bigg) .
\]
\end{Lem}


Coupling this formula with the controls on the real and imaginary parts of $\widetilde{s_N}(z)$ immediately yields the following.

\begin{Prop}[Bound on the anisotropy]\label{lem:HS_A_N}
Under the same assumptions as in Proposition \ref{lem:HS_L_N}, we have:
		\begin{equation} \label{eq:HS_bound_A_N RMT}
			\|A_N^{(k)}(\phi)\|_q\lesssim\mathtt{c}(\phi)\begin{cases}
				\dfrac{q}{N\alpha \eta}+\dfrac{q^{2}N^{\frac\delta2}}{N^{2}\alpha\eta^{2}}  & \text{if }\alpha \eta \geq N^\delta, \\
				\dfrac{q}{N}+\dfrac{q^{2}}{N^{2}\eta}&\text{otherwise.}
			\end{cases} 
		\end{equation}

\end{Prop}


\begin{Pro}
Let $\Phi$ be a degree 2 pseudo-analytic extension in Proposition \ref{prop:Helffer–Sjöstrand} with cutoff at scale $\delta=\eta$. Via a direct computation, we see that
\begin{equation*}
|\overline{\partial}\Phi(x+\ii y)| \lesssim y^2\chi_\eta(y)|\phi^{(3)}(x)|+\sum_{k=0}^2 |(\chi_\eta)'(y)|y^k|\phi^{(k)}(x)|.
\end{equation*}
From this and Lemma \ref{eq:HS_A_N(f)} we obtain:
\begin{equation*}
\|A_N^{(k)}(\phi)\|_{q}\lesssim \int_{\Omega_+}\left(\chi_\eta(y)y^2|\phi^{3}(x)|+|(\chi_\eta)'(y)|\sum_{k=0}^2y^k|\phi^{(k)}(x)|\right) \cdot\|\widetilde{s_N}(x+\ii y)\|_{q}^2\diff x \diff y
\end{equation*}
with $\Omega_+$ as in Lemma \ref{lem:eq}. 

\textbf{1. Case $ \eta \geq \tfrac{N^{\delta}}{\alpha}$.} Let $ \eta \geq \tfrac{N^{\delta}}{\alpha}$, we make use of Theorem \ref{prop:LL2} up to scale $\sqrt{\frac{\eta}{\alpha}}$ and Theorem \ref{prop:LL1} above. Then,
%
 \begin{align*}
 	\|A_N^{(k)}(\phi)\|_{q}&\lesssim \|\phi^{(3)}\|_{L^1(\R)}\left( \int_0^{\frac{N^{\delta
 			/2}}{\alpha}}y^{2}\left(\dfrac{q}{Ny}+\dfrac{q^{2}}{N^{2}y^{2}}\right)\diff y  +\int_{\frac{N^{\delta/2}}{\alpha}}^{2\eta}y^2\frac{q}{N\alpha y^{2}}\diff y\right)
 	\\&\quad+\sum_{k=0}^2\|\phi^{(k)}\|_{L^1(\R)}\int_{\eta}^{2\eta} y^k|(\chi_\eta)'(y)|\frac{q}{N\alpha y^2}\diff y \\
 	&\lesssim  \mathtt{c}(\phi)\eta ^{-2}\left(\frac{qN^{\delta}}{N\alpha^{2}}+\frac{q^{2}N^{\delta/2}}{N^{2}\alpha}+\dfrac{q\eta}{N\alpha}\right)+\sum_{k=0}^2\mathtt{c}(\phi)\eta^{1-k}\frac{q}{N\alpha } \eta^{k-2}\\
 	&\lesssim \mathtt{c}(\phi)\left(\frac{q}{N\alpha \eta}+\frac{q^{2}N^{\delta/2}}{N^{2}\alpha\eta^{2}}\right),
 \end{align*}
where we neglegted the term $\tfrac{qN^{\delta}}{N\alpha^{2}\eta^{2}}$ because $\eta\geq\tfrac{N^{\delta}}{\alpha}$.

\textbf{2. $ \eta \leq \tfrac{N^{\delta}}{\alpha}$.} Let $\eta\leq\tfrac{N^{\delta}}{\alpha}$, we have:
\begin{align*}
\|A_N^{(k)}(\phi)\|_{q}&\lesssim \|\phi^{(3)}\|_{L^1(\R)} \int_0^{2\eta}y^2\left(\dfrac{q}{Ny}+\frac{q^2}{N^2y^2}\right) \diff y 
\\&\quad+\sum_{k=0}^2\|\phi^{(k)}\|_{L^1(\R)}\int_{\eta}^{2\eta} |(\chi_\eta)'(y)|y^k\left(\frac{q}{Ny}+\dfrac{q^{2}}{N^{2}y^{2}}\right) \diff y \\
&\lesssim \mathtt{c}(\phi)\eta ^{-2}\left(\frac{q\eta^{2}}{N}+\frac{q^{2}\eta}{N^{2}}\right)+\sum_{k=0}^2 \mathtt{c}(\phi)\eta^{1-k}\left(\dfrac{q\eta^{k-1}}{N}+\dfrac{q^{2}\eta^{k-2}}{N^{2}}\right)
\\&\lesssim  \mathtt{c}(\phi)\left(\frac{q}{N}+\frac{q^{2}}{N^{2}\eta}\right).
\end{align*}
This concludes the proof.
\end{Pro}

	\section{Proof of Mesoscopic CLT} \label{sec: proof CLT}

The goal of this section is to show the following mesoscopic CLT with explicit rates.
\begin{Thm}\label{thm:meso sec5}
	Let $\beta_N\gg N^{-1+\kappa}$ for some $\kappa>0$. For all $E\in(-1+\varepsilon,1-\varepsilon)$ for some fixed $\varepsilon>0$ and all $q\geq1$, with $f_\eta\defi f\left(\eta^{-1}(\cdot-E)\right)$, the following holds:
	\begin{enumerate}
		\item \textbf{\textit{(Random Matrix regime)}} if $\eta\gg\tfrac{1}{\alpha}$ then for all $f\in\mcal{C}_c^{5}(\R)$:
		\begin{equation*}
			\mathbf{W}_q\left (\sqrt{N\alpha}L_N^{(k)}(f_\eta), \mcal{N}\left(0,\|f\|_{\msf{H}^{1/2}}^{2}\right)\right ) \lesssim \begin{cases}
			\mfrak{R}_1 & \text{if } \eta \geq \tfrac{N^\delta}{\alpha}, \\\mfrak{R}_2
				&\text{otherwise}.\end{cases}
		\end{equation*}
		where $$	\mfrak{R}_1 \defi\dfrac{q^{2}N^{\frac{\delta}{2}}}{(\alpha\eta)^{2}}\sqrt{\dfrac{\alpha^{3}}{N^{3}}}\log(\eta^{-1})+\dfrac{q^{\frac{3}{2}}N^{\frac{2\delta}{3}}}{ N\alpha\eta^{2}}+	\dfrac{q}{\sqrt{N\alpha}\eta}\log(\eta^{-1})+\sqrt{q}\left (\dfrac{1}{\alpha\eta}+\eta\right )+\eta \log \eta^{-1}\sqrt{\frac{\alpha}{N}}$$
		and
		$$	\mfrak{R}_2 \defi\dfrac{q^{2}}{\alpha\eta}\sqrt{\dfrac{\alpha^{3}}{N^{3}}}\log(\eta^{-1})+	\dfrac{q^{\frac{3}{2}}}{N\eta}+q\sqrt{\dfrac{\alpha}{N}}\log(\eta^{-1})+\sqrt{q}\left(\dfrac{1}{\sqrt{\alpha\eta}}+\eta\right)+\eta \log \eta^{-1}\sqrt{\frac{\alpha}{N}}.$$
		\item \textbf{(\textit{Poisson regime})} if $\tfrac{1}{N}\ll\eta\ll\theta\ll\tfrac{1}{\alpha}$ then for all $f\in\mcal{C}_c^{3}(\R)$:
		$$		\mathbf{W}_q\left(\sqrt{\tfrac{N}{\eta \mu_{\mrm{eq}}(E)}}L_N^{(k)}(f_\eta), \mathcal{N}(0, \|f\|_{L^2(\R)}^2)\right) \lesssim \mfrak{R}_3$$
		where 
		$$\mfrak{R}_3\defi q^{2}\dfrac{\alpha\theta}{(N\eta)^{3/2}}+\dfrac{q^{\frac{3}{2}}\eta}{N}+\dfrac{q}{\sqrt{N\eta}}\left(\alpha\theta+\dfrac{\eta}{\theta}\right)+\sqrt{q}\left(\dfrac{1}{\alpha}+\alpha\theta\right)+\alpha\theta\sqrt{\dfrac{\eta}{N}}.$$
		\item \textbf{(\textit{Critical regime})} If $\alpha\eta\mu_\mrm{eq}(E)\tau= 1$ for some $\tau>0$ then  for all $f\in\mcal{C}_c^{3}(\R)$:
		$$\mbf{W}_q\left (\sqrt{\tfrac{N}{\eta \mu_{\mrm{eq}}(E)}}L_N^{(k)}(f_\eta), \mathcal{N}(0, \Sigma_\tau^2(f))\right )\lesssim\mfrak{R}_4$$
		where 
		$$\mfrak{R}_4\defi\dfrac{q^{2}}{(N\eta)^{3/2}}+\dfrac{q^{\frac{3}{2}}}{N\eta}+\dfrac{q}{\sqrt{N\eta}}+\sqrt{q}\left(\dfrac{1}{N\eta}+\dfrac{1}{\alpha}\right)+\sqrt{\frac{\eta}{N}}.$$
	\end{enumerate}
	The semi-norm $\|\cdot\|_{\msf{H}^{1/2}}$ is given in \eqref{def: H12 norm} and
	\begin{equation*}
		\Sigma_\tau^2(\phi):=\frac{\tau}{2\pi}\int_{\R} \frac{|\xi|\cdot|\widehat{\phi}(\xi)|^2}{\pi+\tau|\xi|}\diff \xi.
	\end{equation*}
\end{Thm}

For the proof, we use Stein's method in the spirit of \cite{LambertLedouxWebb,angst2024sharp} as stated in Proposition \ref{prop:bound Wp}.
	We will use again the following differential operator $\mathcal{L}$ defined in \eqref{eq:def L}.
%
The goal will be to write our linear statistic as an approximate image under $\mathcal{L}$, and apply Proposition \ref{prop:bound Wp}. 
%
%
%
%
To set notation, for a function $f$, we will write $F\defi \int_\R f(x) \diff\bm{\mu}_N(x)$ and $\psi \defi f'$, then \eqref{eq:def L} applied to $F$ yields:
\begin{equation}\label{eq:Stein0}
\dfrac{\mc{L}[ F]}{N\alpha}   = -\bigg(\dfrac{1}{\alpha}-\frac{1}{2N}\bigg)\int_{\R} \psi'(x) \diff\bm{\mu}_N(x) -  \frac{1}2 \iint_{\R^{2}} \mcal{D}[\psi](x,y) \diff\bm{\mu}_N(x)\diff\bm{\mu}_N(y)  + \frac12 \int_{\R} V'(x) \psi(x) \diff\bm{\mu}_N(x)
\end{equation}
where $\mcal{D}$ is defined in Definition \ref{def:ope D}.

We will consider each of the Random matrix regime $(\alpha \eta \rightarrow \infty)$, the Poisson regime $(\alpha \eta \rightarrow 0)$ and the transition regime $(\alpha \eta \asymp 1)$ individually. In this section, we take $\alpha\geq N^{\delta}$ for some $\delta>0$ and assume $k$ is large enough.
\subsection{The Random Matrix regime}

We first tailor \eqref{eq:Stein0} to the Random matrix regime. Recall the operator $\Xi$ given in Definition \ref{def:ope D}.


	\begin{Prop}\label{prop: SteinRMT}
	Let $\varepsilon>0$ be fixed and suppose $B \defi(-1+\varepsilon, 1-\varepsilon)$ and let $\psi \in \mcal{C}_c^3(B)$ satisfy \eqref{eq: Cnorm} with $\eta>0$. Let $f$ be such that $f'=\psi$, and write $F= N\bm{\mu}_N(f_\eta)$. Then for $k$ large enough:
		\begin{equation}\label{Stein1} 
\frac{\mc{L} [F]}{N\alpha}   =   L_N^{(k)}(\Xi[\psi])
+  \boldsymbol{\zeta}(\psi)  +O\bigg(\frac{\eta}N\mathtt{c}(\psi) \bigg)
\end{equation}
where the error term is deterministic and the random part is given by:
\begin{equation}\label{def:zetaRM}
	\boldsymbol{\zeta}(\psi)   \defi   -\bigg(\frac1\alpha-\frac1{2N}\bigg)L_N^{(k)}(\psi') +\sum_{1\le \ell\le k} \re\left(\alpha^{-\ell}\int_{\Omega_+} \overline{\partial}\Psi(z) \sigma_\ell(z) \widetilde{s_N}(z) \frac{\diff^2z}{\pi}\right)  +A_N^{(k)}(\psi) .
\end{equation}
	\end{Prop}

	
	\begin{Pro}
	We expand the terms in \eqref{eq:Stein0}. Using a degree $2$ pseudo-analytic extension $\Psi$ of $\psi$ and applying Proposition \ref{prop:Helffer–Sjöstrand} we find:
	\begin{equation}\label{eq: psiexp}
	-\int_\R \psi'(x) \diff\bm{\mu}_N(x) =  \re\int_{\Omega_+} \overline{\partial}\Psi(z) \partial s_N(z) \frac{\diff^2z}{\pi}  = -L_N^{(k)}(\psi') + \re\int_{\Omega_+} \overline{\partial}\Psi(z)  \partial m_k(z) \frac{\diff^2 z}{\pi} 
\end{equation}
where $m_k$ is the Stieltjes transform of $\mu_k$ defined in Definition \eqref{def:mup}. Similarly, we may expand  
\begin{align*}\label{eq: Vpsiexp}
\nonumber \int_\R V'(x) \psi(x) \diff \bm{\mu}_N( x) &=  L_N^{(k)}(V'\psi)  +\mu_k\left(    \bigg( \re\int_{\Omega_+} \frac{\overline{\partial}\Psi(z)}{\cdot-z} \frac{\diff^2z}{\pi}   \bigg)V'  \right) \\
 &= 
L_N^{(k)}(V'\psi)  +\re \int _{\Omega_+}\overline{\partial}\Psi(z) \mu_k\left(\frac{V'}{\cdot-z}\right) \frac{\diff^2z}{\pi} 
\end{align*}
and with $\mcal{D}$ defined in Definition \ref{def:ope D}:
\begin{align}\label{eq: aniexp}
-\frac{1}2 \iint_{\R^{2}} \mcal{D}[\psi](x,y) &\diff\bm{\mu}_N( x)\diff\bm{\mu}_N( y) \nonumber  =\frac12 \re\int_{\Omega_+} \overline{\partial}\Psi(z) s_N^2(z) \frac{\diff^2z}{\pi}\nonumber
\\&=  \re\int_{\Omega} \overline{\partial}\Psi(z) m_k(z) \widetilde{s_N}(z) \frac{\diff^2z}{\pi}  + \frac12 \re\int_{\Omega_+}\overline{\partial}\Psi(z) m_k^2(z) \frac{\diff^2z}{\pi} + A_N^{(k)}(\psi).
\end{align}
using the formula
\begin{equation*}
	\mcal{D}[g](\lambda,\lambda')=\frac{g(\lambda)-g(\lambda')}{\lambda-\lambda'} 
	=   -\re\iint_{\Omega_+} \frac{\overline{\partial}G (z)}{(\lambda - z)(\lambda' - z)} \frac{\diff^2 z}{\pi}
	\end{equation*}
	for general $g$ with pseudo-analytic extension $G$, Lemma \ref{eq:HS_A_N(f)}. Using this formula again, the first term may be written explicitly, using Definition \ref{def:ope D} as:
\begin{align*}
 \re \int_{\Omega_+} \overline{\partial}\Psi(z) m_k(z) \widetilde{s_N}(z) \frac{\diff^2z}{\pi} 
 &=-L_N^{(k)}\left( \int_{-1}^{1} \mcal{D}[\psi](\cdot,u) \diff\mu_{\mrm{eq}}( u)\right)
\\&\quad+\sum_{\ell=1}^{k}  \alpha^{-\ell}\re\int_{\Omega_+} \overline{\partial}\Psi(z) \sigma_{\ell}(z) \widetilde{s_N}(z) \frac{\diff^2z}{\pi} \\
&= -\frac12 L_N^{(k)}(\psi V') + L_N^{(k)}( \Xi[\psi])
+\sum_{\ell=1}^{k} \alpha^{-\ell}\re\int_{\Omega_+} \overline{\partial}\Psi(z) \sigma_{\ell}(z) \widetilde{s_N}(z) \frac{\diff^2z}{\pi}
\end{align*}
Inserting these into \eqref{eq:Stein0}, we obtain:

\begin{align*}
\frac{\mc{L}[ F]}{N\alpha}   & =L_N^{(k)}(\Xi [\psi])+\bigg( -\bigg(\frac1\alpha-\frac1{2N}\bigg)L_N^{(k)}(\psi')  +\sum_{\ell=1}^{k} \alpha^{-\ell}\re\int_{\Omega_+} \overline{\partial}\Psi(z) \sigma_\ell(z) \widetilde{s_N}(z) \frac{\diff^2z}{\pi}  + A_N^{(k)}(\psi) \bigg) \\
&\quad+\frac12  \re\int_{\Omega_+} \overline{\partial}\Psi(z)\left( \frac2\alpha\partial m_k(z)+m_k^2(z)+\mu_k\bigg( \frac{V'}{\cdot-z}\bigg)\right)  \frac{\diff^2z}{\pi} +O\bigg(\frac{\eta}N \mathtt{c}(\psi)  \bigg)
\end{align*}
where the $O$ error term is deterministic and controlled by
\begin{equation*}
\frac{1}{2N} \bigg| \re\int_{\Omega_+} \overline{\partial}\Psi(z) \partial m_k(z) \frac{\diff^2z}{\pi} \bigg| 
\lesssim  \frac1{N} \int_{\Omega_+} \big|\overline{\partial}\Psi(z) \big| \diff^2z \lesssim  \frac1{2N} \sum_{0\le \ell\le 3 } \eta^{\ell} \|\psi^{(\ell)}\|_{L^1(\R)} 
\lesssim \frac{\mathtt{c}(\psi) \eta}{2N} .
\end{equation*}
By coupling Propositions \ref{prop:mun stieltjes}, the second line is $\mathtt{c}(\psi)\eta O\big(\frac{ 1}{\alpha^{k+1}}+\frac{1}N \big)$ which is negligible if $\alpha^k \ge N$. 
	\end{Pro}
We now make a choice of $\psi$ such that $\Xi [\psi] \approx \phi_\eta$ for a mesoscopic test function $\phi_\eta = \phi\big(\frac{\cdot-E}{\eta}\big)$ where $E \in B$ and $\psi$ satisfies \eqref{eq: Cnorm} with $\eta>0$. The key observation we will use is that 
%
%
on mesoscopic scales,
\[
\Xi [\psi_\eta] \approx -\pi\mu_{\mrm{eq}}(E)\mathcal{H} [\psi_\eta]
\]
where $\mathcal{H}$ is the Hilbert transform, given by 
\begin{equation*}
\mcal{H}[f] : x\in\R \mapsto \dfrac{1}{\pi}\fint_\R \dfrac{f(y)}{y-x}\diff y
\end{equation*}
and $\fint$ denotes the Cauchy principal value integration. 
\begin{Lem}\label{lem:hilbert}
	Let $\phi\in\mcal{C}_c^{5}(\R)$ then: $$\mcal{H}[\phi](x)\underset{|x|\to\infty}{=}-\dfrac{1}{\pi x}\int_\R \phi(t)\diff t+o\left(\dfrac{1}{x}\right) ,\quad\quad\quad\mcal{H}[\phi'](x)\underset{|x|\to\infty}{=}O\left(\dfrac{1}{x^{2}}\right)$$
	and $\mcal{H}[\phi]\in\mcal{C}_{\mrm{loc}}^{4}(\R)$.
\end{Lem}

\begin{Pro} We start with the fact that $\phi\in\mcal{W}^{5,2}(\R)$, thus $\mcal{H}\left[\phi\right]\in\mcal{W}^{5,2}(\R)\subset\mcal{C}^{4}_{\mrm{loc}}(\R)$ and for all $k\in\llbracket0,5\rrbracket$, $\mcal{H}[\phi]^{(k)}=\mcal{H}[\phi^{(k)}]$ (see \cite{Hilberttransforms}). Thus, by \cite[Lemma 2.3]{dworaczek2024clt}, the fact that $\phi\in L^1(\R)$ and setting $g:t\mapsto t\phi(t)$, we have for all $x\in\R$:
	$$x\mcal{H}[\phi](x)+\dfrac{1}{\pi}\int_\R\phi(t)\diff t=\mcal{H}[g](x)$$
	Since $g\in\mcal{W}^{1,2}(\R)$ so does $\mcal{H}[g]$ and since functions in $\mcal{W}^{1,2}(\R)$ go to zero at infinity, we conclude that $\mcal{H}[g](x)=o(1)$ and that $\mcal{H}[\phi](x)=\tfrac{-1}{\pi x}\int_\R\phi(t)\diff t+o(x^{-1})$ at infinity. For the second point, setting $g:t\mapsto t\phi'(t)\in L^{1}(\R)$ and $h:t\mapsto t^{2}\phi'(t)$, we have since for all $x\in\R$ :
	$$x^{2}\mcal{H}[\phi](x)+\dfrac{x}{\pi}\int_\R\phi'(t)\diff t+\dfrac{1}{\pi}\int_\R t\phi'(t)\diff t=x\mcal{H}[g](x)+\dfrac{1}{\pi}\int_\R t\phi'(t)\diff t=\mcal{H}[h](x).$$
	Now since $\int_\R\phi'(t)\diff t=0$ and $h\in \mcal{W}^{1,2}(\R)$, we deduce the second point.
\end{Pro}

\begin{Prop}\label{prop: RMT error control}
	Let $\phi\in\mcal{C}_c^{5}(\R)$ and $\phi_\eta \defi \phi\big(\frac{\cdot-E}{\eta}\big)$ where $E \in B$ and $\eta>0$ satisfies $1\gg\eta \gg\tfrac{1}{\alpha}$. Let $\chi(x)\in \mcal{C}_c^\infty(B)$ satisfy \eqref{eq: Cnorm} with $\eta=1$ such that $\chi =1$ on a neighborhood of $E$. Define 
	\begin{equation*}\label{eq:invert MO}
		\psi(x):=\frac{\chi(x)\mcal{H}[\phi_\eta](x)}{\pi \mu_{\mrm{eq}}(x)}, \qquad x\in\R.
	\end{equation*} 
	Then $ \mathtt{c}(\psi) \lesssim \log(\eta^{-1})$ for $k\geq1$ large enough and all $q\geq1$, with $\boldsymbol{\zeta}$ is as in Proposition \ref{prop: SteinRMT} 
	one has:
	\begin{equation*} \label{err2}
		\sqrt{N\alpha} \cdot \|\boldsymbol{\zeta} \|_q \lesssim  \begin{cases}
			\dfrac{q^{2}N^{\frac{\delta}{2}}}{(\alpha\eta)^{2}}\sqrt{\dfrac{\alpha^{3}}{N^{3}}}\log(\eta^{-1})+ \dfrac{q}{\alpha\eta}\sqrt{\dfrac{\alpha}{N}}\log(\eta^{-1})+\dfrac{\sqrt{q}}{\alpha\eta}  &\quad\quad\quad \text{if }\alpha \eta \geq N^\delta, \\
			\dfrac{q^{2}}{\alpha\eta}\sqrt{\dfrac{\alpha^{3}}{N^{3}}}\log(\eta^{-1})+q\sqrt{\dfrac{\alpha}{N}}\log(\eta^{-1})+\sqrt{\dfrac{q}{\alpha\eta}}&\quad\quad\quad\text{otherwise},
		\end{cases}
	\end{equation*}
	and for $\boldsymbol\xi \defi  L_N^{(k)}(\phi_\eta-\Xi[\psi])$, one has $\sqrt{N\alpha} \cdot \|\boldsymbol{\xi} \|_q \lesssim \eta\sqrt{q}$.
\end{Prop}

\begin{Pro}
	\textbf{1. Estimation of $\bm\zeta$.} We first notice that $\mcal{\psi}\in\mcal{C}^{4}_c(B)$ because of Lemma \ref{lem:hilbert} and use that $|\mathcal{H}[\phi_\eta](x)| \lesssim \frac{C(\phi)}{1+|x-E|/\eta}$ for all $x\in\R$ where $C(\phi)>0$ is a bounded constant depending only on $\phi$, which yields:
	\begin{equation*}
		\|\psi\|_{L^1(\R)} \lesssim \int_\R \chi(x)\cdot |\mathcal{H}[\phi_\eta](x)|\diff x \lesssim \eta \log(\eta^{-1}) .
	\end{equation*}
	Let $\varrho:x\mapsto  \dfrac{\chi(x)}{\pi \mu_\mrm{eq}(x)} \in \mathcal{C}^\infty_c(B)$ for notational ease. Then, $\psi'= \varrho'  \cdot \mathcal{H}[\phi_\eta] + \varrho \cdot \mathcal{H}[\phi_\eta]'.$
	A computation reveals that $\mathcal{H}[\phi_\eta]'=\eta^{-1}\mathcal{H}[\phi']_\eta$. Since $\mathcal{H}[\phi']  \in L^1(\R)$ by Lemma \ref{lem:hilbert} then, $\displaystyle \int_\R \mathcal{H}[\phi']_\eta (x)\diff x \lesssim \eta.$ Thus,
	\[
	\|\psi'\|_{L^1(\R)} \lesssim \eta \int_\R  |\varrho'(x\eta+E)|\cdot| \mathcal{H}[\phi](x)|\diff x +  O(1)\lesssim 1.
	\]
	We can show similarly that for all $k\in\llbracket1,4\rrbracket$:
	\[
	\|\psi^{(k)}\|_{L^1(\R)} \lesssim \eta^{1-k}.
	\]
	Consequently, we have $ \mathtt{c}(\psi) \lesssim \log(\eta^{-1})$ and $ \mathtt{c}(\psi') \lesssim \eta^{-1}$. So, we may make use of the estimates given in Proposition \ref{lem:HS_L_N}, Corollary \ref{lem: genHS_L_N} and Proposition \ref{lem:HS_A_N} to deduce with the constant $C>0$ in Theorem \ref{prop:LL1} and any $q \geq1$:
	\[
	\|L_N^{(k)}(\psi')\|_q\lesssim \left(\sqrt{\dfrac{q}{N\alpha\eta^{2}}}+\dfrac{qN^{\frac{2}{3}\delta}}{N\alpha\eta^{2}}\right)\ind{\alpha\eta\geq N^{\delta}}+\left(\sqrt{\dfrac{q}{N\eta}}+\dfrac{q}{N\eta}\right)\ind{C\leq\alpha\eta\leq N^{\delta}},
	\]
	$$\|A_N^{(k)}(\psi)\|_q \lesssim \log(\eta^{-1})\left(\dfrac{q}{N\alpha\eta}+\dfrac{q^{2}N^{\frac{\delta}{2}}}{N^{2}\alpha\eta^{2}}\right)\ind{\alpha\eta\geq N^{\delta}}+\log(\eta^{-1})\left(\dfrac{q}{N}+\dfrac{q^{2}}{N^{2}\eta}\right)\ind{C\leq\alpha\eta\leq N^{\delta}},$$
	$$\dfrac{1}{\log(\eta^{-1})}\bigg\|\re\int_{\Omega_+} \overline{\partial}\Psi(z) \sigma_{\ell}(z) \widetilde{s_N}(z) \frac{\diff^2z}{\pi} \bigg\|_q  \lesssim   \left(\sqrt{\dfrac{q}{N\alpha}}+\dfrac{qN^{\frac{2}{3}\delta}}{N\alpha\eta}\right)\ind{\alpha\eta\geq N^{\delta}}+\left(\sqrt{\dfrac{q\eta}{N}}+\dfrac{q}{N}\right)\ind{C\leq\alpha\eta\leq N^{\delta}} $$
	Combining all the estimates above with \eqref{def:zetaRM}, we obtain for $\alpha\eta\geq N^{\delta}$:
	\begin{align*} 
		\|\boldsymbol{\zeta} \|_q &\lesssim\left(\dfrac{1}{\alpha\eta}+\dfrac{\log(\eta^{-1})}{\alpha}\right) \left(\sqrt{\dfrac{q}{N\alpha}}+\dfrac{qN^{\frac{2}{3}\delta}}{N\alpha\eta}\right)+\log(\eta^{-1})\left(\dfrac{q}{N\alpha\eta}+\dfrac{q^{2}N^{\frac{\delta}{2}}}{N^{2}\alpha\eta^{2}}\right)
		\\&\lesssim \dfrac{1}{\alpha\eta}\sqrt{\dfrac{q}{N\alpha}}+\dfrac{q}{N\alpha\eta}\log(\eta^{-1})+\dfrac{q^{2}N^{\frac{\delta}{2}}}{N^{2}\alpha\eta^{2}}\log(\eta^{-1})\end{align*}
	while for $1\ll\alpha\eta\leq N^{\delta}$:
	\begin{align*}
		\|\boldsymbol{\zeta} \|_q &\lesssim\left(\dfrac{1}{\alpha\eta}+\dfrac{\log(\eta^{-1})}{\alpha}\right) \left(\sqrt{\dfrac{q\eta}{N}}+\dfrac{q}{N}\right)+\log(\eta^{-1})\left(\dfrac{q}{N}+\dfrac{q^{2}}{N^{2}\eta}\right)
		\\&\lesssim \dfrac{1}{\alpha\eta}\sqrt{\dfrac{q\eta}{N}}+\dfrac{q}{N}\log(\eta^{-1})+\dfrac{q^{2}}{N^{2}\eta}\log(\eta^{-1}).\end{align*}
	We thus obtain:
	\begin{equation*}
		\sqrt{N\alpha}\|\boldsymbol{\zeta} \|_q \lesssim\begin{cases} \dfrac{q^{2}N^{\frac{\delta}{2}}}{(\alpha\eta)^{2}}\sqrt{\dfrac{\alpha^{3}}{N^{3}}}\log(\eta^{-1})+ \dfrac{q}{\alpha\eta}\sqrt{\dfrac{\alpha}{N}}\log(\eta^{-1})+\dfrac{\sqrt{q}}{\alpha\eta} & \text{if } \eta \alpha \ge N^{\delta} \\
			\dfrac{q^{2}}{\alpha\eta}\sqrt{\dfrac{\alpha^{3}}{N^{3}}}\log(\eta^{-1})+q\sqrt{\dfrac{\alpha}{N}}\log(\eta^{-1})+\sqrt{\dfrac{q}{\alpha\eta}} & \text{if } 1\ll\eta\alpha\le N^{\delta}. \\
		\end{cases}
	\end{equation*}
	\textbf{2. Estimation of $\bm\xi$.} Now, with this choice of $\psi$, we have for all $x\in \R$, using that $\mc{H}\left[\mcal{H}[\phi]\right]=-\phi$, \eqref{eq:effpotentialV}, Definition \ref{def:ope D} and the fact that $\psi=0$ outside of a neighborhood of $E$:
	\begin{align*}
		\Xi[\psi] (x) = \fint_{\R} \dfrac{\psi(y)}{x-y} \diff\mu_{\mrm{eq}} ( y)
		&= \frac1\pi \fint_{\R}  \dfrac{\chi(y)}{x-y} \mathcal{H}[\phi_\eta](y) \diff y 
		\\
		& = \chi(x)  \fint_{\R}\dfrac{ \mathcal{H}[\phi_\eta](y)}{x-y}  \frac{\diff y}{\pi}  - \frac1\pi\int_{\R} \dfrac{\chi(x)-\chi(y)}{x-y} \mathcal{H}[\phi_\eta](y) \diff y \\
		&=:  \phi_\eta(x) +\eta\mfrak{R}(x)
	\end{align*}
	where 
	\begin{equation*}
		x\in\R \mapsto \mfrak{R}(x)=\frac{-1}{\pi\eta}\int \dfrac{\chi(x)-\chi(u)}{x-u} \mathcal{H}[\phi_\eta](u) \diff u
	\end{equation*} 
	is a smooth function. Using the Fourier transform, we may compute explicitly:
	\begin{equation*}
		\frac{\chi(x)-\chi(u)}{x-u}   = \int_{\R}\widehat{\chi}(u) \frac{e^{\ii u x}-e^{\ii u y}}{x-y}  \frac{\diff u}{2\pi}=
		\int_\R\int_0^1  \ii u \widehat{\chi}(u) e^{\ii u((1-t) x+ ty)} \diff t \frac{\diff u}{2\pi}
	\end{equation*}
	and for all $t\in(0,1)$, $u\in\R$:
	\begin{equation*}
		\int_\R \mathcal{H}[\phi_\eta](y) e^{\ii y tu}  \diff y 
		= \widehat{ \mathcal{H}[\phi_\eta]}(-u t)  =- \ii \mathrm{sgn}(u) \widehat{\phi_\eta}(-u t)
		= - \eta\ii \mathrm{sgn}(u)  e^{\ii u t E} \widehat{\phi}(-u t \eta). 
	\end{equation*}
	Thus, we may write by Fubini:
	\begin{align*}
		\mfrak{R}(x) &= \frac{-1}{\pi\eta} \int_{\R} \dfrac{\chi(x)-\chi(u)}{x-u} \mathcal{H}[\phi_\eta](u) \diff u 
		\\&=\dfrac{-1}{2\pi^{2}\eta}\int_{0}^{1}\diff t\int_{\R}\diff u \ii u \widehat{\chi}(u) e^{\ii u(1-t) x}\int_{\R}\mathcal{H}[\phi_\eta](y) e^{\ii y tu}  \diff y 
		\\&=  \dfrac{1}{2\pi^{2}}\int_\R |u| \widehat{\chi}(u)  \int_0^1  \widehat{\phi}(-u t \eta) e^{\ii u((1-t) x+ t E)} \diff t \diff u
	\end{align*}
	Differentiating in $x$ and using that $\| \widehat{\phi} \|_{\mathcal{C}^0(\R)} \le \|\phi\|_{L^1(\R)} $ and that $\chi\in\mcal{S}(\R)$, we have for any $p\geq1$ and $x\in\R$ that
	\[
	|\mfrak{R}^{(p)}(x)| \le \|\phi\|_{L^1(\R)}  \int_\R  |u|^{p+1} |\widehat{\chi}(u)| \int_0^1 t^p     \diff t \diff u
	\lesssim C_k.
	\]
	While we cannot use the mesoscopic fluctuation results from Section \ref{sec: concentration} as with the error terms above since $\phi_\eta-\Xi[\psi]$ is not supported inside of the bulk, we may still use the moment estimates  (Proposition \ref{prop: moment estimates}) to obtain
	\[
	\|\bm{\xi}\|_{q}\leq\eta\| L_N^{(k)}(\mfrak{R}) \|_q \lesssim \eta\sqrt{\frac{q}{N\alpha}} 
	\]
	with $\boldsymbol\xi =  L_N^{(k)}(\phi_\eta-\Xi[\psi])$.
\end{Pro}

Finally, we identify the variance term and close the proof of the mesoscopic CLT in the random matrix regime.
\begin{Prop}\label{prop: RMTCLT}
Make the same conventions as in Proposition~\ref{prop: RMT error control}. Then for $k$ large enough and all $q\geq1$ and for $G\sim \mcal{N}\big (0, \|\phi\|_{\msf{H}^{1/2}}^2\big )$, one has:\begin{itemize}
	\item For $ \eta \geq \tfrac{N^\delta}{\alpha}$:
	$$\mathbf{W}_q(\mbf X, G) \lesssim 
		\dfrac{q^{2}N^{\frac{\delta}{2}}}{(\alpha\eta)^{2}}\sqrt{\dfrac{\alpha^{3}}{N^{3}}}\log(\eta^{-1})+
		\dfrac{q^{\frac{3}{2}}N^{\frac{2\delta}{3}}}{ N\alpha\eta^{2}}+	\dfrac{q}{\alpha\eta}\sqrt{\dfrac{\alpha}{N}}\log(\eta^{-1}) +\sqrt{q}\left(\dfrac{1}{\alpha\eta}+\eta\right) +\eta \log \eta^{-1}\sqrt{\frac{\alpha}{N}}.$$
			\item For $\tfrac{1}{\alpha}\ll\alpha \eta \leq \tfrac{N^\delta}{\alpha}$:
			\begin{equation*}
				\mathbf{W}_q(\mbf X, G)\lesssim
				\dfrac{q^{2}}{\alpha\eta}\sqrt{\dfrac{\alpha^{3}}{N^{3}}}\log(\eta^{-1}) +\dfrac{q^{\frac{3}{2}}}{N\eta}+q\sqrt{\dfrac{\alpha}{N}}\log(\eta^{-1})+\sqrt{q}\left(\dfrac{1}{\sqrt{\alpha\eta}}+\eta\right)+\eta \log \eta^{-1}\sqrt{\frac{\alpha}{N}}.
			\end{equation*}
\end{itemize}
\end{Prop}
\begin{Pro}We first begin with the fact that since $\phi\in\mcal{C}_c^{5}(\R)$, $ \|\phi\|_{\msf{H}^{1/2}}<+\infty$. Propositions \ref{prop: SteinRMT} and \ref{prop: RMT error control} show that we may rewrite \eqref{Stein1}, with $\mbf X = \sqrt{N\alpha} L_N^{(k)}(\phi_\eta) $, as:
\[
\sqrt{N\alpha}\frac{\mc{L} [F]}{N\alpha} = \mbf X+ \sqrt{N\alpha}( \boldsymbol{\zeta} +\boldsymbol{\xi}) + O\bigg(\eta \log \eta^{-1}\sqrt{\frac{\alpha}{N}}\bigg)
\]
using that $\mathtt{c}(\psi)\lesssim \log\eta^{-1}$ (as shown in the proof of Proposition \ref{prop: RMT error control}). Applying Proposition \ref{prop:bound Wp} with $G\sim  \mcal{N}\big (0, \|\phi\|_{\msf{H}^{1/2}}^2\big )$ and $\lambda=\sqrt{N\alpha}$, we find: 
\begin{equation*}
\mathbf{W}_q(\mbf X,G) \lesssim
\dfrac{\sqrt{q}}{\sqrt{N\alpha}}\Big\| \nabla \mbf{X} \cdot \nabla F -\sqrt{N\alpha}\|\phi\|_{\msf{H}^{1/2}}^2\Big\|_q + \sqrt{N\alpha}\|\boldsymbol{\zeta}\|_q +\sqrt{N\alpha}\|\boldsymbol{\xi}\|_q+\eta\log\eta^{-1}\sqrt{\dfrac{\alpha}{N}}.
\end{equation*}
\textbf{1. Estimation of quadratic term.} Since $\nabla F = \big(\psi(\lambda_i) \big)_{i=1}^N$ and $\nabla\mbf X =\tfrac{\sqrt{N\alpha}}{N\eta}\big(\,(\phi')_\eta(\lambda_i) \big)_{i=1}^N$, we have $\nabla \mbf{X} \cdot \nabla F=\eta^{-1}\sqrt{N\alpha}\bm{\mu}_N(\psi (\phi')_\eta) $ and thus:
\[
\sqrt{\dfrac{q}{N\alpha}}\Big\| \nabla \mbf{X} \cdot \nabla F -\sqrt{N\alpha}\|\phi\|_{\msf{H}^{1/2}}^2\Big\|_q\leq\eta^{-1}\sqrt{q}\| L_N^{(0)}(\psi (\phi')_\eta)\|_q+\sqrt{q}\left|\|\phi\|_{\msf{H}^{1/2}}^2-\eta^{-1}\mu_{\mrm{eq}}( (\phi')_\eta\psi) \right| 
\]
Since, $\mathtt{c}(\psi (\phi')_\eta)=O(1)$ we have by Proposition \ref{lem:HS_L_N} (and noticing that $\|m_0-m_k\|_{\mcal{C}^{0}(B)}=O(\alpha^{-1})$ for all $q\geq1$),
$$\eta^{-1}\sqrt{q}\|L_N^{(0)}(\psi (\phi')_\eta)\|_{q}
\lesssim\begin{cases}
	\dfrac{q^{\frac{3}{2}}N^{\frac{2\delta}{3}}}{ N\alpha\eta^{2}}+\dfrac{q}{\sqrt{N\alpha}\eta}+\dfrac{\sqrt{q}}{\alpha\eta}  & \quad\quad\text{if }\alpha \eta \geq N^\delta, \\
	\dfrac{q^{\frac{3}{2}}}{N\eta}+\dfrac{q}{\sqrt{N\eta}}+\dfrac{\sqrt{q}}{\alpha\eta}&\quad\quad\text{otherwise.}
\end{cases}$$
For the second term, by Plancherel's formula, one has using that with our conventions $\widehat{\mcal{H}[g]}(x)=\ii\cdot\mrm{sgn}(x)\widehat{f}(x)$:
\begin{align*}
	\eta^{-1}\mu_{\mrm{eq}}(\psi (\phi')_\eta)&=\frac{1}{\pi} \int_{\eta^{-1}(-1+\varepsilon-E)}^{\eta^{-1}(1-\varepsilon-E)} \chi(u\eta+E)\phi'(u)\cdot \mathcal{H}[\phi](u) \diff u 
	\\&=\frac{1}{\pi} \int_{-\infty}^{+\infty} \phi'(u)\cdot \mathcal{H}[\phi](u) \diff u +O(\eta)
	\\ &=  \frac{1}{2\pi^2} \int_\R \ii\xi  \widehat{\phi}(\xi) \cdot \big(-\ii\mathrm{sgn}(\xi)\overline{\widehat{\phi}(\xi)}\big) \diff\xi +O(\eta) \\
	&=\frac{1}{2\pi^2} \int_\R |\xi| |\widehat{\phi}(\xi)|^2 \diff\xi  +O(\eta)
	\\&= \|\phi\|_{\msf{H}^{1/2}}^2 +O(\eta).
\end{align*}
We thus obtain:
\begin{equation*}
	\dfrac{\sqrt{q}}{\sqrt{N\alpha}}\Big\| \nabla \mbf{X} \cdot \nabla F -\sqrt{N\alpha}\|\phi\|_{\msf{H}^{1/2}}^2\Big\|_q\lesssim\begin{cases}
		\dfrac{q^{\frac{3}{2}}N^{\frac{2\delta}{3}}}{ N\alpha\eta^{2}}+\dfrac{q}{\sqrt{N\alpha}\eta} +\sqrt{q}\left(\dfrac{1}{\alpha\eta}+\eta\right) & \quad\quad\text{if }\alpha \eta \geq N^\delta, \\
		\dfrac{q^{\frac{3}{2}}}{N\eta}+\dfrac{q}{\sqrt{N\eta}}+\sqrt{q}\left(\dfrac{1}{\alpha\eta}+\eta\right) &\quad\quad\text{otherwise.}
	\end{cases}
\end{equation*}
\textbf{2. Conclusion.} So, with Proposition \ref{prop: RMT error control}, we obtain using for $\eta\geq\tfrac{N^{\delta}}{\alpha}$:
\begin{equation*}
\mathbf{W}_q(\mbf X, G) \lesssim\dfrac{q^{2}N^{\frac{\delta}{2}}}{(\alpha\eta)^{2}}\sqrt{\dfrac{\alpha^{3}}{N^{3}}}\log(\eta^{-1})+
\dfrac{q^{\frac{3}{2}}N^{\frac{2\delta}{3}}}{ N\alpha\eta^{2}}+	\dfrac{q}{\alpha\eta}\sqrt{\dfrac{\alpha}{N}}\log(\eta^{-1}) +\sqrt{q}\left(\dfrac{1}{\alpha\eta}+\eta\right) +\eta \log \eta^{-1}\sqrt{\frac{\alpha}{N}}
\end{equation*}
while
\begin{equation*}
	 \mathbf{W}_q(\mbf X, G)\lesssim
\dfrac{q^{2}}{\alpha\eta}\sqrt{\dfrac{\alpha^{3}}{N^{3}}}\log(\eta^{-1}) +\dfrac{q^{\frac{3}{2}}}{N\eta}+q\sqrt{\dfrac{\alpha}{N}}\log(\eta^{-1})+\sqrt{q}\left(\dfrac{1}{\sqrt{\alpha\eta}}+\eta\right)+\eta \log \eta^{-1}\sqrt{\frac{\alpha}{N}}
\end{equation*}
which is the desired result.
\end{Pro}

\subsection{The Poisson regime}
We now tailor \eqref{eq:Stein0} to the Poisson regime.
	\begin{Prop}\label{prop: SteinPoisson}
	Let $\varepsilon>0$ be fixed and suppose $B \defi(-1+\varepsilon, 1-\varepsilon)$ and let $\psi \in \mcal{C}_c^3(B)$ satisfy \eqref{eq: Cnorm} with $\eta>0$. Let $f$ such that $f'=\psi$, and write $F= N\bm{\mu}_N(f)$. Then for $k$ large enough:
	\begin{equation}\label{eq:Stein2}
 \dfrac{\mc{L} [F]}{N}   =   -L_N^{(k)}(\psi')
+  \boldsymbol{\zeta}   +O\bigg(\frac{\eta \mathtt{c}(\psi)\alpha}N\bigg)
\end{equation}
where the error term is deterministic and the random part is given by 
\begin{align*}
\boldsymbol{\zeta}  &=   \frac\alpha{2N} L_N^{(k)}(\psi') +  \frac{\alpha}{2} \re\int_{\Omega_+} V'(z) \overline{\partial}\Psi(z) \widetilde{s_N}(z) \frac{\diff^2z}{\pi} +\frac{\alpha}{2}\re\int_{\Omega_+} \overline{\partial}\Psi(z)\Delta_k(z)\frac{\diff^2z}{\pi}\\&\quad+ \alpha \re\int_{\Omega_+} \overline{\partial}\Psi(z) m_k(z) \widetilde{s_N}(z) \frac{\diff^2z}{\pi} +\alpha A_N^{(k)}(\psi).
\end{align*}
	\end{Prop}
	\begin{Pro}
	Starting from \eqref{eq:Stein0}, using \eqref{eq: psiexp}-\eqref{eq: aniexp}
\begin{align*}
\dfrac{\mc{L}[ F]}{N}\   & = -L_N^{(k)}(\psi')+\frac\alpha2\frac2\alpha \re\int_{\Omega_+} \overline{\partial}\Psi(z)\left( \partial m_k(z)+m_k^2(z)+\mu_k\Big(  \frac{V'(\cdot)}{\cdot-z} \Big)\right)  \frac{\diff^2z}{\pi}  
\\
&\quad + \frac\alpha{2N}L_N^{(k)}(\psi') + \frac{\alpha}{2} \re\int_{\Omega_+} V'(z) \overline{\partial}\Psi(z) \widetilde{s_N}(z) \frac{\diff^2z}{\pi} +\frac{\alpha}{2}\re\int_{\Omega_+} \overline{\partial}\Psi(z)\Delta_k(z)\frac{\diff^2z}{\pi}\\\
&\quad+  \alpha\re \int_{\Omega_+} \overline{\partial}\Psi(z) m_k(z) \widetilde{s_N}(z) \frac{\diff^2z}{\pi} + \alpha A_N^{(k)}(\psi)+O\bigg(\frac{ \mathtt{c}(\psi) \eta\alpha}N\bigg)
\end{align*}
where the error term is deterministic and controlled by
\[
\frac{\alpha}{N}\bigg|\re \int_{\Omega_+} \overline{\partial}\Psi(z) \partial m_k(z) \diff^2z \bigg| 
\lesssim  \frac{\alpha}{N} \int_{\Omega_+} \big|\overline{\partial}\Psi(z) \big| \diff^2z \lesssim  \frac{\alpha}{N} \sum_{0\le k\le 3 } \eta^{k} \|\psi^{(k)}\|_{L^1(\R)} 
\lesssim \frac{\mathtt{c}(\psi) \alpha\eta}{N} 
\]
and we have made use of the expansion
\begin{align*}
L_N^{(k)}(V'\psi)&=L_N^{(k)}\left( V'(\cdot)\int \frac{\overline{\partial}\Psi(z)}{\cdot-z}\frac{\diff^2 z}{\pi}\right)\\
&=\re\int_{\Omega_+} V'(z)\overline{\partial}\Psi(z)\widetilde{s_N}(z)\frac{\diff^2z}{\pi}+\re\int_{\Omega_+} \overline{\partial}\Psi(z)\Delta_k(z)\frac{\diff^2 z}{\pi}
\end{align*}
with $\Delta_k$ as in Definition \ref{def: delta}. The second line of our $\mc{L}[F]$ expansion is deterministic and is a $O(\alpha^{-(k+1)})$ by Proposition \ref{prop:mun stieltjes} which is negligible for $k$ large enough.
	\end{Pro}

Unlike in the random matrix regime, the dominant term is $L_N^{(k)}(\psi')$ and it is this term that we seek to cast as the fluctuations of $\phi_\eta$.
\begin{Prop}\label{prop: Poissonerror}
Let $\phi\in \mcal{C}_c^{3}(\R)$, $\tfrac{1}{N}\ll\eta\ll\theta\ll\tfrac{1}{\alpha}$, $\chi\in \mcal{C}_c^\infty(\R)$ such that $\ind{(-1,1)}\leq\chi\leq\ind{(-2,2)}$ and set $\chi_\theta\defi\chi(\tfrac{\cdot-E}{\theta})$ and define
\begin{equation*}\label{eq: defpsiP}
\psi(x) \defi\chi_\theta(x) \int_x^{\infty} \phi_\eta(t)\diff t.
\end{equation*}
Then, $ \mathtt{c}(\psi) \lesssim \theta$, with $\boldsymbol{\zeta}$ is as in Proposition \ref{prop: SteinPoisson} and $\boldsymbol{\xi}\defi  -L_N^{(k)}( \psi' +\phi_\eta )$,
we have for $k$ large enough and all $q\geq1$:
\begin{equation*}
\sqrt{\dfrac{N\eta^{-1}}{\mu_{\mrm{eq}}(E)}}\|\boldsymbol{\zeta} \|_q  \lesssim   \alpha\theta\left(\dfrac{q^{2}}{(N\eta)^{3/2}}+\dfrac{q}{\sqrt{N\eta}}+\sqrt{q}\right) , \qquad \sqrt{\dfrac{N\eta^{-1}}{\mu_{\mrm{eq}}(E)}}\|\boldsymbol{\xi} \|_q  \lesssim\frac{\eta}{\theta} \left(\dfrac{q}{\sqrt{N\eta}}+\sqrt{q} \right).
\end{equation*}
\end{Prop}
\begin{Pro}
We first notice that $\psi\in\mcal{C}_c^{4}(B)$ and that we have for all $x\in\R$:
 \begin{equation*}
	\psi'(x) =-\phi_\eta(x) + \gamma (\chi')_{\theta}(x)\int_{\eta^{-1}(x-E)}^{+\infty}\phi(u)\diff u
\end{equation*}
 with $\gamma \defi \tfrac{\eta}{\theta}\ll1$. Inserting this into our expansion of $\mc{L}[F]$ in Proposition \ref{prop: SteinPoisson}, we find 
 \begin{equation*}
 \dfrac{\mc{L} [F] }{N}  =   L_N^{(k)}(\phi_\eta)
+  \boldsymbol{\zeta} +\boldsymbol{\xi}  +O\bigg(\frac{ \mathtt{c}(\psi)\eta\alpha}N\bigg)
\end{equation*}
with $\boldsymbol{\zeta}$ as in Proposition \ref{prop: SteinPoisson}, and:
\begin{equation*}
\boldsymbol{\xi}\defi  -L_N^{(k)}( \psi' +\phi_\eta )  =  - \gamma L_N^{(k)}\left ( (\chi')_\theta\int_{\eta^{-1}(\cdot-E)}^{+\infty}\phi(u)\diff u\right ).
\end{equation*}
Now, computing directly,
\begin{equation*}
\|\psi\|_{L^1(\R)} \lesssim \|\chi_\theta\|_{L^1(\R)} \|\phi_\eta\|_{L^1(\R)} \lesssim \eta \theta.
\end{equation*}
Using Leibniz differentiation formula, we find for all $p\in\llbracket1,4\rrbracket$:
\begin{align*}
	\psi^{(p)} &= -\sum_{\ell=0}^{p-2}\binom{p-1}{\ell}(\chi^{(\ell)})_\theta\theta^{-\ell}(\eta)^{-(p-1-\ell)}(\phi^{(p-1-\ell)})_\eta  + \gamma \theta^{-(p-1)} (\chi^{(p)})_\theta\int_{\eta^{-1}(\cdot-E)}^{+\infty}\phi(u)\diff u
	\\&\quad+\gamma \sum_{\ell=0}^{p-2}\binom{p-1}{\ell}(\chi^{(\ell+1)})_\theta\theta^{-\ell}(\eta)^{-(p-1-\ell)}(-\phi^{(p-2-\ell)})_\eta
 \end{align*}we have 
\begin{equation*}
\|\psi^{(p)} \|_{L^1(\R)} \lesssim  \theta \eta^{1-p}  + \gamma\theta^{2-p} +\gamma\theta\eta^{1-p}\lesssim  \theta \eta^{1-p},
 \end{equation*}
 which allows to deduce that $ \mathtt{c}(\psi) \lesssim \theta$ and  $\mathtt{c}(\psi') \lesssim \theta$. 
Thus, applying Propositions \ref{lem:HS_L_N} and \ref{lem:HS_A_N} we find, noting that $\mathtt{c}((\chi')_\theta)\lesssim1$:
$$ \|  L_N^{(k)}\left( (\chi')_\theta \right) \|_q\lesssim\sqrt{\dfrac{q\theta}{N}}+\dfrac{q}{N},\quad\quad\quad\| L_N^{(k)}(\psi') \|_q \lesssim\theta\sqrt{\dfrac{q\eta}{N}}+\dfrac{q\theta}{N},\quad\quad \| A_N^{(k)}(\psi)\|_q \lesssim\theta
\left(	\dfrac{q}{N}+\dfrac{q^{2}}{N^{2}\eta}\right).$$
Furthermore, by Corollary \ref{lem: genHS_L_N},
 \[
  \bigg\| \re \int_{\Omega_+} V'(z) \overline{\partial}\Psi(z) \widetilde{s_N}(z) \diff^2z \bigg\|_q +\bigg\| \re \int_{\Omega_+} m_k(z) \overline{\partial}\Psi(z) \widetilde{s_N}(z) \diff^2z \bigg\|_q\lesssim \theta  \left(	\sqrt{\dfrac{q\eta}{N}}+\dfrac{q}{N}\right) ,
 \]
 and by Corollary \ref{cor: glocon},
 \begin{equation*}
 \left\|\re\int_{\Omega_+} \overline{\partial}\Psi(z)\Delta_k(z)\frac{\diff^2z}{\pi}\right\|_q \lesssim \sqrt{\frac{q}{N\alpha}}\int_{\Omega_+}  |\overline{\partial}\Psi(z)|\diff^2 z\lesssim \sqrt{\frac{q}{N\alpha}}\sum_{\ell=0 } ^{3}\eta^{\ell} \|\psi^{(\ell)}\|_{L^1(\R)} 
\lesssim \frac{\theta\eta\sqrt{q}}{\sqrt{N\alpha}}
 \end{equation*}
Recalling the definitions of $\boldsymbol{\zeta}$ and $\boldsymbol{\xi}$, this implies that if $\eta \ge N^{-1}$, $\alpha \ll N$ and $\eta \ll \theta\ll \alpha^{-1}$, the anisotropy is negligible and
\[
\sqrt{\dfrac{N\eta^{-1}}{\mu_{\mrm{eq}}(E)}}\|\boldsymbol{\zeta} \|_q \lesssim    \alpha\theta\left(\sqrt{q}+\dfrac{q}{\sqrt{N\eta}}+\dfrac{q^{2}}{(N\eta)^{3/2}}\right), \qquad \sqrt{\dfrac{N\eta^{-1}}{\mu_{\mrm{eq}}(E)}}\|\boldsymbol{\xi} \|_q  \lesssim \gamma \left(\sqrt{q}+\dfrac{q}{\sqrt{N\eta}} \right)   .
\]
This concludes the proof.
\end{Pro}

We now recover the variance term and close the CLT in the Poisson regime.
\begin{Prop}\label{prop: PoissonCLT}
Take the same conventions as in Proposition \ref{prop: Poissonerror}. Then for $G\sim\mathcal{N}\left(0, \|\phi\|_{L^2(\R)}^2\right)$, $k$ large enough and any $q\geq1$:
	\begin{equation*}
\mathbf{W}_q\left(\sqrt{\tfrac{N}{\eta \mu_{\mrm{eq}}(E)}}L_N^{(k)}(\phi_\eta), G\right) \lesssim q^{2}\dfrac{\alpha\theta}{(N\eta)^{3/2}}+\dfrac{q^{\frac{3}{2}}\eta}{N}+\dfrac{q}{\sqrt{N\eta}}\left(\alpha\theta+\dfrac{\eta}{\theta}\right)+\sqrt{q}\left(\dfrac{1}{\alpha}+\alpha\theta\right)+\alpha\theta\sqrt{\dfrac{\eta}{N}}.
	\end{equation*}
 \end{Prop}
\begin{Pro}
Let $\mbf X = \sqrt{\frac{N}{\eta \mu_{\mrm{eq}}(E)}} L_N^{(k)}(\phi_\eta)$, using \eqref{eq:Stein2} and the control on $\boldsymbol{\zeta}$ and $\boldsymbol{\xi}$ from Proposition \ref{prop: Poissonerror} we obtain:
\[
\mbf X = \sqrt{\frac{N}{\eta \mu_{\mrm{eq}}(E)}}\left(\dfrac{\mc{L} [F]}{N} -  \boldsymbol{\zeta}-\boldsymbol{\xi}\right)  +O\bigg(\alpha\theta\sqrt{\dfrac{\eta}{N}}\bigg).
\]
Applying Proposition \ref{prop:bound Wp} with $G \sim  \mcal{N}(0,\|\phi\|_{L^{2}(\R)}^2)$, we find, with $\Gamma = \nabla \mbf X\cdot \nabla F$
\begin{equation}\label{eq:WpPoissonproof}
	\mbf{W}_q(\mbf X, G) \lesssim \sqrt{q}\Big\| \tfrac\Gamma{\sqrt{N\eta\mu_{\mrm{eq}}(E)}} -\|\phi\|_{L^{2}(\R)}^2\Big\|_q +\sqrt{\tfrac{N}{\eta\mu_{\mrm{eq}}(E)}}\left(\|\bm{\zeta}\|_q+\|\bm\xi\|_q\right)+\alpha\theta\sqrt{\dfrac{\eta}{N}}.
\end{equation}
Since $\nabla F = \big(\psi(\lambda_i)\big)_{i=1}^N$ and $\nabla\mbf X = \tfrac{\eta^{-1}}{\sqrt{N\eta\mu_{\mrm{eq}}(E)}}\big( (\phi')_\eta(\lambda_i)\big)_{i=1}^N$, we have:
\[
\dfrac{1}{\sqrt{N\eta\mu_{\mrm{eq}}(E)}}\Gamma = \frac{N\eta^{-1}\bm{\mu}_N(\psi (\phi')_\eta) }{N\eta\mu_{\mrm{eq}}(E)} =  \frac{\mu_{\mrm{eq}}(\psi (\phi')_\eta)}{\mu_{\mrm{eq}}(E)\eta^2}+\frac{(\mu_k-\mu_{\mrm{eq}})(\psi (\phi')_\eta)}{\mu_{\mrm{eq}}(E)\eta^2}
+\frac{L_N^{(k)}(\psi (\phi')_\eta)}{\mu_{\mrm{eq}}(E)\eta^2}.
\]
The fluctuation term is an error term; by definition of $\psi$ we have $\psi (\phi')_\eta  =\chi_\theta(\phi')_\eta  \int_\cdot^{\infty} \phi_\eta(t)\diff t$, so that :
\[
\|\psi (\phi')_\eta\|_{L^1(\R)} \lesssim \|(\phi')_\eta\|_{L^1(\R)} \|\phi_\eta\|_{L^1(\R)} \lesssim \eta^2 , \qquad
\|(\psi (\phi')_\eta)^{(k)}\|_{L^1(\R)} \lesssim  \eta^{2-k} \text{ for $k\in\llbracket1,3\rrbracket$}.
\]
Hence, $\mathtt{c}(\psi (\phi')_\eta) \lesssim \eta$ and we may conclude from Proposition \ref{lem:HS_L_N} that 
\[
\sqrt{q}\|L_N^{(k)}(\psi (\phi')_\eta)\|_q \lesssim  \eta q\sqrt{\dfrac{\eta}{N}}+\dfrac{q^{\frac{3}{2}}\eta}{N}. 
\]
Furthermore, we find 
\begin{equation*}
\left|\frac{(\mu_k-\mu_{\mrm{eq}})(\psi (\phi')_\eta)}{\mu_{\mrm{eq}}(E)\eta^2}\right|\lesssim \frac{1}{\alpha\eta^2}\|\psi (\phi')_\eta\|_{L^1(\R)}\lesssim \frac{1}{\alpha}.
\end{equation*}
Using changes of variables, the fact that $\chi(x)=1$ on $(-1,1)$, that there exists $M>0$ such that supp $\phi\subset(-M,M)$ and an integration by part, we have:
\begin{align*}
\mu_{\mrm{eq}}(\psi (\phi')_\eta)  &= \eta\int_{E-2\theta}^{E+2\theta} (\phi')_\eta(x)  \int_{\eta^{-1}(x-E)}^{+\infty}\phi(t)\diff t\mu_{\mrm{eq}}(x)\diff x
\\&=\eta^{2}\int_{-M}^{M}\mu_{\mrm{eq}}(u\eta+E) \phi'(u)\int_{u}^{+\infty}\phi(t)\diff t \diff  u
\\&=\eta^{2}\int_{-M}^{M}\Big(\mu_{\mrm{eq}}(E)+\int_{E}^{E+u\eta}\rho_{\mrm{eq}}'(t)\diff t\Big) \phi'(u)\int_{u}^{+\infty}\phi(t)\diff t \diff  u
\\&=\eta^{2}\mu_{\mrm{eq}}(E)\int_{\R} \phi'(u)\int_{u}^{+\infty}\phi(t)\diff t \diff  u+O(\eta^{3})
\\&=  \eta^{2} \, \mu_{\mrm{eq}}(E)\|\phi\|_{L^{2}(\R)}^{2} +O(\eta^{3}).
\end{align*}
We thus conclude that:
$$\sqrt{q}\Big\|\dfrac{1}{\sqrt{N\eta\mu_{\mrm{eq}}(E)}}\Gamma- \|\phi\|_{L^{2}(\R)}^{2}\Big\|_{q}\lesssim \dfrac{q^{\frac{3}{2}}\eta}{N}+ \eta q\sqrt{\dfrac{\eta}{N}}+\sqrt{q}\left(\dfrac{1}{\alpha}+\eta\right) ,
$$
and thus by Proposition \ref{prop: Poissonerror} and \eqref{eq:WpPoissonproof}:
$$ \mbf{W}_q(\mbf X, G) \lesssim q^{2}\dfrac{\alpha\theta}{(N\eta)^{3/2}}+\dfrac{q^{\frac{3}{2}}\eta}{N}+\dfrac{q}{\sqrt{N\eta}}\left(\alpha\theta+\dfrac{\eta}{\theta}\right)+\sqrt{q}\left(\dfrac{1}{\alpha}+\alpha\theta\right)+\alpha\theta\sqrt{\dfrac{\eta}{N}}.$$\qedsymbol{symbol}
\end{Pro}

\subsection{The critical regime}
In this section, we restrict ourselves to the regime $\eta \tau \mu_{\mrm{eq}}(E)=\alpha^{-1}$ for some fixed $\tau>0$. As we will see, an interesting competition emerges between the behavior in the random matrix and Poisson regimes on the level of the variance. The formula we will derive for the variance interpolates between the $\msf{H}^{1/2}$ semi-norm that is characteristic of the Random matrix regime, and the $L^2$ norm that is characteristic of the Poisson regime. We first have the following slight modification of Proposition \ref{prop: SteinRMT}, which is obtained merely by regrouping some error terms which in this regime are no longer negligible.
\begin{Prop}\label{prop: Steintrans}	Let $\varepsilon>0$ be fixed and suppose $B \defi(-1+\varepsilon, 1-\varepsilon)$ and let $\psi \in \mcal{C}_c^3(B)$ satisfy \eqref{eq: Cnorm} with $\eta>0$. Let $f$ such that $f'=\psi$, and write $F= N\bm{\mu}_N(f)$. Then for $k$ large enough:
		\begin{equation}\label{eq:Stein3}
\frac{\mc{L} [F]}{N\alpha}   =   L_N^{(k)}\bigg(\Xi[\psi]-\frac{\psi'}{\alpha}\bigg)
+  \boldsymbol{\zeta}(\psi)  +O\bigg(\frac{\eta}N \mathtt{c}(\psi) \bigg)
\end{equation}
where the error term is deterministic and the random part is given by 
\begin{equation}\label{eq:zetaCrit}
	\boldsymbol{\zeta}(\psi)   \defi  \frac1{2N}L_N^{(k)}(\psi') +\sum_{\ell=1}^{k} \re\left(\alpha^{-\ell}\int_{\Omega_+} \overline{\partial}\Psi(z) \sigma_\ell(z) \widetilde{s_N}(z) \frac{\diff^2z}{\pi}\right)  +A_N^{(k)}(\psi) .
\end{equation}
\end{Prop}
The goal is now to invert the operator $\psi \mapsto \Xi [\psi]-\frac{\psi'}{\alpha}$, so we may write the dominant term as $L_N^{(k)}(\phi_\eta)$. To motivate our choice, notice that for a mesoscopic test function $\psi$ we may write
\[
\Xi[\psi](x) = \int_{\R} \frac{\psi(y)}{x-y} \diff\mu_{\mrm{eq}}( y) \simeq \mu_{\mrm{eq}}(E)  \int_{\R} \frac{\psi(y)}{x-y}\diff y  =- \pi \mu_{\mrm{eq}}(E) \mathcal{H}[\psi](x)
\]
so that, if we are looking to invert only approximately, it suffices to solve $  -\pi \mathcal{H}[\varphi_\eta]-\eta\tau (\varphi_\eta)'  = \phi_\eta $, where we recall that $\eta \tau \mu_{\mrm{eq}}(E)=\alpha^{-1}$. Now, in Fourier space one has 
\begin{equation*}
-\ii(\tau\eta\xi+\pi\mathrm{sgn}(\xi))\widehat{\varphi_\eta}(\xi) =\widehat{\phi_\eta}(\xi) 
\end{equation*}
which rescaling yields
\begin{equation*}
-\ii(\tau\eta\xi+\pi\mathrm{sgn}(\xi))\eta \widehat{\varphi}(\eta \xi)e^{-\ii E \xi}=\eta  \widehat{\phi}(\eta\xi) e^{-\ii E \xi}.
\end{equation*}
Solving for $\varphi$, this yields 
\begin{equation*}
\widehat{\varphi}(\eta \xi)  = \frac{\ii  \mathrm{sgn}(\xi)\widehat{\phi}(\eta \xi)}{\pi+\tau|\eta\xi|}.
\end{equation*}
Thus, if we define
\begin{equation}\label{eq: defvarphi}
\varphi(x) \defi \ii  \int_\R e^{\ii x\xi} \frac{\mathrm{sgn}(\xi)\widehat{\phi}(\xi)} {\pi+\tau|\xi|}  \diff \xi = -2  \im\bigg( \int_{\R_+}e^{\ii x\xi} \frac{\widehat{\phi}(\xi)}{\pi+\tau|\xi|} \diff \xi\bigg) 
\end{equation}
then
\[
-\pi\mathcal{H}[\varphi_\eta] -\tau (\varphi')_\eta  = \phi_\eta.
\]

\begin{Prop}\label{prop: transerror}
Let $\phi\in \mcal{C}_c^3(\R)$. Assume $\eta \tau \mu_{\mrm{eq}}(E)=\alpha^{-1}$ for some fixed $\tau>0$. Define
\begin{equation}\label{eq: deftranspsi}
\psi \defi  \frac{\chi \cdot  \varphi_\eta}{\mu_{\mrm{eq}}}
\end{equation}
for some $\chi\in\mcal{C}^{\infty}_c(\R)$ and $\varphi$ as in \eqref{eq: defvarphi}. Then, with $\boldsymbol{\zeta}$ is as in Proposition \ref{prop: Steintrans} and $\boldsymbol{\xi}$ defined in \eqref{eq: defxitrans}, we have:
\begin{equation*}
\frac{\mc{L} [F]}{N \alpha}   =   L_N^{(k)}(\phi_\eta)
-  \boldsymbol{\zeta} -\boldsymbol{\xi}  +O\bigg(\frac{\eta}N\bigg)
\end{equation*}
with
\begin{equation*}
\sqrt{\dfrac{N\eta^{-1}}{\mu_{\mrm{eq}}(E)}}\|\boldsymbol{\zeta}\|_q \lesssim\dfrac{q^{2}}{(N\eta)^{3/2}}+\dfrac{q}{\sqrt{N\eta}}+\sqrt{q}\left(\dfrac{\eta^{-1}}{N}+\dfrac{1}{\alpha}\right), \qquad \sqrt{\dfrac{N\eta^{-1}}{\mu_{\mrm{eq}}(E)}}\|\boldsymbol{\xi}\|_q\lesssim \dfrac{q}{\alpha\sqrt{N\eta}}+\dfrac{\sqrt{q}}{\alpha}.
\end{equation*}
\end{Prop}
\begin{Pro}
%
%
%
\textbf{1. Approximate inversion.}
Using our definition of  $\psi$, we have for all $x\in\R$, with $\varrho = \chi/\mu_{\mrm{eq}}$ :
\[
\Xi[\psi](x)    - \alpha^{-1}\psi'(x) =  \fint_\R \frac{\chi(y)\varphi_\eta(y)}{x-y} \diff y - \frac{\tau\mu_{\mrm{eq}}(E)}{\mu_{\mrm{eq}}(x)} \chi(x)  (\varphi')_\eta(x) - \frac{\varrho'(x)\varphi_\eta(x)}{\alpha},
\]
with $\Xi$ as in Definition \ref{def:ope D}, and hence:
\begin{align*}
\Xi[\psi](x)    -\dfrac{\psi'(x)}{\alpha} &= \pi \chi(x) \fint_{\R} \frac{\varphi_\eta(y)}{x-y} \frac{\diff y}\pi - \int_{\R} \frac{\chi(x)-\chi(y)}{x-y} \varphi_\eta(y) \diff y - \tau  \chi(x) (\varphi')_\eta(x)
\\&\quad+\bigg(1- \frac{\mu_{\mrm{eq}}(E)}{\mu_{\mrm{eq}}(x)} \bigg) \tau\chi(x) (\varphi')_\eta(x) - \frac{\varrho'(x)\varphi_\eta(x)}{\alpha} \\
&= \phi_\eta(x)- \int_\R \mcal{D}[\chi](x,y) \varphi_\eta(y) \diff y +\bigg(1- \frac{\mu_{\mrm{eq}}(E)}{\mu_{\mrm{eq}}(x)} \bigg) \tau\chi(x) (\varphi')_\eta(x) - \frac{\varrho'(x)\varphi_\eta(x)}{\alpha}
\end{align*}
which can be written as
\begin{equation*}
\Xi[\psi](x)    - \dfrac{\psi'(x)}{\alpha}   =\phi_\eta(x) - \eta  \Delta (x) - \Upsilon(x) 
\end{equation*}
where we have defined
\begin{equation}\label{eq: def Delta Upsilon}
	\Delta(x) \defi \frac1\eta  \int_\R \frac{\chi(x)-\chi(y)}{x-y} \varphi_\eta(u) \diff u , \qquad 
	\Upsilon(x) \defi  -\bigg(1- \frac{\mu_{\mrm{eq}}(E)}{\mu_{\mrm{eq}}(x)} \bigg) \tau\chi(x) (\varphi')_\eta(x) +\frac{\varrho'(x)\varphi_\eta(x)}{\alpha}.
\end{equation}
Inserting these into \eqref{eq:Stein3} yields:
\begin{equation}\label{eq:Stein3new}
\frac{\mc{L}[ F]}{N\alpha}   =   L_N^{(k)}(\phi_\eta)
-  \boldsymbol{\zeta}-\boldsymbol{\xi}   +O\bigg(\frac{\eta}N\mathtt{c}(\psi) \bigg)
\end{equation}
where $\boldsymbol{\zeta}$ is as in Proposition \ref{prop: Steintrans} and
\begin{equation}\label{eq: defxitrans}
\boldsymbol{\xi}   \defi \eta L_N^{(k)}(\Delta) +L_N^{(k)}(\Upsilon)
\end{equation}
\textbf{2. Control of $\boldsymbol{\zeta}$ and $\boldsymbol{\xi}$.}
First, as a rescaled test function, $\mathtt{c}(\psi)\lesssim 1$ and $\mathtt{c}(\psi')\lesssim \eta^{-1}$. Hence, using Proposition \ref{lem:HS_L_N}, Corollary \ref{lem: genHS_L_N} and Proposition \ref{lem:HS_A_N} we deduce
\[
\bigg\|\re\int_{\Omega_+} \overline{\partial}\Psi(z) \sigma_{\ell}(z) \widetilde{s_N}(z) \frac{\diff^2z}{\pi} \bigg\|_q +\eta\left\|L_N^{(k)}(\psi')\right\|_q \lesssim\sqrt{\dfrac{q\eta}{N}}+\dfrac{q}{N} , \qquad\qquad 
\|A_N^{(k)}(\psi)\|_q \lesssim 	\dfrac{q}{N}+\dfrac{q^{2}}{N^{2}\eta}.
\]
Thus by \eqref{eq:zetaCrit}:
\begin{align*}
\sqrt{\dfrac{N\eta^{-1}}{\mu_{\mrm{eq}}(E)}}\|\boldsymbol{\zeta}\|_q &\lesssim\sqrt{N\eta^{-1}}\left(\dfrac{\eta^{-1}}{N}+\dfrac{1}{\alpha}\right) \left(\sqrt{\dfrac{q\eta}{N}}+\dfrac{q}{N}\right) +	\sqrt{N\eta^{-1}}\left(\dfrac{q}{N}+\dfrac{q^{2}}{N^{2}\eta}\right) 
\\&\lesssim\dfrac{q^{2}}{(N\eta)^{3/2}}+\dfrac{q}{\sqrt{N\eta}}+\sqrt{q}\left(\dfrac{\eta^{-1}}{N}+\dfrac{1}{\alpha}\right).
\end{align*}
We can control $\Delta$ using a Fourier trick as we did in the proof of Proposition \ref{prop: RMT error control}. Indeed,  for all $x,y\in\R$
\begin{equation*}
\frac{\chi(x)-\chi(y)}{x-y}   = \int_{\R}\widehat{\chi}(\xi) \frac{e^{\ii\xi x}-e^{\ii\xi y}}{x-y}  \frac{\diff\xi}{2\pi}=
\int_\R\int_0^1  \ii \xi \widehat{\chi}(\xi) e^{\ii\xi((1-t) x+ ty)} \diff t \frac{\diff\xi}{2\pi}
\end{equation*}
and thus integrating over $y$:
\begin{align*}
\Delta(x) = \frac{1}{\eta} \int_\R \dfrac{\chi(x)-\chi(y)}{x-y} \varphi_\eta(y) \diff y
&= \frac{\ii}{\eta} \int_\R\int_0^1   \xi \widehat{\chi}(\xi) \widehat{\varphi_\eta}(-t\xi) e^{\ii\xi(1-t) x} \diff t \frac{\diff\xi}{2\pi}\\
&= \ii \int_\R \int_0^1   \xi \widehat{\chi}(\xi) \widehat{\varphi}(-t\xi \eta)e^{\ii t\xi E} e^{\ii\xi(1-t) x} \diff t \frac{\diff\xi}{2\pi}
\end{align*}
Now, from the definition of $\varphi$ \eqref{eq: defvarphi}, we have $\widehat{\varphi}(\xi)  = \frac{\ii  \mathrm{sgn}(\xi)\widehat{\phi}(\xi)}{\pi+\tau|\xi|}$. Thus, we obtain:
\[
\Delta(x) = - \int_\R \int_0^1 \frac{|\xi|}{\pi+\tau|\xi|t\eta} \widehat{\phi}(-t\eta\xi)  \widehat{\chi}(\xi) e^{\ii t\xi E} e^{\ii\xi(1-t) x} \diff t \frac{\diff\xi}{2\pi} . 
\]
Now, using $\| \widehat{\phi} \|_{\mathcal{C}^0(\R)} \le \|\phi\|_{L^1(\R)}$ and the fact that $\Delta\in\mcal{C}^{\infty}(\R)$, so for all $k\geq1$:
\[
|\Delta^{(k)}(x)| \le \|\phi\|_{L^1(\R)}  \int_\R |\xi|^{k+1}  |\widehat{\chi}(\xi)| \int_0^1 t^k   \diff t \diff\xi
\lesssim C_k.
\]
Hence, even though $\Delta$ has support that extends outside of the bulk, we may apply the moment estimates from Proposition \ref{prop: moment estimates} and \ref{prop: moment estimates infinity}, to find for $k$ large enough:
\begin{equation}\label{eq: q moment Delta}
\sqrt{\dfrac{N\eta^{-1}}{\mu_{\mrm{eq}}(E)}}\|\eta L_N^{(k)}(\Delta)\|_{q}\lesssim \sqrt{\frac{q\eta}{\alpha }}.
\end{equation}
Finally, we analyze $\Upsilon$. The fluctuations of the term $\frac{\varrho'\varphi_\eta}{\alpha}$ is exactly amenable to Proposition \ref{lem:HS_L_N}; as a rescaled test function, $\mathtt{c}(\varphi_\eta)\lesssim 1$ and so 
\begin{equation}\label{eq: q moment varrho}
\sqrt{\dfrac{N\eta^{-1}}{\mu_{\mrm{eq}}(E)}}\left\|L_N^{(k)}\left(\frac{\varrho'\varphi_\eta}{\alpha}\right)\right\|_q \lesssim \dfrac{q}{\alpha\sqrt{N\eta}}+\dfrac{\sqrt{q}}{\alpha}.
\end{equation}
Next, we compute $\mathtt{c}\left(\bigg(1- \frac{\mu_{\mrm{eq}}(E)}{\mu_{\mrm{eq}}} \bigg) \tau\chi (\varphi')_\eta\right)\lesssim \eta$. The function $\varphi\in\mcal{C}^{4}(\R)$ and using mean-value theorem:
\begin{align*}
\int_B \left|\bigg(1- \frac{\mu_{\mrm{eq}}(E)}{\mu_{\mrm{eq}}(x)} \bigg) \tau\chi(x) (\varphi')_\eta(x)\right|\diff x&\lesssim \int_B |\mu_{\mrm{eq}}(x)-\mu_{\mrm{eq}}(E)|(\varphi')_\eta(x)|\diff x
\\&\lesssim  \int_B|x-E|\cdot |(\varphi')_\eta(x)|\diff x
\\&\lesssim \eta^2\int_\R |u|\cdot|\phi'(u)|\diff u\lesssim\eta^{2}.
\end{align*}
Furthermore, we also have for all $k\in\llbracket1,4\rrbracket$ by Leibniz formula:
\begin{align*}
\left(\Big(1- \frac{\mu_{\mrm{eq}}(E)}{\mu_{\mrm{eq}}} \Big) \tau\chi (\varphi')_\eta\right)^{(k)}&=\Big(1- \frac{\mu_{\mrm{eq}}(E)}{\mu_{\mrm{eq}}} \Big) \tau\chi \big(\varphi^{(k+1)}\big)_\eta\eta^{-k}
\\&\quad+\sum_{\ell=1}^k c_\ell\Big(\chi^{(\ell)}-\mu_{\mrm{eq}}(E)\varrho^{(\ell)}\Big)\left(\varphi^{(k+1-\ell)}\right)_\eta \eta^{\ell-k}
\end{align*}
for some constants $c_k$, from which we may deduce for all $k\in\llbracket1,4\rrbracket$ that, as before, by the mean-value theorem:
\begin{equation*}
\left\|\left(\bigg(1- \frac{\mu_{\mrm{eq}}(E)}{\mu_{\mrm{eq}}} \bigg) \tau\chi (\varphi')_\eta\right)^{(k)}\right\|_{L^1(\R)}\lesssim \eta^{2-k}+\sum_{\ell=1}^k\eta^{1+\ell-k} \lesssim \eta^{2-k}.
\end{equation*}
Hence, $\mathtt{c}\left(\bigg(1- \frac{\mu_{\mrm{eq}}(E)}{\mu_{\mrm{eq}}} \bigg) \tau\chi (\varphi')_\eta\right)\lesssim \eta$ and by Proposition \ref{lem:HS_L_N} we find:
\begin{equation*}
\sqrt{\dfrac{N\eta^{-1}}{\mu_{\mrm{eq}}(E)}}\left\|L_N^{(k)}\left(\bigg(1- \frac{\mu_{\mrm{eq}}(E)}{\mu_{\mrm{eq}}} \bigg) \tau\chi (\varphi')_\eta\right)\right\|_{q}\lesssim q\sqrt{\dfrac{\eta}{N}}+\sqrt{q}\eta.
\end{equation*}
So, by the previous estimate, \eqref{eq: def Delta Upsilon} and \eqref{eq: q moment varrho}, we obtain using that $\alpha\eta\asymp1$:
\begin{equation*}
\sqrt{\dfrac{N\eta^{-1}}{\mu_{\mrm{eq}}(E)}}\|L_N^{(k)}(\Upsilon)\|_q\lesssim \dfrac{q}{\alpha\sqrt{N\eta}}+\dfrac{\sqrt{q}}{\alpha}
\end{equation*}
and thus by the previous estimate, \eqref{eq: defxitrans} and \eqref{eq: q moment Delta} that: 
\begin{equation*}
\sqrt{\dfrac{N\eta^{-1}}{\mu_{\mrm{eq}}(E)}}\|\boldsymbol{\xi}\|_q\lesssim \dfrac{q}{\alpha\sqrt{N\eta}}+\dfrac{\sqrt{q}}{\alpha}.
\end{equation*}
This concludes the proof.
\end{Pro}
We now identify the variance term and close the proof of the CLT in the critical regime.
\begin{Prop}
Make the same conventions as in Proposition \ref{prop: transerror}. Then, setting 
\begin{equation}\label{def: transvar}
\Sigma_\tau^2(\phi):=\frac{\tau}{2\pi}\int_{\R} \frac{|\xi||\widehat{\phi}(\xi)|^2}{\pi+\tau|\xi|}\diff \xi
\end{equation}
we have for $k$ large enough and all $q\geq1$:
$$\mbf{W}_q(\sqrt{\frac{N\eta^{-1}}{ \mu_{\mrm{eq}}(E)}}L_N^{(k)}(\phi_\eta), \mcal{N}(0, \Sigma_\tau^2(\phi)))\lesssim\dfrac{q^{2}}{(N\eta)^{3/2}}+\dfrac{q^{\frac{3}{2}}}{N\eta}+\dfrac{q}{\sqrt{N\eta}}+\sqrt{q}\left(\dfrac{\eta^{-1}}{N}+\dfrac{1}{\alpha}\right)+\sqrt{\frac{\eta}{N}}.$$
\end{Prop}

\begin{Pro}We set $\mbf X \defi  \sqrt{\frac{N\eta^{-1}}{ \mu_{\mrm{eq}}(E)}}L_N^{(k)}(\phi_\eta)$ and use Proposition \ref{prop: transerror} to obtain:
\[
\mbf X =  \sqrt{\frac{N\eta^{-1}}{ \mu_{\mrm{eq}}(E)}}\left(\mc{L} \left[\frac{F}{N\alpha}\right] +  \boldsymbol{\zeta} +\boldsymbol{\xi}\right)  +O\bigg(\sqrt{\frac{\eta}{N}}\bigg)
\]
Applying Proposition \ref{prop:bound Wp} with $G \sim  \mcal{N}(0, \Sigma_\tau^2(\phi))$, we find using $\alpha\eta\tau\mu_{\mrm{eq}}(E)=1$:
\[
\mbf{W}_q(\mbf X, G) \lesssim  \sqrt{\dfrac{q}{N\alpha}}\| W\|_{q}+\sqrt{N\eta^{-1}}\left(\|\bm\zeta\|_{q}+\|\bm\xi\|_{q}\right) +\sqrt{\dfrac{\eta}{N}}, \qquad  W \defi \nabla \mbf X \cdot \nabla F - \sqrt{\dfrac{N\alpha}{\tau}} \Sigma_\tau^2(\phi).
\]
Since $\nabla F = \big(\psi(\lambda_i) \big)_{i=1}^N$ and $\nabla\mbf X = \eta^{-1}\tfrac{1}{\sqrt{N\eta\mu_{\mrm{eq}}(E)}}\big( (\phi')_\eta(\lambda_i)\big)_{i=1}^N= \eta^{-1}\sqrt{\tfrac{\alpha\tau}{N}}\big( (\phi')_\eta(\lambda_i)\big)_{i=1}^N$, we have:
\[
\dfrac{1}{\sqrt{\tau N\alpha}}W + \dfrac{1}{\tau}\Sigma_\tau^2(\phi)= \eta^{-1}\bm{\mu}_N(\psi (\phi')_\eta) =   \eta^{-1}\left(\mu_{\mrm{eq}}(\psi (\phi')_\eta)+(\mu_k-\mu_{\mrm{eq}})(\psi (\phi')_\eta)
+L_N^{(k)}(\psi (\phi')_\eta)\right) 
\]
Since $\psi$ and $(\phi')_\eta$ are both rescaled at scale $\eta$, $\mathtt{c}(\psi (\phi')_\eta)\lesssim 1$ and we may conclude from Proposition \ref{lem:HS_L_N} that:
\[
\frac{\sqrt{q}}{\eta}\|L_N^{(k)}(\psi (\phi')_\eta)\|_q \lesssim\dfrac{q^{\frac{3}{2}}}{N\eta}+\dfrac{q}{\sqrt{N\eta}}. 
\]
Furthermore, since $\|\mu_k-\mu_{\mrm{eq}}\|\lesssim \alpha^{-1}$,
\begin{equation*}
\left|\frac{\tau (\mu_k-\mu_{\mrm{eq}})(\psi (\phi')_\eta)}{\eta}\right|\lesssim \frac{\tau}{\alpha \eta}\|\psi(\phi')_\eta\|_{L^1(\R)}\lesssim \eta.
\end{equation*}
Using the definition of $\psi$ in \eqref{eq: deftranspsi} along with Plancherel and \eqref{eq: defvarphi}, we have:  
\begin{align*}
\frac{\tau \mu_{\mrm{eq}}(\psi (\phi')_\eta)}{\eta}=\frac{\tau}{\eta}\int_{\R} \varphi_\eta(x)(\phi')_\eta(x)\diff x=\tau \int_{\R} \varphi(x)\phi'(x)\diff x&=\frac{\tau}{2\pi}\int_{\R} \widehat{\varphi}(\xi)\overline{\ii\xi\widehat{\phi}(\xi)}\diff \xi \\
&=\frac{\tau}{2\pi}\int_{\R} \frac{\ii \mathrm{sgn}(\xi)\widehat{\phi}(\xi)}{\pi+\tau|\xi|}\overline{\ii \xi \widehat{\phi}(\xi)}\diff \xi \\
&=\frac{\tau}{2\pi}\int_{\R} \frac{|\xi|\cdot|\widehat{\phi}(\xi)|^2}{\pi+\tau|\xi|}\diff \xi=\Sigma_\tau^2(\phi)
\end{align*}
Thus,
\begin{equation*}
 \sqrt{\dfrac{q}{N\alpha}}\| W\|_{q}\lesssim\dfrac{q^{\frac{3}{2}}}{N\eta}+\dfrac{q}{\sqrt{N\eta}}+ \sqrt{q}\eta.
\end{equation*}
Finally using Proposition \ref{prop: transerror}, we obtain:
$$\mbf{W}_q(\mbf X, G)\lesssim\dfrac{q^{2}}{(N\eta)^{3/2}}+\dfrac{q^{\frac{3}{2}}}{N\eta}+\dfrac{q}{\sqrt{N\eta}}+\sqrt{q}\left(\dfrac{\eta^{-1}}{N}+\dfrac{1}{\alpha}\right)+\sqrt{\frac{\eta}{N}}$$
which yields the conclusion.
\end{Pro}	
		\begin{Rem}
			Notice that this variance interpolates between the Random matrix and Poisson regimes. For the sake of this comparison, define:
			\begin{equation*}
				\epsilon^{-1}:=\sqrt{N\min\left(\tau \alpha, \frac{1}{ \mu_{\mrm{eq}}(E)\eta}\right)},
			\end{equation*}
			If we take $\tau\downarrow 0$, then $\alpha \eta \rightarrow \infty$ (Random matrix regime) and the normalization of $L_N^{(k)}$ in the above theorem is given by $\sqrt{N\alpha \tau}$. Renormalizing by $\frac{1}{\sqrt{\tau}}$ to agree with the Random matrix regime, the variance becomes:
			\begin{equation*}
				\frac{1}{2\pi}\int_{\R}\frac{|\xi|\cdot|\widehat{\phi}(\xi)|^2}{\pi+\tau|\xi|}\diff \xi\underset{\tau \downarrow 0}{\longrightarrow} \frac{1}{2\pi^2}\int_{\R} |\xi|\cdot|\widehat{\phi}(\xi)|^2\diff\xi=\|\phi\|_{\msf{H}^{1/2}}^2
			\end{equation*}
			Analogously, if we take $\tau \rightarrow \infty$ then $\alpha \eta \rightarrow 0$ (Poisson regime), and the normalization of $L_N^{(k)}$ in the above theorem is the correct $\sqrt{\frac{\mu_{\mrm{eq}}(E)\eta}{N}}$. The variance then becomes 
			\begin{equation*}
				\frac{\tau}{2\pi}\int_{\R}\frac{|\xi|\cdot|\widehat{\phi}(\xi)|^2}{\pi+\tau|\xi|}\diff \xi\underset{\tau \uparrow \infty}{\longrightarrow} \frac{1}{2\pi}\int_{\R} |\widehat{\phi}(\xi)|^2\diff\xi=\|\phi\|_{L^2(\R)}^2
			\end{equation*}
			which agrees with the variance we found in the Poisson regime.
		\end{Rem}

		\appendix{
			\section{A priori bound and loop equations}\label{app:apriori}
			﻿\subsection{A priori bound}\label{subsec2:apriori}
			Before that we need, some objects defined first in \cite{maida2014free}.
			\begin{Def}\label{def:spacedconfig} Let $\bm{\lambda}\in\R^{N}$ such that $\lambda_1<\dots<\lambda_N$, we define the configuration $\widetilde{\bm{\lambda}}$ by:
				$$\widetilde{\lambda}_{i+1}=\widetilde{\lambda}_{i}+\max\left(\lambda_{i+1}-\lambda_i,\dfrac{1}{N^{3}}\right),$$
				and denote $\widetilde{\bm{\mu}}_N\defi\dfrac{1}{N}\sum_{i=1}^{N}\delta_{\widetilde{\lambda}_i}$. We also define
				$\widetilde{\bm{\mu}}_{N,u}\defi\widetilde{\bm{\mu}}_N\star\mcal{U}_N$
				where $\mcal{U}_N$ is the uniform distribution on $(0,N^{-5})$.
				For an arbitrary vector $\bm{\lambda}\in\R^{N}$, we first order by applying a permutation $\sigma$. We then obtain the same objects $\widetilde{\bm{\lambda}}$, $\widetilde{\bm{\mu}}_N$ and $\widetilde{\bm{\mu}}_{N,u}$ and can apply $\sigma^{-1}$.
			\end{Def}
			
			\begin{Def}[Distance]
				We define the following squared-distance $\mfrak{D}$ for all measures $\mu,\mu'\in\mcal{M}_1(\R)$ by:
				\begin{align}
					\label{def:distanceMesures}
					\mfrak{D}(\mu,\mu') &\defi \left( -\int_{\R^2} \log|x-y|\diff(\mu-\mu')(x)\diff(\mu-\mu')(y) \right)^{1/2}\\
					&= \left( \int_0^{+\infty} \frac{1}{t}\big|\widehat{(\mu-\mu')}(t) \big|^2 \diff t \right)^{1/2}.\nonumber
				\end{align}
			\end{Def}
			﻿		\begin{Prop}[Concentration inequality]
				\label{thm:concentration}
				There exists $K\in \R$ depending on $V$, such that for any $N\geq 1$, $\bm{\lambda}\in\R^{N}$, $\eta\in(0,1)$,
				\begin{equation}
					\label{ineq:concentration}
					\log p_N (\bm{\lambda})\leq -\dfrac{N\alpha}{2}\mfrak{D}^{2}[\widetilde{\bm{\mu}}_{N,u},\mu_{\mrm{eq}}]-\dfrac{N\alpha}{2}(1-\eta)\int_\R V_{\mrm{eff}}(x)\diff \widetilde{\bm{\mu}}_{N,u}(x)+CN(\log N)^{2}.
				\end{equation}
			\end{Prop}

	\begin{Pro}
	We first show that 
	\begin{equation}\label{eq:lowerboundZN}
		\log\mcal{Z}_N[V]\geq-\dfrac{N\alpha}{2} \mcal{E}(\mu _{\mrm{eq}})+O\left(N\log N\right). 
	\end{equation} Since $\rho _{\mrm{eq}}$ is bounded by a constant $M$, we have:
	$$|\gamma_i-\gamma_j|\geq\dfrac{M}{N}$$
	where the $\gamma_i$'s are the classical positions, uniquely defined for all $i\in\llbracket1,N\rrbracket$ by:
	$$\int_{-1}^{\gamma_i}\diff\mu _{\mrm{eq}}(t)=\dfrac{i}{N}.$$
	Restricting the space of integrations to
	$$\Omega=\left\{\bm{\lambda}\in\R^{N},\,\max\limits_{i}|\lambda_i-\gamma_i|\leq\dfrac{1}{4MN}\right\},$$
	noticing that for all $\bm{\lambda}\in\Omega$, $i\in\llbracket1,N\rrbracket$,
	$$-V(\lambda_i)\geq-V(\gamma_i)-\dfrac{\|V'\|_{\infty,[-2,2]}}{4MN}$$
	and using the change of variables $\lambda_i\rightarrow u_i+\gamma_i$, we obtain:
	\begin{multline*}
		\mcal{Z}_N[V]\geq e^{-\alpha\frac{\|V'\|_{\infty,[-2,2]}}{8M}}\cdot\prod_{i=1}^{N}e^{-\frac{\alpha}{2} V(\gamma_i)}\cdot\int_{[\pm\frac{1}{4MN}]^{N}}\prod_{i<j}^{N}\left|u_i-u_j+\gamma_i-\gamma_j\right|^{\beta}\cdot \diff^{N}\bm{u}
		\\\geq e^{-\alpha\frac{\|V'\|_{\infty,[-2,2]}}{8M}}\cdot\prod_{i=1}^{N}e^{-\frac{\alpha}{2} V(\gamma_i)}\cdot\prod_{i<j}^{N}\left(\gamma_j-\gamma_i\right)^{\beta}\int_{v_1<\dots<v_N}\prod_{i=1}^{N}\ind{|u_i|\leq\frac{1}{4MN}}\cdot \diff^{N}\bm{v}
		\\\geq e^{-\alpha\frac{\|V'\|_{\infty,[-2,2]}}{8M}}\cdot\prod_{i=1}^{N}e^{-\frac{\alpha}{2} V(\gamma_i)}\cdot\prod_{i<j}^{N}\left(\gamma_j-\gamma_i\right)^{\beta}\dfrac{1}{N!}\left(\dfrac{1}{2MN}\right)^{N}.
	\end{multline*}
	We know use, that
	$$\dfrac{1}{N}\sum_{i=1}^{N}V(\gamma_i)=\int_{-1}^{1}V(x)\diff\mu _{\mrm{eq}}(x)+O(1/N)$$
	and a similar argument, see \cite[Proof Theorem 4.4]{guionnet2019asymptotics},
	$$\sum_{i<j}^{N}\log(\gamma_j-\gamma_i)=\dfrac{N^{2}}{2}\int_{-1}^{1}\log|x-y|\diff\mu _{\mrm{eq}}(x)\diff\mu _{\mrm{eq}}(y)+O(N).$$
	This concludes the proof for \eqref{eq:lowerboundZN}.
	
	We now prove \eqref{ineq:concentration}. Using \eqref{eq:lowerboundZN} and the spaced configuration $\bm{\lambda}$ defined in Definition \ref{def:spacedconfig}, we obtain the existance of a constant $C>0$ such that:
	\begin{equation}
		\label{ineq:boundDensity}
		\log p_N (\bm{\lambda})\leq -\dfrac{N\alpha}{2}\mathcal{E}(\mu_{\mrm{eq}})+CN\log N+\frac{\beta}{2}\sum_{i\neq j}\log|\widetilde{\lambda}_i-\widetilde{\lambda}_j|-\dfrac{\alpha}{2} \sum_{i=1}^N V(\lambda_i)\, . 
	\end{equation}
	We now show the following estimate, for a constant $C>0$:
	\begin{equation}
		\label{ineq:estimeeEnergie}
		\sum_{i\neq j}\log|y_i-y_j|\leq2+N^2\iint_{\R^2} \log|x-y|\diff\widetilde{\bm{\mu}}_N(x)\diff\widetilde{\bm{\mu}}_N(y) + CN\log N\, .
	\end{equation}
	Let $i< j$ and $u,v\in [0,N^{-5}]$. Since for $x\neq 0$ and $|h|\leq \frac{|x|}{2}$, we have $$\big|\log|x+h| - \log|x|\big|\leq \frac{2|h|}{|x|},$$
	we deduce that:
	$$\big|\log|\widetilde{\lambda}_i-\widetilde{\lambda}_j+u- v|-\log|y_i-y_j|\big| \leq \frac{2|u-v|}{|\widetilde{\lambda}_i-\widetilde{\lambda}_j|} \leq \frac{2N^{-5}}{N^{-3}} =\dfrac{2}{N^{2}}.$$
	Thus, summing over $i\neq j$ and integrating with respect to $u$ and $v$, we get
	\begin{align*}
		\sum_{i\neq j}&\log|\widetilde{\lambda}_i-\widetilde{\lambda}_j|\leq 2+ N^{10}\sum_{i\neq j}\iint_{\R^2}\log|\widetilde{\lambda}_i-\widetilde{\lambda}_j+u-v| \ind{0<u,v<N^{-5}}\diff u\diff v \\
		&= 2+ N^2\iint_{\R^2} \log|x-y|\diff \widetilde{\bm{\mu}}_{N,u}(x)\diff\widetilde{\bm{\mu}}_{N,u}(y) - N^{11}\iint_{\R^2} \log|u-v|\ind{0<u,v<N^{-5}}\diff u\diff v\,
		\\&\leq2+ N^2\iint_{\R^2} \log|x-y|\diff \widetilde{\bm{\mu}}_{N,u}(x)\diff\widetilde{\bm{\mu}}_{N,u}(y)+CN\log N
	\end{align*}
	for a constant $C>0$ independent of $\bm{\lambda}$, which allows us to deduce \eqref{ineq:estimeeEnergie}. We can now combine \eqref{ineq:boundDensity} and \eqref{ineq:estimeeEnergie}, we get:
	\begin{equation*}
		\log p_N (\bm{\lambda})\leq \dfrac{N\alpha}{2}\mathcal{E}(\mu_{\mrm{eq}})-\dfrac{N\alpha}{2}\mathcal{E}(\widetilde{\bm{\mu}}_{N,u}) + \dfrac{N\alpha}{2}\int_\R V(x)\diff\widetilde{\bm{\mu}}_{N,u}(x) -\dfrac{\alpha}{2} \sum_{i=1}^N V(\lambda_i)+CN\log N.
	\end{equation*}
	By using,
	$$\mathcal{E}(\widetilde{\bm{\mu}}_{N,u})=\mathcal{E}(\mu_{\mrm{eq}})+\mfrak{D}^{2}[\widetilde{\bm{\mu}}_{N,u},\mu_{\mrm{eq}}]+\int_\R V_{\mrm{eff}}(x)\diff \widetilde{\bm{\mu}}_{N,u}(x),$$
	we obtain, since there exists $\varepsilon,C>0$ such that for any $|t|<\varepsilon$, $|V'(x+t)|\leq C(|V(x)|+1)$
	\begin{multline*}
		\Big|\int_\R V(x)\diff(\widetilde{\bm{\mu}}_{N,u}-\bm{\mu}_N)(x)\Big|\leq\dfrac{1}{N^{4}}\sum_{i=1}^{N}(i-1)\max_{i}\left\{V'(\lambda_i+|t|), t\in[0,\tfrac{(i-1)}{N^{3}}]\right\}
		\\\leq\dfrac{C}{N^{2}}\left(\int_\R |V(x)|\diff x+1\right), 
	\end{multline*}
	we obtain:
		\begin{equation*}
		\log p_N (\bm{\lambda})\leq \dfrac{-N\alpha}{2}\mfrak{D}^{2}[\widetilde{\bm{\mu}}_{N,u},\mu_{\mrm{eq}}]-\dfrac{N\alpha}{2}(1-\dfrac{C}{N^{2}})\int_\R V_{\mrm{eff}}(x)\diff\widetilde{\bm{\mu}}_{N,u}(x)+CN\log N.
	\end{equation*}
	 This concludes the proof.
	\end{Pro}

			\begin{Lem}[Truncation lemma]\label{lemma:trunc}There exists $\mfrak{R}_N^{(n)}$ defined on 
			$\mcal{A}_{\kappa,0}$ such that for all $f_n\in\mcal{A}_{\kappa,0}:$
				$$\braket{f_n}=\braket{f_{n|c}}+\mfrak{R}_N^{(n)}[f_n],\hspace{1cm}\left|\mfrak{R}_N^{(n)}[f_n]\right|\leq c_n\cdot e^{-c\alpha}\|f\|_{\kappa,0} $$
				where $f_{n|c}(x_1,\dots,x_n)\defi f_{n}(x_1,\dots,x_n)\cdot\prod_{i=1}^{n}\phi(x_i)$ and $\phi\in\mcal{C}_c^{\infty}(\R)$ such that:
				$$\mbf{1}_{[-1-\varepsilon,1+\varepsilon]}(x)\leq\phi(x)\leq\mbf{1}_{[-2,2]}(x)$$
				for $\varepsilon<1$.
			\end{Lem}
			
			\begin{Pro}[of Lemma \ref{lemma:trunc}] For simplicity, we denote $f_n$ (resp. $f_{n|c}$) by $f$  (resp. $f_{|c}$) and because:
			$$\braket{f}=\braket{f_{\mrm{sym}}},$$
			where $f_{\mrm{sym}}(\bm{x})\defi\dfrac{1}{n!}\sum_{\sigma\in\mfrak{S}_n}f\left(x_{\sigma(1)},\dots,x_{\sigma(n)}\right)$, we notice that we can assume to $f$ to be symmetric function. Now,
			$$\braket{f}=\braket{f_{|c}}+\mfrak{A}_N(f)$$
			where $\mfrak{A}_N(f)\defi\sum_{i=1}^{n}\binom{n}{i}\mfrak{A}_{N,i}(f)$ and
			$$\mfrak{A}_{N,i}(f)\defi\mbb{E}_N\left[\int_{\left([-1-\varepsilon,1+\varepsilon]^{c}\right)^{i}}\prod_{\ell=1}^{i}(1-\phi(x_\ell))\diff\bm{\mu}_N(x_\ell)\int_{[-2,2]^{n-i}}f(\bm{x})\cdot\prod_{\ell=i+1}^{n}\phi(x_\ell)\diff L_N^{(0)}(x_\ell)\right].$$
			Thus a straightforward bound gives:
			$$|\mfrak{A}_{N,i}(f)|\leq4\cdot2^{n-i}\cdot\|f\|_{\kappa,0}\cdot\|e^{\kappa V}\|_{\infty,[-2,2]}^{n-i}\cdot A_i$$
			where $A_i=\braket{\otimes_{l=1}^{i}(e^{\kappa V}\ind{[-1-\varepsilon,1+\varepsilon]^c})}$. Now using the exact same proof as in \cite[Lemma 3.1.9]{borot2016asymptotic} together with \eqref{ineq:concentration}, we obtain that $\mfrak{A}_N(f)=O\left(\|f\|_{\kappa,0} e^{-c\alpha}\right)$ for some universal constant $c>0$.
		\end{Pro}
		
			\begin{Pro}[of Proposition \ref{a priori bound}]
			The proof can be found in \cite[Corollary 3.1.10]{borot2016asymptotic} in a similar setting by considering the event $\Omega_{M,N}\defi\left\{\bm{\lambda}\in\R^{N},\quad\mfrak{D}^{2}(\widetilde{\bm{\mu}}_{N,u},\mu _{\mrm{eq}})>\dfrac{M\log N}{\alpha}\right\}$ for $M>0$ large enough. This event satisfies: $\mbb{P}_N\left(\Omega_{M,N}\right)=O\left(e^{-cN\log N}\right) $ for some $c>0$ and a $M$ large enough by Proposition \ref{thm:concentration}.
		\end{Pro}
		
		\subsection{Loop equations}
	\label{subsec:app loop}
		
			\begin{Pro}[of Proposition \ref{prop:DSequations}]Let $\big(\psi^{(i)}\big)_{a=1}^{n+1}\in \mc{C}_c^{\infty}(\R)^{n+1}$ and set $\phi^{(1)}\defi\Xi^{-1}[\phi^{(1)}]$
			, set $$V_{\bm{\varepsilon}}(\lambda)\defi V(\lambda)+\displaystyle\sum_{i=2}^{n+1}\epsilon_i\left(\psi^{(i)}(\lambda)-\int_{-1}^{1}\psi^{(i)}(y)\diff\mu_{\mrm{eq}}(y)\right) $$ and then define
			\begin{equation*}
				p_N^{(\bm{\varepsilon})}(\bm{\lambda})\defi \dfrac{1}{\mc{Z}_N[V_{\bm{\varepsilon}}]}\prod_{i<j}^N|\lambda_i-\lambda_j|^{\beta }\cdot\prod_{i=1}^{N}e^{-\frac{\alpha}{2} V_{\bm{\varepsilon}}(\lambda_i)}
			\end{equation*}
			The new normalisation constant $\mc{Z}_N[V_{\bm{\varepsilon}}]$ is such that $p_N^{(\bm{\varepsilon})}$ is still a pdf on $\R^N$. We then define $G_t(\mu)=\mu+t\phi^{(1)}(\mu)$. Since $\partial_1\phi^{(1)}$ is bounded, for $t$ small enough $G_t$ is a diffeomorphism over $\R$. By a change of variable $\lambda_i=G_t(\mu_i)$, we obtain :
			\begin{equation*}
				1=\int_{\R^N}p_N^{(\bm{\varepsilon})}(\lambda_1,\dots,\lambda_N)\diff^N\bm{\lambda}=\int_{\R^N}p_N^{(\bm{\varepsilon})}\Big(G_t(\mu_1),\dots,G_t(\mu_N)\Big)\prod_{i=1}^NG_t'(\mu_i)\diff\mu_i.
			\end{equation*}
			Making an asymptotic expansion up to the first order in $t$ of the right-hand side of this equality yields:
			\begin{multline*}
				1=\int_{\R^N}\diff^N\bm{\lambda}\cdot p_N^{(\bm{\varepsilon})}(\bm{\lambda})\cdot\Bigg\{1+t\sum_{i=1}^{N}\partial_1\phi^{(1)}(\lambda_i)\Bigg\}\cdot\Bigg\{1+t\dfrac{\beta}{2}\sum_{i\neq j}\dfrac{\phi^{(1)}(\lambda_i)-\phi^{(1)}(\lambda_j)}{\lambda_i-\lambda_j}\Bigg\}\\\cdot\Bigg\{1-t\frac{\alpha}{2} \sum_{i=1}^NV_{\bm{\varepsilon}}'(\lambda_i)\phi^{(1)}(\lambda_i)\Bigg\}+O(t^2).
			\end{multline*}
			Identifying the terms linear in $t$ leads to, with $\braket{g_n}_{\otimes_{\ell=1}^{n}\mu_\ell}\defi\E_N^{(\bm{\varepsilon})}\left[\otimes_{\ell=1}^{n}\mu_\ell(g_n)\right] $	where $\E_N^{(\bm{\varepsilon})}$ is the expectation value with density $	p_N^{(\bm{\varepsilon})}$:	\begin{equation*}
				0=-\dfrac{N\alpha}{2}\Braket{V_{\bm{\varepsilon}}'\phi^{(1)}}_{\bm{\mu}_N }^{(\bm{\varepsilon})}+\dfrac{N\alpha}{2}\Braket{\mcal{D}[\phi^{(1)}]}_{\bm{\mu}_N \otimes \bm{\mu}_N }^{(\bm{\varepsilon})}+(N-\dfrac{\alpha}{2} )\Braket{\partial_1\phi^{(1)}}_{\bm{\mu}_N }^{(\bm{\varepsilon})}.
			\end{equation*}
			By definition of $V_{\bm{\varepsilon}}$, we get:
			\begin{equation*}
				\dfrac{1}{2}\Braket{V'\phi^{(1)}}_{\bm{\mu}_N }^{(\bm{\varepsilon})}+\sum_{i=2}^{n+1}\dfrac{\varepsilon_i}{2}\Braket{\phi^{(1)}\partial_1\psi^{(i)}}_{\bm{\mu}_N }^{(\bm{\varepsilon})}-\dfrac{1}{2}\Braket{\mcal{D}[\phi^{(1)}]}_{\bm{\mu}_N \otimes \bm{\mu}_N }^{(\bm{\varepsilon})}-\left (\dfrac{1}{\alpha }-\dfrac{1}{2N}\right )\Braket{\partial_1\phi^{(1)}}_{\bm{\mu}_N }^{(\bm{\varepsilon})}=0.
			\end{equation*}
			It becomes, after recentring the empirical measures against $\mu _{\mrm{eq}}$ we obtain:
			\begin{multline*}
			\dfrac{1}{2}	\Braket{V'\phi^{(1)}}_{ L_N}^{(\bm{\varepsilon})}+\dfrac{1}{2}\Braket{V'\phi^{(1)}}_{\mu _{\mrm{eq}}}-\dfrac{1}{2}\Braket{\mcal{D} [\phi^{(1)} ]}_{ L_N\otimes  L_N}^{(\bm{\varepsilon})}-\Braket{\mcal{D} [\phi^{(1)} ]}_{ L_N\otimes \mu _{\mrm{eq}}}^{(\bm{\varepsilon})}-\dfrac{1}{2}\Braket{\mcal{D} [\phi^{(1)} ]}_{\mu _{\mrm{eq}}\otimes \mu _{\mrm{eq}}}
				\\-\left (\dfrac{1}{\alpha }-\dfrac{1}{2N}\right )\left (\Braket{\partial_1\phi^{(1)}}_{ L_N}^{(\bm{\varepsilon})}+\Braket{\partial_1\phi^{(1)}}_{\mu _{\mrm{eq}}}\right)+\sum_{i=2}^{n+1}\varepsilon_i\Big(\Braket{\phi^{(1)}\partial_1\psi^{(i)}}_{ L_N}^{(\bm{\varepsilon})}+\Braket{\phi^{(1)}\partial_1\psi^{(i)}}_{\mu _{\mrm{eq}}}\Big)=0.
			\end{multline*}
			Using \eqref{eq:contrainst mu infinity} and Definition \ref{def:ope D}, we get:
			\begin{multline}\label{DSbuilding}
				-\Braket{\Xi\left[\phi^{(1)}\right] }_{ L_N}^{(\bm{\varepsilon})}=\dfrac{-1}{2}\Braket{\mcal{D}[\phi^{(1)}]}_{ L_N\otimes  L_N}^{(\bm{\varepsilon})}
				-\left (\dfrac{1}{\alpha }-\dfrac{1}{2N}\right )\left (\Braket{\partial_1\phi^{(1)}}_{ L_N}^{(\bm{\varepsilon})}+\Braket{\partial_1\phi^{(1)}}_{\mu _{\mrm{eq}}}\right)
				\\+\sum_{i=2}^{n+1}\varepsilon_i\Big(\Braket{\phi^{(1)}\partial_1\psi^{(i)}}_{ L_N}^{(\bm{\varepsilon})}+\Braket{\phi^{(1)}\partial_1\psi^{(i)}}_{\mu _{\mrm{eq}}}\Big).
			\end{multline}
			Taking $\varepsilon_i$ go to 0 for all $i$, $\phi^{(1)}=\Xi^{-1}[\psi]$ leads to:
			\begin{equation*}
				\Braket{\psi}_{ L_N}=\left (\dfrac{1}{\alpha }-\dfrac{1}{2N}\right )\Braket{\Xi^{-1}[\psi]'}_{\mu _{\mrm{eq}}}
				+\left (\dfrac{1}{\alpha }-\dfrac{1}{2N}\right)\Braket{\Xi^{-1}[\psi]'}_{ L_N}+\dfrac{1}{2}\Braket{\mcal{D} \Xi^{-1}[\psi]}_{ L_N\otimes  L_N}.
			\end{equation*}
			Hence we obtain \eqref{DSlvl1}.
			
			Furthermore, in \eqref{DSbuilding} multiplying by $\dfrac{\mc{Z}_N[V_{\bm{\varepsilon}}]}{\mcal{Z}_N[V]}$, applying $\partial_{\varepsilon_2}\dots\partial_{\varepsilon_{n+1}}$ and evaluate at $\varepsilon_i=0$ leads to
			\begin{multline}\label{DSbuildinglvl}
				\Braket{\Xi[\phi^{(1)}] (\xi_1).\prod_{i=2}^{n+1}\phi^{(i)}(\xi_i)}_{\overset{n+1}{\bigotimes} L_N}=\left (\dfrac{1}{\alpha }-\dfrac{1}{2N}\right )\Braket{\partial_1\phi^{(1)}(\xi_1)\prod_{i=2}^{n+1}\phi^{(i)}(\xi_i)}_{\mu _{\mrm{eq}}\overset{n}{\bigotimes} L_N}
				\\+\left (\dfrac{1}{\alpha }-\dfrac{1}{2N}\right )\Braket{\partial_1\phi^{(1)}(\xi_1)\prod_{i=2}^{n+1}\phi^{(i)}(\xi_i)}_{\overset{n+1}{\bigotimes} L_N}+\dfrac{1}{2}\Braket{\mcal{D}[\phi^{(1)}](\xi_1,\xi_2)\prod_{i=2}^{n+1}\phi^{(i)}(\xi_{i+1})}_{\overset{n+2}{\bigotimes} L_N}
				\\+\dfrac{1}{N\alpha}\sum_{i=2}^{n+1}\Bigg(\Braket{\phi^{(1)}(\xi_1)\partial_1\phi^{(i)}(\xi_1)\prod_{\substack{j=2\\j\neq i}}^{n+1}\phi^{(j)}(\xi_j)}_{\overset{n}{\bigotimes} L_N}+\Braket{\phi^{(1)}(\xi_1)\partial_1\phi^{(i)}(\xi_1)\prod_{\substack{j=2\\j\neq i}}^{n+1}\phi^{(j)}(\xi_j)}_{\mu _{\mrm{eq}}\overset{n-1}{\bigotimes} L_N}\Bigg)
			\end{multline}
			Using that $ L_N$ has zero mass and defining $\phi_{n+1}(\xi_1,\dots,\xi_{n+1})\defi \prod_{i=1}^{n+1}\phi^{(i)}(\xi_i)$ and using the operators defined earlier leads to:
			\begin{multline}\label{DSbuiltl}
				\Braket{\Xi[\phi_{n+1}]}_{\overset{n+1}{\bigotimes} L_N}=\left (\dfrac{1}{\alpha }-\dfrac{1}{2N}\right )\Braket{\partial_1\phi_{n+1}}_{\mu _{\mrm{eq}}\overset{n}{\bigotimes} L_N}+\left (\dfrac{1}{\alpha }-\dfrac{1}{2N}\right )\Braket{\partial_1\phi_{n+1}}_{\overset{n+1}{\bigotimes} L_N}
				\\+\dfrac{1}{2}\Braket{\mcal{D}[\phi_{n+1}]}_{\overset{n+2}{\bigotimes} L_N}+\dfrac{1}{N\alpha}\Bigg(\Braket{\Theta\Xi[\phi_{n+1}]\ }_{\overset{n}{\bigotimes} L_N}+\Braket{\theta\Xi[\phi_{n+1}] }_{\overset{n-1}{\bigotimes} L_N}\Bigg).
			\end{multline}
			Using that $\overset{n+1}{\bigotimes}\mcal{C}_c^{\infty}(\R)$ is dense in $\mcal{C}_c^{\infty}(\R^{n+1})$ for $\|.\|_{\mcal{C}^{1}(\R^{n+1})}$ and that $\Xi$ and $\mcal{D}$ are continuous for $\|.\|_{\mcal{C}^{1}(K)}$ for any compact $K\subset\R$, we deduce that \eqref{DSbuiltl} is valid for any $\phi_{n+1}\in\mcal{C}_c^{\infty}(\R^{n+1})$.
			Now using that for all $0<\kappa<\kappa'$,
			$$\|\Xi_1\left[\phi_{n+1}\right] \|_{\kappa',1}\leq C\|\phi_{n+1}\|_{\kappa,1}$$ and that $\mcal{C}_c^{\infty}(\R^{n+1})$ is dense in $\mcal{C}^1(\R^{n+1})$ for $\|.\|_{\kappa,1}$, we obtain \eqref{DSlvlnthm} for general function $\phi_{n+1}\in\mcal{A}_{\kappa,1}$.
		\end{Pro}

\section{Large deviation principles}\label{app LDP}
We recall that the definition of the space $\mcal{A}_{\kappa,r}$ was introduced in \eqref{def:norm kappa}. The notation $\mcal{A}_{\kappa,r}^{W}$ denotes the same space when changing the potential $V\rightsquigarrow W$ in \eqref{def:norm kappa}.
\begin{Lem}[Replacement lemma]\label{lem:replacement}
	Let $V$ be a potential satisfying Assumptions \ref{assumptions}, $a>0$ the largest numbers such that $[-1-3a,1+3a]\subset\mcal{U}\cap\R$. Set $V_\phi(x)=\phi(x)V(x)+(1-\phi(x))(x^{2}+C)$ for $\phi\in\mcal{C}_c^{\infty}(\R)$ such that $\ind{[-1-a,1+a]}\leq\phi\leq \ind{[-1-2a,1+2a]}$ and $C>0$ large enough. The potential $V_\phi\in\mcal{C}_{\mrm{loc}}^{\infty}(\R)$ and the corresponding equilibium measures $\mu_{\mrm{eq}}^{V}=\mu_{\mrm{eq}}^{V_\phi}$ and for all $k\geq1$ and $f_1,\dots,f_q\in\mcal{A}_{\kappa,0}^{V_\phi}$ for some $\kappa>0$:
	$$\E_N^{V}\left [\prod_{i=1}^{q}\bm{\mu}_N(f_i)\right ]=\E_N^{V_\phi}\left [\prod_{i=1}^{q}\bm{\mu}_N(f_i)\right ]+\prod_{i=1}^{q}\|f_i\|_{\kappa,0}O\left (e^{-c\alpha}\right )$$
	where the remainder $O\left (e^{-c\alpha}\right )$ depends on $q$.
\end{Lem}

\begin{Pro}
	The claim that $V_\phi\in\mcal{C}_{\mrm{loc}}^{\infty}(\R)$ is clear from Assumption \textit{\ref{assumption1}.} To check that the two equilibrium measures coincide, we show that $\mu_{\mrm{eq}}^{V_\phi}$ satisfies \eqref{eq:effpotentialV}. For all $x\in[-1,1]$, since $V(x)=V_\phi(x)$:
	$$ C_\mrm{eq}^{V}=V(x)-2\int_\R\log|x-y|\diff\mu _{\mrm{eq}}^{V}(x)= V_\phi(x)-2\int_\R\log|x-y|\diff\mu _{\mrm{eq}}^{V}(x).$$
	Thus $x\in[-1,1]\mapsto V_\phi(x)-2\int_\R\log|x-y|\diff\mu _{\mrm{eq}}^{V}(x)$ is constant. Furthermore, for $|x|>1$:
	\begin{align*}
	V_\phi(x)-2\int_\R\log|x-y|\diff\mu _{\mrm{eq}}^{V}(x)-C_\mrm{eq}^{V}&=\phi(x)\left(V(x)-2\int_\R\log|x-y|\diff\mu _{\mrm{eq}}^{V}(x)-C_\mrm{eq}^{V}\right)\\&\quad+(1-\phi(x))\left(x^{2}+C-2\int_\R\log|x-y|\diff\mu _{\mrm{eq}}^{V}(x)-C_\mrm{eq}^{V}\right) 
	\end{align*}
	In the RHS, the first term is positive by \eqref{eq:effpotentialV}. For the second one, since $$x\mapsto x^{2}+C-2\int_\R\log|x-y|\diff\mu _{\mrm{eq}}^{V}(x)-C_\mrm{eq}^{V}$$ is continuous and behaves like $x^{2}+o(x^{2})$ at infinity, it is bounded by below and by taking $C>0$ large enough:
	$$x^{2}+C-2\int_\R\log|x-y|\diff\mu _{\mrm{eq}}^{V}(x)-C_\mrm{eq}^{V}>0.$$
	Thus $\mu_{\mrm{eq}}^{V}$ satisfies the equilibrium condition \eqref{eq:effpotentialV} for $V_\phi$ and thus $\mu_{\mrm{eq}}^{V}=\mu_{\mrm{eq}}^{V_\phi}$ and $C_{\mrm{eq}}^{V}=C_{\mrm{eq}}^{V_\phi}$.
	
	Now, let $W\in\{V,V_\phi\}$ and set $J=[-1-a,1+a]$ and $\mbb{P}_N^{W,J}$ the probability measure  defined by:
	$$\diff\mbb{P}_N^{W,J}(\bm{\lambda})=\dfrac{\mbf{1}_{\forall j, \lambda_j\in J}}{\mcal{Z}_N^{J}[W]}\prod_{i<j}^{N}\left |\lambda_i-\lambda_j\right|^{\beta }.\prod_{i=1}^{N}e^{-\frac{\alpha}{2} W(\lambda_i)}\diff\lambda_i.$$
	It is clear that $\mcal{Z}_N^{J}[W]/\mcal{Z}_N[W]=\mbb{P}_N^{W}\left(\forall j, \lambda_j\in J\right)=1+O(e^{-c\alpha})$ for some $c>0$ by Proposition \ref{appprop: LDP lambdamax}. Thus,
	\begin{align*}\label{key}
		\E_N^{W}\left [\prod_{i=1}^{q}\bm{\mu}_N(f_i)\right ]&=\dfrac{\mcal{Z}_N^{J}[W]}{\mcal{Z}_N[W]}	\E_N^{W,J}\left [\prod_{i=1}^q\bm{\mu}_N(f_i)\right ]+\E_N^{W}\left [\prod_{i=1}^{q}\bm{\mu}_N(f_i)\cdot\mbf{1}_{\exists j, \lambda_j\notin J}\right ]
		\\&=	\E_N^{W,J}\left [\prod_{i=1}^{q}\bm{\mu}_N(f_i)\right ]+O\left (e^{-c\alpha}\prod_{i=1}^{q}\|f_i\|_{\kappa,0}\right )
		\\&\quad+\E_N^{W}\left [\prod_{i=1}^{q}\bm{\mu}_N(e^{\kappa W})\cdot\mbf{1}_{\exists j, \lambda_j\notin J}\right ]O\left(\prod_{i=1}^{q}\|f_i\|_{\kappa,0}\right) .
	\end{align*}
	where the norm $\|\cdot\|_{\kappa,0}$ has to be understood with respect to the potential $W$.
	Since: $$\E_N^{W}\left [\prod_{i=1}^{q}\bm{\mu}_N(e^{\kappa W})\cdot\mbf{1}_{\exists j, \lambda_j\notin J}\right ]\leq C_qe^{-c\alpha}$$ 
	for some $C_q>0$ and $c>0$ (this bound can be shown using the same bounds in the proof of Lemma \ref{lemma:trunc} for $A_i$).
	Since, $\E_N^{V_\phi,J}\left [\prod_{i=1}^{q}\bm{\mu}_N(f_i)\right ]=\E_N^{V,J}\left [\prod_{i=1}^{q}\bm{\mu}_N(f_i)\right ]$, it concludes the proof.
\end{Pro}

The following Proposition is the LDP for the largest particle $\lambda_{\mrm{max}}$. A corollary of the next proposition is that with probability greater than $1-e^{-\alpha c}$ for some $c>0$, there is no particle outside of $[-1-\varepsilon,1+\varepsilon]$.

\begin{Prop}[LDP for $\lambda_{\max}$]\label{appprop: LDP lambdamax}
	Let $V:\R\rightarrow\R$ be a continuous potential satisfying:
	$$\liminf_{x\rightarrow +\infty}\dfrac{V(x)}{2\log|x|}>1.$$
	The laws of the largest particle $\lambda_{\max}=\max\limits_{i\in\llbracket1,N\rrbracket}\lambda_i$ under $\P_N$ satisfies a large deviaion principle at speed $\alpha$ with good rate function $\mcal{I}_V$ defined by:
	$$\mcal{I}_V(x)=\begin{cases}
		\dfrac{1}{2}V_\mrm{eff}(x)\hspace{1cm}&\text{if}\quad x\geq1,
		\\+\infty&\text{if}\quad x<1,
	\end{cases}$$
	where the effective potential $V_{\mrm{eff}}$ has been defined in \eqref{eq:effpotentialV}.
\end{Prop}

\begin{Pro} We follow \cite[Appendix A]{BoG1}, the proof is almost the same so we only stress the differences and follow the same structure of the proof. From \cite[\textbf{A.1}]{BoG1}, $\mcal{I}_V$ is a good rate function for the same reasons. To get the exponential tightness of the laws of $\lambda_{\max}$, following \cite[\textbf{A.2}]{BoG1}, we obtain similarly:
	\begin{equation*}
		\dfrac{\mc{Z}_N^{\beta_N}[V]}{\mcal{Z}_{N-1}^{\beta_{N}}[V]}
		\geq\int_\R\diff\lambda_Ne^{-\alpha \frac{V(\lambda_N)}{2}}\mbb{E}_{N-1}^{\beta_N}\left[\exp\left(\beta_N\sum_{i=1}^{N-1}\log|\lambda_N-\lambda_i|-(\alpha-(N-1)\beta_N)\sum_{i=1}^{N-1}\dfrac{V(\lambda_i)}{2}\right) \right].
	\end{equation*}
	Introducing $\diff\chi(x)\defi \dfrac{e^{-\frac{V(x)}{2}}}{\kappa}\diff x$ and $\kappa\defi\int_\R e^{-\frac{V(x)}{2}}\diff x$, we get by Jensen's inequality:
		\begin{align*}
		\dfrac{\mc{Z}_N^{\beta_N}[V]}{\mcal{Z}_{N-1}^{\beta_{N}}[V]}
		&\geq\kappa\exp\Bigg(\chi\otimes\mbb{P}_{N-1}^{\beta_N}\left(-(\alpha-1)\frac{V(\lambda_N)}{2}+(\alpha-\beta_N)\int_\R\log|\lambda_N-x|\diff\mu_{N-1}(x)\right)
		\\&\quad-\beta_N\sum_{i=1}^{N-1}\frac{V(\lambda_i)}{2}\Bigg)
		\\&\geq\kappa\exp\left(\alpha c+(\alpha-\beta_N)\mbb{E}_{N-1}^{\beta_N}\left[\int_\R\left(\int_{\R}\log|\lambda_N-x|\diff\chi(\lambda_N)-\frac{V(x)}{2}\right)\diff\mu_{N-1}(x)\right]\right),
	\end{align*}
	for some $c>0$. Now since $x\mapsto\int_{\R}\log|\lambda_N-x|\diff\chi(\lambda_N)$ is bounded by below (as a continuous function growing like $\log|x|$ at infinity) and by exponential tightness of \cite[Eq (2.6.21)]{anderson2010introduction}, we get:
		\begin{equation*}
		\dfrac{1}{\alpha}\log\dfrac{\mc{Z}_N^{\beta_N}[V]}{\mcal{Z}_{N-1}^{\beta_{N}}[V]}>-\infty.
	\end{equation*}
	For the \textbf{A.3} part, we first recall that the LDP for the empirical measure holds at speed $N\alpha$ \cite{Garcia}. Now, we set:
	$$Y_N\defi\dfrac{\mcal{Z}_{N-1}^{\beta_N}[\frac{N}{N-1}V]}{\mcal{Z}_N^{\beta_N}[V]},\hspace{0.8cm}\Xi_N(\xi)\defi\E_{N-1}^{\beta_N,\frac{N }{N-1}V}\left[\exp\left(\beta_N\sum_{i=1}^{N-1}\log|\xi-\lambda_i|-(\alpha-1)\frac{V(\xi)}{2}\right)\ind{\forall i, \lambda_i<\xi} \right].$$
	We thus have for every closed set $F\subset [1,+\infty)$,
	$$\P_N\left(\lambda_{\max}\in F\right)=Y_N\int_Fe^{-\frac{V(\xi)}{2}}\Xi_N(\xi)\diff\xi.$$ 
	The upper bound on $\Xi_N$:
	$$\sup_{\xi\in F}\Xi_N(\xi)\leq\exp\left(\alpha\Big(\eta'-\inf_{\xi\in F}(\mcal{I}_V(\xi)+C_{\mrm{eq}})\Big)\right) $$
	 valid for all $\eta'>0$ and the bound:
	$$\dfrac{1}{Y_N}\geq \exp\left(-\alpha\Big(\eta''+\inf_{\xi\in\R}(\mcal{I}_V(\xi)+C_{\mrm{eq}})\Big)\right) $$ that holds for any $\eta''>0$ are obtained in the exact same way. Noticing that because of \eqref{eq:effpotentialV}, $\inf_{\xi\in\R}\mcal{I}_V(\xi)=0$, this leads to the bound:
	$$\limsup_{N \rightarrow +\infty}\dfrac{1}{\alpha}\log\P_N\left(\lambda_{\max}\in F\right)\leq-\inf_{\xi\in F}\mcal{I}_V(\xi)$$
	For the \textbf{A.4} part, the bounds are obtained in the same way which leads to the lower bound valid for all $x>1$:
		$$\liminf_{\varepsilon\rightarrow0}\liminf_{N \rightarrow +\infty}\dfrac{1}{\alpha}\log\P_N\left(\lambda_{\max}\in (x-\varepsilon,x+\varepsilon)\right)\geq-\mcal{I}_V(x).$$
		\qedsymbol{symbol}
	\end{Pro}

As discussed in Subsection \ref{subsec:connection}, the thermal equilibrium measure $\mu_{\alpha}$ converges to the equilibrium measuer $\mu_{\mrm{eq}}$ as $\alpha\rightarrow\infty$.

\begin{Lem}\label{lem: stieltjes transform approx}
	The following convergence in law holds:
	$$\mu_\alpha\tend{\alpha\rightarrow\infty}\mu_{\mrm{eq}}.$$
\end{Lem}

\begin{Pro}

	The proof is similar to ideas used in \cite[Lemma 3.6]{padilla2024emergence} \cite[Lemma 2.1]{armstrong2022thermal}. First, notice that $\mcal{E}$ defined in \eqref{eq:energy functional}	is lower semicontinuous under our assumptions by \cite[Lemma 2.4]{serfaty26}. The lack of positivity of the entropy prevents us from concluding immediately; we circumvent this issue by a trick used in \cite{armstrong2022thermal}; namely, we can rewrite the function $\mathcal{E}_\alpha$ as
	\begin{align*}
		\mathcal{E}_\alpha(\mu)&=\left(1-\frac{\alpha_\ast}{\alpha}\right)\int_{\R} V(x)\, \mathrm d \mu(x)-\iint_{\R^{2}} \log|x-y|~\diff\mu(x)\diff\mu(y)
		\\&\quad+\left(\int_\R \frac{\alpha_\ast}{\alpha}V(x)\, \mathrm d \mu(x) +\frac{2}{\alpha}\int_\R \log \dfrac{\diff\mu}{\diff x}(x)\diff\mu(x)\right),
	\end{align*}
	where $\alpha_\ast$ is such that $\int_{|x|\geq 1} e^{-\frac{\alpha_\ast}{2} V(x)}~\diff x$ is finite (we have a good deal of freedom over what to choose for $\alpha_\ast$ since we assume that $V=x^2+C$ outside of a compact set $K$). Now, let $(\mu_N)$ be a sequence of probability measures converging weakly to some probability measure $\mu$. By \cite[Lemma 2.4]{serfaty26}, we have for any $\epsilon>0$ that 
	\begin{multline*}
		\liminf_{N \rightarrow \infty}\left(1-\frac{\alpha_\ast}{\alpha_N}\right)\int_{\R} V\, \mathrm d \mu_N-\iint_{\R^{2}} \log|x-y|~\diff\mu_N(x)\diff\mu_N(y)  \\\geq \left(1-\epsilon\right)\int_{\R} V\, \mathrm d \mu-\iint_{\R^{2}} \log|x-y|~\diff\mu(x)\diff\mu(y).
	\end{multline*}
	Next, 
		\begin{equation*}
		\liminf_{N \rightarrow \infty}\left(\int_{\R} \frac{\alpha_\ast}{\alpha_N}V(x)\, \mathrm d \mu_N(x) +\frac{2}{\alpha_N}\int_\R  \log \diff\mu_N(x)\diff\mu_N(x)\right) 
		\geq \liminf_{N \rightarrow \infty}\frac{-2}{\alpha_{N}}\int_{\R}e^{-\frac{\alpha_\ast}{2}V(x)-1}~\diff x=0
	\end{equation*}
	where we have used the inequality
	\begin{equation*}
		\frac{c\gamma}{\alpha_N}+\frac{\gamma \log \gamma}{\alpha_N} \geq -\frac{1}{\alpha_N}e^{-c-1},
	\end{equation*} 
	with $c=\alpha_\ast V$ and $\gamma=\mu_{N}$,
	which can be seen by minimizing the above as a function of $\gamma$. 
	
	It follows that
	\begin{equation*}
		\liminf_{N \rightarrow \infty}\mathcal{E}_{\alpha_N}(\mu_N)\geq  \mathcal{E}(\mu)+\int_{\R} (1-\epsilon)V(x)~ \diff\mu(x)
	\end{equation*}
	and $\epsilon>0$ was arbitrary so we are done.
	
	Now, let $\{\mu_{\alpha_N} \}$ be a sequence of minimizers of $\mathcal{E}_{\alpha_N}$; it is sufficient for us to show that this sequence is tight. Indeed, if it is then by Prokhorov's theorem it has some weak limit $\mu$ and we can conclude as in \cite[Theorem 3.8]{serfaty26}. Namely, we have by lower semicontinuity of the functionals that 
		\begin{equation}
		\mathcal{E}(\mu_{\mrm{eq}}) = \limsup_{N \rightarrow +\infty}\mathcal{E}_{\alpha_N}(\mu_{\mrm{eq}}) \geq \limsup_{N \rightarrow +\infty}\mathcal{E}_{\alpha_N}(\mu_{\alpha_N}) \geq \mathcal{E}(\mu). 
	\end{equation}
	Uniqueness of the equilibrium measure implies $\mu=\mu_{\mrm{eq}}$ and hence that $\mu_{\alpha_N}$ converges weakly to $\mu_V$ under the topology of weak convergence.

	To show that $\{\mu_{\alpha_N}\}$ is a tight sequence of probability measures, we first show that the sequence $\mathcal{E}_{\alpha_N}(\mu_{\alpha_N})$ is bounded. Without loss of generality, assume $\alpha_1=\inf \alpha_N$. The boundedness follows from $\frac{1}{\alpha_N}\rightarrow 0$ and the minimality of $\mu_{\alpha_N}$:
	\begin{equation}
		\mathcal{E}_{\alpha_N}(\mu_{\alpha_N}) \leq \mathcal{E}_{\alpha_N}(\mu_{\alpha_1}) \leq \mathcal{E}_{\alpha_1}(\mu_{\alpha_1}):=C_1<+\infty
	\end{equation}
	for all $N$.  Now, for any $C_2>0$, there is a compact set $K \times K$ outside of which 
	\begin{equation}
		-\log|x-y|+\gamma_N\frac{V(x)}{2}+\gamma_N\frac{V(y)}{2}>C_2
	\end{equation}
	with $\gamma_N=\left(1-\frac{\alpha_\ast}{\alpha_N}\right)$ by our coercivity assumptions on $V$, where the set $K \times K$ is independent of $N$. Moving some of the potential onto the entropy as above, we have
	\begin{align*}
		C_1 &\geq \mathcal{E}(\mu_{\alpha_N})+\left(1-\frac{\alpha_\ast}{\alpha_N}\right)\int_{\R} V(x)\, \mathrm d \mu_{\alpha_N}(x)+\int_{\R} \left(\frac{\alpha_\ast}{\alpha_N}V+\frac{2}{\alpha_N}\log \mu_{\alpha_N}\right)\, \mathrm d \mu_{\alpha_N} \\
		&\geq \iint_{\R^2} \left(-\log|x-y|+\gamma_n\frac{V(x)}{2}+\gamma_N\frac{V(y)}{2}\right)\, \mathrm d \mu_{\alpha_N}(x)\, \mathrm d \mu_{\alpha_N}(y)-\frac{1}{\alpha_N}\int_{\R} e^{-\frac{\alpha_\ast}{2}V(x)-1}\, \mathrm \diff x \\
		& \geq -C_3+C_2\mu_{\alpha_N}\otimes \mu_{\alpha_N}((K \times K)^c) \\
		& \geq -C_3+C_2\mu_{\alpha_N}(K^c),
	\end{align*}
	where we have used that $(x,y)\in K\times K\mapsto-\log|x-y|+\gamma_N\frac{V(x)}{2}+\gamma_N\frac{V(y)}{2}$ is everywhere bounded below  independently of $N$, since $V$ is bounded below, $K$ is bounded and
	$$\log|x-y|\leq\log(1+|x|)+\log(1+|y|).$$
	We have also used that $-\frac{1}{\alpha_N}\int_{\R} e^{-\frac{\alpha_\ast}{2}V(x)-1}\, \mathrm \diff x$ is bounded below as well since $\alpha_N \rightarrow \infty$. Since $C_2$ can be made arbitrarily large, this gives us the tightness of the sequence $\{\mu_{\alpha_N}\}$, and we conclude the desired convergence.
\end{Pro}

\section{The Gaussian case}\label{app: Gaussian}
	\subsection{Convergence of the series}\label{subsec2:convergence}

The goal of this subsection is to discuss the definition a linear map $\nu_\alpha:\psi\mapsto\sum_{j=0}^{+\infty}\tfrac{\nu_j(\psi)}{\alpha^j}$ acting on the space $\mcal{C}^{\infty}([-1,1])\defi\cap_{j\geq0}\,\mcal{C}^{j}([-1,1])$ endowed with the norms:
$$\|\psi\|_{K,\,\mcal{C}^{\infty}}\defi\sum_{j\geq0}\dfrac{\|\psi\|_{\mcal{C}^{j}([-1,1])}}{K^{\sqrt{j}}},\quad\quad\quad K>0.$$

\begin{Prop}[Convergence of the series]\label{appthm:conv muk}Take $V(x)=2x^{2}$, then for all  $\psi\in\mcal{C}^{\infty}([-1,1])$, for all $i\geq0$, there exists $C,r>0$ (both independent of $i$) such that
	$$|\nu_i(\psi)|\leq C r^{i}\cdot \|\psi\|_{\mcal{C}^{2i^{2}}([-1,1])}.$$ 
	Thus for all $\psi\in\mcal{C}^{\infty}([-1,1])$ such that $\|\psi\|_{K,\,\mcal{C}^{\infty}}<\infty$ for $K>0$, then for all $\alpha>K^{\sqrt{2}}r$ the serie $\sum_{j\geq0}\dfrac{\nu_j(\psi)}{\alpha^{j}}$ converges absolutely (or equivalently the sequence $\mu_k(\psi)$ converges). Furthermore, defining $\mu_{\infty}(\psi)\defi\sum_{j\geq0}\dfrac{\nu_j(\psi)}{\alpha^{j}}$, we have for all $\alpha>K^{\sqrt{2}}r$:
	$$\mu_{\infty}(\psi)\leq\dfrac{C}{1-\dfrac{K^{\sqrt{2}}r}{\alpha}}\cdot\|\psi\|_{K,\,\mcal{C}^{\infty}}$$
\end{Prop}

\begin{Pro}
	By \eqref{eq:inverse Xi formula} and \cite[Lemma A.1 \textit{(i)}]{dworaczekguera2025clt},
	$$\|\Xi^{-1}[\psi]\|_{\mcal{C}^{i}([-1,1])}\leq \dfrac{1}{\pi}\|\psi\|_{\mcal{C}^{i+1}([-1,1])}$$ for all $i\geq0$.
	We now prove the estimate on $\nu_{i}$ by induction. For $n=1$, we have: $$|\nu_1(\psi)|\leq\|\Xi^{-1}[\psi]\|_{\mcal{C}^{1}([-1,1])}\leq\dfrac{1}{\pi}\|\psi\|_{\mcal{C}^{2}([-1,1])}.$$
	We first obtain that $\nu_{i}$ is continuous with respect to the $\mcal{C}^{2i^{2}}$-norm. For that, we recall that if for all $i\geq1$, $\nu_i$ continuous for the $\mcal{C}^{k_i}$-norm with constant $C_i>0$ then, for all $i,j\geq1$
	\begin{align*}
		|\nu_i\otimes\nu_j\left(\mcal{D} \Xi^{-1}[\psi]\right)|\leq C_iC_j\|\mcal{D} \Xi^{-1}[\psi]\|_{\mcal{C}^{k_i+k_j}([-1,1]^{2})}&\leq C_iC_j\|\Xi^{-1}[\psi]\|_{\mcal{C}^{k_i+k_j+1}([-1,1])}
		\\&\leq \dfrac{C_iC_j}{\pi}\|\psi\|_{\mcal{C}^{k_i+k_j+2}([-1,1])}.
	\end{align*}
	We thus obtain, using \eqref{eq:constrait nu i}, that $k_n$ must satisfy $k_n\geq \max\left(k_{n-1}+2, \max\limits_{1\leq p\leq n-1} k_p+k_{n-p}+2\right)$ and $k_1\geq2$, therefore it is easy to see that choosing $k_n=2n^{2}$ works. We now compute the constant. By the same reasoning, using \eqref{eq:constrait nu i}, it can be shown that the constants $C_n$ must satisfy $C_1\geq1$ and these conditions:
	$$C_n\geq \dfrac{C_{n-1}}{\pi}+\dfrac{1}{2\pi}\sum_{p=1}^{n-1}C_pC_{n-p}.$$
	Looking for a sequence $(C_n)$ such that these bounds are equalities and defining $F(z)=\sum_{i\geq1}C_iz^{i}$, we obtain that:
	$$F(z)-z=\dfrac{z}{\pi}F(z)+\dfrac{1}{2\pi}F(z)^{2}\iff F(z)^{2}+\left(z-\pi\right)F(z) +2\pi z=0.$$
	Thus, $F(z)=\pi-z-\sqrt{(z-\pi)^{2}-2\pi z}$ where we used the fact that $F(0)=0$. The radius of convergence of $F$ is given by the singularity of the square root, thus one has to solve: $$(z-\pi)^{2}-2\pi z=z^{2}-4\pi z+\pi^{2}=0$$ whose solutions are $2\pi\pm\pi\sqrt{3}$. The closest solution to 0 is $z_0=2\pi-\pi\sqrt{3}\in(0,1)$ and thus since $F(z)\sim c+ c'\sqrt{1-\dfrac{z}{z_0}}$, we obtain by standard singularity analysis that such sequence $(C_n)$ exists and admits the following asymptotc bound $C_n= O\left(\dfrac{1}{z_0^{n}n^{3/2}}\right)$. This concludes the proof.
\end{Pro}

\begin{Rem}
	The proof of Proposition \ref{appthm:conv muk} only works because in the quadratic case, the function $S$ defined in \eqref{eq: S} is constant.
\end{Rem}

\begin{Prop}[Constraint on $\mu_\infty$]
	For any $f\in\mcal{C}^{\infty}([-1,1])$ such that $\|f\|_{K,\mcal{C}^{\infty}}<\infty$ for some $K>0$, then for all $\alpha$ large enough and $V(x)=2x^{2}$:
	$$\mu_\infty(V'f)-\mu_\infty\otimes\mu_\infty\left(\mcal{D}[f]\right)-\dfrac{2}{\alpha}\mu_\infty(f')=0.$$
\end{Prop}

\begin{Pro}
	Using and summing the constraints \eqref{eq:constrait nu i}, the definiton of the master operator Definition \ref{def:ope D} and using the constraint satisfied by $\mu_{\mrm{eq}}$ \eqref{eq:contrainst mu infinity}, we obtain the conlusion noticing that $$\|f\|_{K,\,\mcal{C}^{\infty}}<\infty\Rightarrow\|f'\|_{K,\,\mcal{C}^{\infty}}<\infty.$$\qedsymbol{symbol}
\end{Pro}

\begin{Rem}
	This constraint has the same form as the one satisfied by the thermal equilibrium measure $\mu_\alpha$ except that $\mu_{\infty}$ is supported on $[-1,1]$. This latter fact confirms that $\mu_{\infty}\neq\mu_\alpha$. If we were able to take $f(x)=(x-z)^{-1}$ for $z\in\C_+$, we would obtain the equation \cite[Eq. (12)]{AllezBouchaudGuionnet} whose unique solution is the Stieltjes transform of the termal equilbrium measure among functions which behave like $-z^{-1}$ at infinity. We would conclude that $\mu_{\infty}=\mu_{\alpha}$ which is wrong.
\end{Rem}

\subsection{Numerical simulations}
\begin{figure}[h]
	\centering
	\includegraphics[width=14cm]{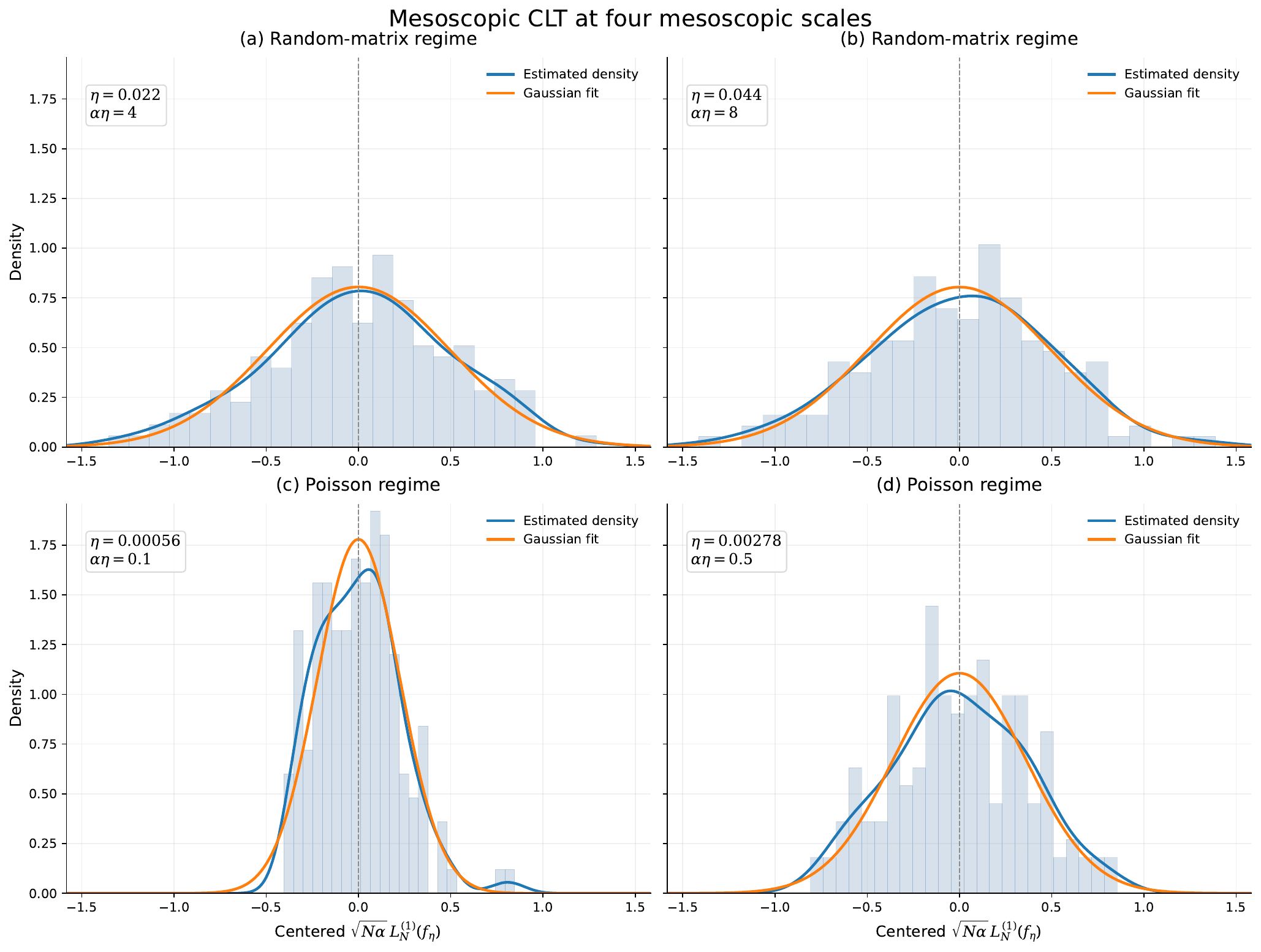}
	\caption{What we can notice on this plot is that the estimated density of $X=\sqrt{N\alpha}L_N^{(1)}(f_\eta)$ for $E=0$ and $f(x)=(1+x^{2})^{-1}$ is asymptotically Gaussian as predicted by Theorem \ref{thm:mesoc clt bulk}. Furthermore in the Poisson regime, we can see that the variance of $X$ is $\eta$-dependent and goes to zero (indeed the correct rescaling is $\sqrt{N\eta^{-1}}$ and not $\sqrt{N\alpha}$). $N=3000$, $\beta=N^{-0.35}=0.06$ and $\alpha=N^{0.65}=182$.}\label{fig:plot}
\end{figure}
Using the tridiagonal model of \cite{dumitriu2002matrix}, we present simulations of the Gaussian $\beta$-ensembles in the intermediary temperature regime (see Figure \ref{fig:plot} on the following page).

\section{Stein's method}\label{app:Stein}
The following result applies to any Gibbs measure $\Q_N = Z_N^{-1} e^{-\beta H}$ on $\R^N$ with a generator $$\mc{L} = -\Delta_{\R^N} + \beta \nabla_{\R^N} H \cdot \nabla_{\R^N} $$ such that $\mcal{C}^{2}(\R^{N}) \subset {\operatorname{Dom}}(\mc L) $ -- meaning that $\mc{L} [F] \in L^2(\P_N)$ for $F \in \mcal{C}^{2}(\R^{N})$. We work in the one-dimensional setting for simplicity, but there is also a multi-variate version of this result. 
Let $\mathbf{W}_q$ be the Wasserstein $q$-distance for probability measures (random variables) on~$\R$.

\begin{Prop}\label{prop:bound Wp}
	Let $q\geq1$, and $G \sim \mcal{N}(0,\sigma^{2})$ for $\sigma>0$. 
	Let $X \in \mcal{C}^{1}(\R^{N})$ and suppose that $X= \lambda^{-1}\big( \mc{L} [F] +Z \big)$ where $\lambda>0$,  $F\in\mcal{C}^{2}(\R^{N})$ and $Z \in L^q(\Q_N)$.
	Then:
	\begin{equation*}\label{Steinbound-gen}
		\mathbf{W}_q(X,G)\leq \lambda^{-1}\big( \sqrt{q/\sigma^{2}}\| W\|_{q}+\|Z\|_{q} \big), \qquad  W = \nabla X \cdot \nabla F - \lambda \sigma^{2} . 
	\end{equation*}

\end{Prop}

\begin{Prop}\label{prop:bound_normq}
	Under the same conditions as in Proposition~\ref{prop:bound Wp}, one has for  all $ q\geq2$,
	\[
	\|X\|_q \le \lambda^{-\frac12} \sqrt{q \|\Gamma \|_{q/2}} + \lambda^{-1} \|Z\|_q , \qquad \Gamma := \nabla F \cdot\nabla X . 
	\]
\end{Prop}

\begin{Pro}
	Integrating by parts and using Holder's inequality
	\[\begin{aligned}
		\|X\|_q^q  &= \lambda^{-1}\E_N[(\mc{L} [F] +Z)X^{q-1}] \\
		& = \lambda^{-1}  \E_N[ ((q-1) \Gamma + ZX)  X^{q-2} ] \\
		&\le \lambda^{-1}  \|(q-1)\Gamma + ZX\|_{q/2} \|X\|_q^{q-2}
	\end{aligned}\]
	Dividing both sides by $ \|X\|_q^{q-2}$, using that $\|(q-1)\Gamma + ZX\|_{q/2}  \le q\|\Gamma \|_{q/2} + \|Z\|_q\|X\|_q$ and rearranging, we obtain 
	\[
	\|X\|_q^2 \le  \lambda^{-1} q \|\Gamma \|_{q/2} +  \lambda^{-1} \|Z\|_q\|X\|_q
	\]
	Using that if $x^2 \le \alpha +\gamma x $ for $x>0$, then   \(x\leq \sqrt{\alpha}+\gamma\) (by solving the quadratic equation), we conclude that 
	\[
	\|X\|_q \le \lambda^{-\frac12} \sqrt{q \|\Gamma \|_{q/2}} + \lambda^{-1} \|Z\|_q.
	\]
\end{Pro}

	\bibliographystyle{alpha}
	\bibliography{Meso_CLT}
\end{document}